\pdfoutput=1                 

\documentclass[10pt,leqno]{ijnam}

\def\currentvolume{1}
\def\currentissue{1}
\def\currentyear{2016}
\def\month{December}
\pagespan{1}{18}
\copyrightinfo{2016}{} 

\usepackage[T1]{fontenc}
\usepackage[utf8]{inputenc}    
\usepackage{amsmath,amssymb}
\usepackage{upgreek}

\usepackage{graphicx}
\usepackage{longtable}

\usepackage{subcaption}

\usepackage{url}               
\usepackage[hidelinks]{hyperref}

\newcommand{\iwp}{in the water phase }
\def\jumpl{\left[\kern-0.15em\left[}
\def\jumpr{\right]\kern-0.15em\right]}
\def\Jumpl{\left[\kern-0.25em\left[}
\def\Jumpr{\right]\kern-0.25em\right]}

\begin{document}


\title[A Flow and Transport Model for MEOR at laboratory conditions]
{A Flow and Transport Model for\\
 Simulation of Microbial Enhanced\\
  Oil Recovery Processes at Core Scale\\
   and Laboratory Conditions}


\author[D\'iaz \and Ortiz \and Hern\'andez \and Castorena\and Rold\'an\and Olgu\'in]{Mart\'in A. D\'iaz-Viera \and Arturo Ortiz-Tapia\and Joaqu\'in R. Hern\'andez-P\'erez\and Gladys Castorena-Cort\'es\and Teresa Rold\'an-Carrillo\and Patricia Olgu\'in-Lora}
\address{
  Instituto Mexicano del Petr\'{o}leo
  Eje Central Norte L\'{a}zaro C\'{a}rdenas 152
  07730, CDMX
  M\'{e}xico
}
\email{mdiazv@imp.mx \and aortiztapia2013@gmail.com}

\dedicatory{}


\date{January 1, 2004 and, in revised form, March 22, 2004.}

\thanks{This research was supported by IMP-D.00417 Project ``Recuperaci\'on Mejorada de Hidrocarburos V\'ia Microbiana'' of the Instituto Mexicano del Petr\'oleo.}

\subjclass[2000]{91.60.Np, 91.60.Tn, 91.65.My, 76D45, 76D55, 49Q12, 93C20, 86-04}

\abstract{A general 3D  flow and transport model in porous media was derived applying an axiomatic continuum modeling approach, which was implemented using the finite element method to numerically simulate, analyze and interpret microbial enhanced oil recovery (MEOR) processes under laboratory conditions at core scale.  From the methodological point of view the development stages (conceptual, mathematical, numerical and computational) of the model are shown. This model can be used as a research tool to investigate the effect on the flow behavior, and consequently the impact on the oil recovery, due to clogging/declogging phenomena by biomass production, and interfacial tension changes because of biosurfactant production. The model was validated and then applied to a case study. 
The experimental results were accurately predicted by the simulations. Due to its generality, the model can be easily extended and applied to other cases.}

\keywords{axiomatic continuum modeling approach, microbial enhanced oil recovery, clogging/declogging, interfacial tension, wettability change, trapping number.}

\maketitle

\section{Introduction}
\label{intro}
 The oil fields at the initial stage of operation produce using basically its natural energy, which is known as \textit{primary recovery}. As the reservoir loses energy it requires the injection of gas or water in order to restore or maintain the pressure of the reservoir. This stage is called \textit{secondary recovery}. When the secondary recovery methods become ineffective it is necessary to apply other more sophisticated methods such as steam injection, chemicals, microorganisms, etc. These are known as \textit{tertiary or enhanced oil recovery} (EOR). Some important oil fields in Mexico are entering the third stage.
 
 For the optimal design of enhanced oil recovery methods it is required to perform a variety of laboratory tests under controlled conditions to understand what are the fundamental recovery mechanisms for a given EOR method in a specific reservoir.The laboratory tests commonly have a number of drawbacks, which include among others, that they are very sophisticated, expensive and largely unrepresentative of the whole range of phenomena involved. A proper modeling of the laboratory tests would be decisive in the interpretation and understanding  of recovery mechanisms and in obtaining the relevant parameters for the subsequent implementation of enhanced recovery processes at the well and the reservoir scale.

In this work, a very general 3D  flow and transport model in porous media was obtained to numerically simulate, analyze and interpret microbial enhanced oil recovery (MEOR) processes under laboratory conditions at core scale, such as clogging/declogging and interfacial tension changes because of biosurfactant production. The model was validated with experimental data from \cite{Hendry1997} and then applied to a case study using a sandstone Berea core, while the oil and the microbial culture are from Agua Fr\'ia field, reported in the work of \cite{Castorena2012a}.  
 
 From the methodological point of view the development stages (conceptual, mathematical, numerical and computational) of the model are shown. In particular, for mathematical modeling was used the axiomatic continuum modeling approach, for numerical modeling a finite element method and COMSOL Multiphysics\circledR $\;$ Software for computational implementation.  
 
 In section \ref{StateOftheArt}, we present a review of the State of the Art leading to this work. In section \ref{modeling_methodology}, we present the core idea of this paper, which is to describe the systematic methodology used to define a conceptual model, then how to derive a mathematical model, then how to discretize this later numerically, and finally how to make computational implementations for validation and application to a case study. In section  \ref{model_derivation} a broad description of the axiomatic modeling of continuum systems is given, with the necessary details for the description of flow and transport of the MEOR model. In section \ref{SeccValidationFlowModel} the validation of the flow model is described. In section \ref{validation}, this work was compared with the experimental data and numerical simulations of other authors, to validate the clogging and declogging process. In section \ref{caseDP1}, this MEOR model was applied to our own experimental data. In section \ref{discussion}, the results of the case study are analyzed, and in section \ref{conclusion} the main contributions are summarized.

\section{Review of the State of the Art}\label{StateOftheArt}
 The modeling of the microorganisms behavior influencing the enhanced oil recovery through microorganisms (MEOR) and their activities in the reservoir, has attracted a strong interest from the beginning of the research about MEOR. Some models describe transport equations in one or at most two dimensions, describing the clogging process, without solving the flow equation, as in  \cite{Corapcioglu1984,Corapcioglu1985} and  \cite{Knapp1988}. Others, beside the transport equations and clogging process, do include a flow formulation, coupled with the transport equations, as in  \cite{Islam1990,Islam1993} and  \cite{Sarkar1990}. Still others investigated the mechanisms of the change in wettability, besides the described process, as in  \cite{Chang1995, Chang1993,Chang1992,Chang1991}.\\
 
  One of the problems for modeling, is obtaining the parameters for the equations, so experiments in laboratory conditions have been performed, as in \cite{Bryant1992}. The rate of growth or chemical reactions coming from, or acting on the microbial activity has also been investigated, as in  \cite{Zhang1995}, together with their action on recovery process, due, for example, because of the surfactant produced by the microorganisms. Others, concentrated on some or all parts of the growth and decaying process of the microorganisms, as in  \cite{Desai1999}.\\
  
  The problem of clogging and declogging has been addressed through transport equations as in  \cite{Thullner2004} or  \cite{Kim2006}. Others, give more importance to the effects of surfactant produced by the microorganisms, as in  \cite{Nielsen2010}. Of course there are some models which include both clogging/declogging and surfactant, as in \cite{Li2011}, although most of them are one-dimensional, with exceptions like   \cite{Huang2012}, and this work, which are fully 3-D.


\section{Modeling Methodology}\label{modeling_methodology}
\begin{figure}[htb]
    \centering
    \includegraphics[width=0.8\textwidth]{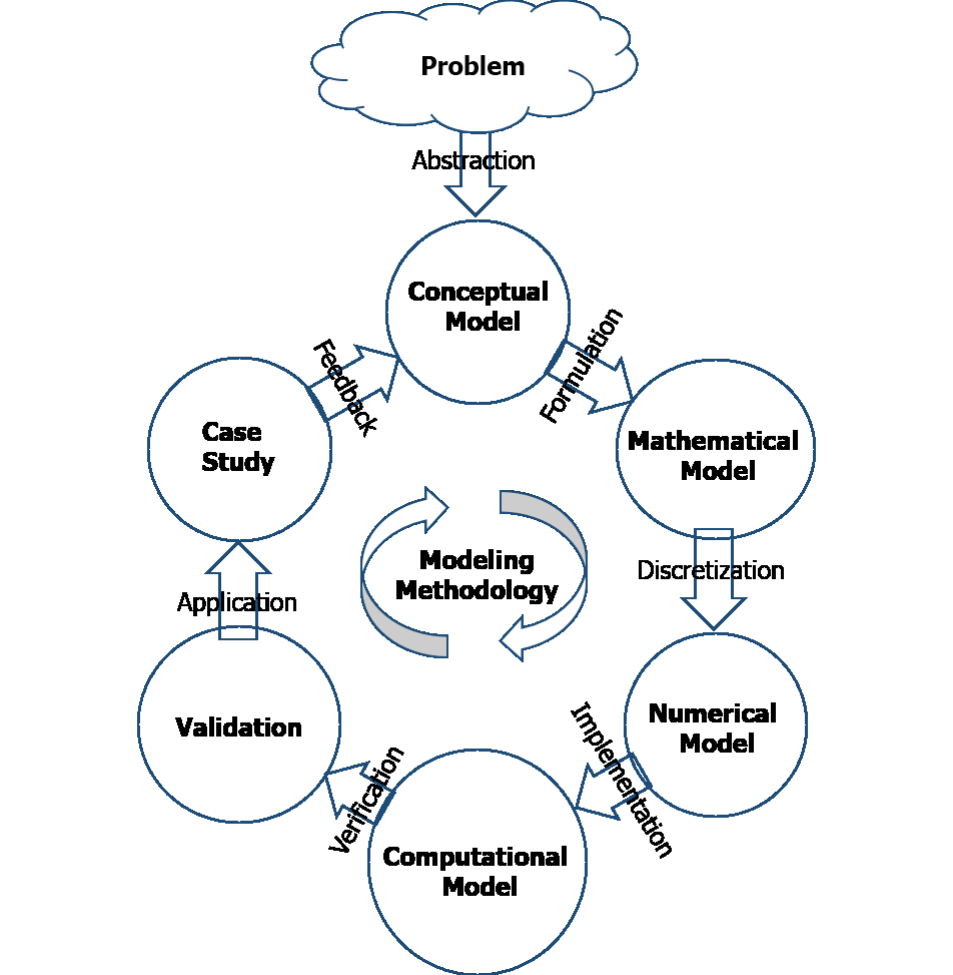}
  \caption[Flow diagram of the Modeling Methodology.]{Flow diagram of the Modeling Methodology.}
\label{Figmethodology}           
\end{figure}

A fundamental issue that is pursued by this paper is to illustrate how to develop  a flow and transport models in porous media applying a general systematic modeling methodology, see Fig. \ref{Figmethodology}. In the literature, to the extent of our knowledge, there are a few publications where this approach is roughly described \cite{helmig1997multiphase,Dietrich2005}.\\

This methodology could be universally applied to develop models in scientific, technological and industrial areas and basically consists on four modeling procedures, namely conceptual, mathematical, numerical and computational models.\\

The conceptual model is an abstraction from a real life problem which comprises the most relevant hypothesis, postulations, assumptions, conditions, restrictions, scope, etc, to be satisfied by the mathematical model. The mathematical model is  a mathematical formulation of the conceptual model, that could be expressed in terms of equations \cite{Herrera2012}. While the numerical model is a discretization version of the mathematical model by the application of the appropriate numerical methods. And finally, the computational model is the implementation of the numerical model in a specific computing platform.\\
 
Two additional stages that are part of the model natural life cycle are its validation and application. The validation step is imperative for the model to be consistent with the conceptual model and normally is a comparison with the output of simplified or referenced problems. Whereas the final step is the application of the model to different  case studies of interest, which are the motivation of model development.\\

The results obtained in the case studies may serve as feedback to the conceptual model and eventually modify its requirements, and consequently another cycle begins in the modeling workflow. Understandably, the change in any of the stages of the above procedure may lead to perform a new loop of model development (Fig. \ref{Figmethodology}).

\clearpage
\section{Model derivation}\label{model_derivation}

\subsection{Conceptual model}\label{conceptual_model}

The following assumptions are considered:

\begin{enumerate}
\item There are two fluid phases, \emph{water} ($w$) and \emph{oil} ($o$),  and a solid ($s$) one. 
\item There are five components: \emph{water} ($w$), only in the water phase, \emph{oil} ($o$), only in the oil phase, \emph{microorganisms} ($m$) in the water phase (\textit{planktonic}) and in the solid phase (\textit{sessile}),  \emph{nutrients} ($n$) in the water phase and \emph{surfactant} ($\mathit{surf}$) produced by the microorganisms, in the water phase. 
\item The porous medium and the fluids are incompressible.
\item The porous medium is homogeneous, isotropic and fully saturated.
\item The fluid phases are separated in the pores and there is no diffusion between them.
\item All phases are in thermodynamical equilibrium.
\item Capillary pressure and relative permeability curves are taken in account. 
\item   Dynamical porosity and permeability variation  due to clogging$/$declogging processes is allowed.

\item Microorganisms, surfactant and nutrients dispersive fluxes follow the Fick's law.

\item Microorganisms and nutrients have biological interaction, as growth Monod equation, also known as  Michaelis-Menten kinetics \cite{Monod1949, Gaudy1980, Michaelis1913}.
\item Microorganisms' death is a first-order irreversible reaction \cite{Chang1991}.
\item Microorganisms have physico-chemical interaction with the porous medium, assuming reversible and irreversible clogging and reversible declogging \cite{Chang1991}.
\item The modification of porosity due to the process of clogging is expressed as in \cite{Corapcioglu1985} and \cite{Chang1992}.
\item 	The oil-water interfacial tension ($\sigma_{ow}$) is assumed to be a function of the surfactant concentration.

\end{enumerate}

After making such restrictions, it is made, for each model, a table of phases, their components and their extensive and intensive properties (see Tables \ref{Tabla_Prop_IntensivasF} and \ref{Tabla_Prop_IntensivasT}).

\begin{table}[h]
\begin{center}
\begin{tabular}{p{2.0cm} p{2.7cm} p{2.0cm} p{2.0cm} } \hline 
\textbf{Phase} & \textbf{Component} & \textbf{Extensive\newline property} & \textbf{Intensive\newline property} \\ \hline 
Water ($w$) & Water ($w$) & $M_{w}^{w} \left(t\right)$ & $\phi S_{w}\rho_{w}^{w}$ \\ 
Oil ($o$) & Oil ($o$) & $M_{o}^{o} \left(t\right)$ & $\phi S_{o}\rho_{o}^{o}$   \\ \hline 
\end{tabular}
\end{center}
\caption{Intensive properties associated with the mass of the components by phases for the flow model.}
\label{Tabla_Prop_IntensivasF}
\end{table}
 
 \begin{table}[h]
 \begin{center}
 \begin{tabular}{p{2.0cm} p{3.8cm} p{2.0cm} p{1.8cm} } \hline 
 \textbf{Phase} & \textbf{Component} & \textbf{Extensive\newline property} & \textbf{Intensive\newline property} \\ \hline 
 Water ($w$) & Microorganisms ($m$) & $M_{w}^{m} \left(t\right)$ & $\phi S_{w} c_{w}^{m} $ \\
  & Nutrients ($n$) & $M_{w}^{n} \left(t\right)$ & $\phi S_{w} c_{w}^{n} $ \\ 
   & Surfactants ($surf$) & $M_{w}^{surf} \left(t\right)$ & $\phi S_{w} c_{w}^{surf} $  \\ 
 Solid ($s$) & Microorganisms ($m$) & $M_{s}^{m} \left(t\right)$ & $\phi \sigma \rho _{m} $ \\
  \hline 
 \end{tabular}
 \end{center}
 \caption{Intensive properties associated with the mass of the components by phases for the transport model.}
 \label{Tabla_Prop_IntensivasT}
 \end{table}

\subsection{Mathematical model}

For deriving the equations of the mathematical model, the \emph{axiomatic formulation of continuum mechanics systems} \cite{Allen1988,Herrera2012} is applied. This axiomatic formulation adopts a macroscopic approach, which considers that the material systems are fully occupied by particles. A continuum system is made of a particle set known as material \textit{body}. The continuum system approach works with the volume average of the body properties and consequently there is a volume called \textit{representative elementary volume} (REV)  over which the property averages are valid.

Using the aforementioned approach, the equations of flow and transport were derived, under the assumptions specified in the conceptual model (section \ref{conceptual_model}).

\subsubsection[Flow model equations]{Flow model equations}

According to the axiomatic formulation of continuum mechanics systems applying the mass balance equations for oil and water components a general biphasic flow model is obtained 

\begin{equation} \label{OilEq} 
\left(\phi S_{w} \rho _{w} \right)_{t} +\nabla \cdot \left(\rho _{w} \underline{u}_{\,w} \right) = g_{w} ;
\end{equation} 
\begin{equation} \label{WaterEq} 
\left(\phi S_{o} \rho_{o} \right)_{t} +\nabla \cdot \left(\rho _{o} \underline{u}_{\,o} \right)= g_{o} ; 
\end{equation}  

where $\phi$ is the porosity, $S_{\alpha}$ is the saturation, $\rho _{\alpha }$ is the density, $\underline{u}_{\,\alpha}$ is the Darcy velocity and $g_\alpha$ is the mass source term, for $\alpha =w,o$.

Substituting in Eqns.(\ref{OilEq}-\ref{WaterEq}) the velocities with the Darcy law expression Eq. \ref{EqDarcyGravity} and dividing by densities, since the fluids are considered incompressible, results the following 

\begin{equation} \label{OilEq1} 
\left(\phi S_{w} \right)_{t} - \nabla \cdot \left( \frac{k_{rw } }{\mu _{w} } \underline{\underline{k}} \cdot \left(\nabla p_{w} +\rho _{w} {\mathrm \gamma_g}\nabla z\right) \right) = q_{w} ;
\end{equation} 
\begin{equation} \label{WaterEq1} 
\left(\phi S_{o} \right)_{t} - \nabla \cdot \left( \frac{k_{ro } }{\mu _{o} } \underline{\underline{k}} \cdot \left(\nabla p_{o} +\rho _{o} {\mathrm \gamma_g}\nabla z\right) \right)= q_{o} ; 
\end{equation} 

where the Darcy law can be expressed as follows 

\begin{equation} \label{EqDarcyGravity} 
\underline{u}_{\alpha } =-\frac{k_{r\alpha } }{\mu _{\alpha } } \underline{\underline{k}} \cdot \left(\nabla p_{\alpha } +\rho _{\alpha } {\mathrm \gamma_g}\nabla z\right);{\mathrm \; \; \; }\alpha ={\mathrm \; }o{\mathrm  , \; }w 
\end{equation} 

and  $\underline{\underline{k}}$ is the absolute permeability tensor, $\mu _{\alpha }$ is the viscosity,  $p_\alpha$ is the phase pressure, $k_{r\alpha }$ is the relative permeability, $q_\alpha=g_\alpha/\rho_\alpha$ is the volumetric source term,   $\gamma_g$ is the gravitational acceleration constant and $z$ is the elevation.\\
As can be observed, the equation system \ref{OilEq1}-\ref{WaterEq1} consists of only two equations and it contains four unknowns $(p_o, p_w,S_o,S_w)$, therefore the following two additional equations are required for the system to be determined:
\begin{eqnarray}
S_o+S_w &=& 1\\\nonumber
p_{cow} &=& p_o-p_w\\\nonumber
\end{eqnarray}
where $p_{cow}$ is the oil-water capillary pressure.\\

For convenience, the following notation is introduced: for $\alpha=w,o$, $\lambda_\alpha =k_{r\alpha}/\mu_\alpha$ are the phase mobility functions, $\lambda = \sum\lambda_\alpha$ is the total mobility and $f_\alpha =\lambda_\alpha/\lambda$ are the fractional flow functions. So that $\sum f_\alpha =1$. The total velocity is defined as $\underline{u}=\sum\underline{u}_\alpha$.
Since the oil is a continuous phase and consequently its pressure is well behaved, we are going to define the oil phase pressure $p_o$ as the pressure variable $p$.\\ 

Rewriting the equations (\ref{OilEq1}-\ref{WaterEq1}) applying the oil phase pressure and total velocity formulation given in \cite{Chen2000}, in which the capillary pressure and relative permeability curves are taken in account, results the biphasic flow model \cite{Diaz2008a,Diaz2015b}: 

\begin{itemize} 

\item Pressure equation

\begin{eqnarray} \label{EqOilPressure} 
 -\nabla \cdot 
 \left\{\lambda \underline{\underline{k}} \cdot \nabla p_{o} -\left(\lambda _{w} \frac{dp_{cow} }{dS_{w} } \right)\underline{\underline{k}} \cdot \nabla S_{w} \right\} \\
-\nabla \cdot 
 \left\{\left(\lambda _{o} \rho _{o} +\lambda _{w} \rho _{w} \right){\mathrm \gamma_g }\underline{\underline{k}} \cdot \nabla z\right\}
+\frac{\partial \phi }{\partial t}=q_{o} +q_{w} \nonumber  
\end{eqnarray}

\item Saturation equation

\begin{eqnarray} \label{EqWaterSaturation} 
\phi \frac{\partial S_{w} }{\partial t} -\nabla \cdot \left\{\lambda _{w} \underline{\underline{k}} \cdot \nabla p_{o} -\left(\lambda _{w} \frac{dp_{cow} }{dS_{w} } \right)\underline{\underline{k}} \cdot \nabla S_{w} \right\}  \\ \nonumber
 -\nabla \cdot \left\{\left(\lambda _{w} \rho _{w} {\mathrm \gamma_g }\right)\underline{\underline{k}} \cdot \nabla z\right\} +\left(\frac{\partial \phi }{\partial t} \right)S_{w} =q_{w} \nonumber  
\end{eqnarray} 
 
\item Velocity equations

\begin{equation} \label{velocities} 
\begin{array}{l} {\underline{u}=-\lambda \underline{\underline{k}} \cdot \nabla p_{o} +\lambda f_{w} \left(\frac{dp_{cow} }{dS_{w} } \right)\underline{\underline{k}} \cdot \nabla S_{w} -\lambda \left(f_{o} \rho _{o} +f_{w} \rho _{w} \right){\mathrm \gamma_g }\underline{\underline{k}} \cdot \nabla z;} \\ {\underline{u}_{w} =-\lambda f_{w} \underline{\underline{k}} \cdot \nabla p_{o} +\lambda f_{w} \left(\frac{dp_{cow} }{dS_{w} } \right)\underline{\underline{k}} \cdot \nabla S_{w} -\lambda f_{w} \rho _{w} {\mathrm \gamma_g }\underline{\underline{k}} \cdot \nabla z;} \\ {\underline{u}_{o} =-\lambda f_{o} \underline{\underline{k}} \cdot \nabla p_{o} -\lambda f_{o} \rho _{o} {\mathrm \gamma_g}\underline{\underline{k}} \cdot \nabla z;} \end{array} 
\end{equation} 

 where $p_{cow}$ is the oil-water capillary pressure,  $\lambda$ is the total mobility, and $\lambda_\alpha$ is the phase mobility.
 
 \end{itemize}
 
Notice that phase velocities can be known once the pressure and saturation equations are solved. The modification of the porosity due to the clogging/declogging processes is also taken in account. 
 
\subsubsection[Transport model equations]{Transport model equations}
Following the aforementioned axiomatic formulation of continuum mechanics systems applying the mass balance equation for a flowing component (${\eta}$) in a water phase ($w$) a generic transport equation was derived (Eq. \ref{Eq_D4}), from which can be instantiated the specific transport equations that were used 

\begin{equation} \label{Eq_D4} 
\begin{array}{l} {\left(\phi S_{w} c_{w}^{\eta} \right)_{t} -\nabla \cdot \left(\phi S_{w} \underline{\underline{D}}_{\, w}^{\eta}  \cdot \nabla  c_{w}^{\eta}\right)+\nabla \cdot \left(\underline{u}_{w} c_{w}^{\eta} \right)=} {g_{w}^{\eta}}
  \end{array} 
\end{equation} 
where in the advection term $\underline{u}_{w}$ is Darcy's velocity for the water phase (Eq. \ref{EqDarcyGravity}); $\phi S_{w} \underline{\underline{D}}_{\, w}^{\eta}  \cdot \nabla  c_{w}^{\eta}$ is the dispersive flux term for the ${\eta}^\mathrm{th}$ component (Fick's Law) (Eq. \ref{Eq_Ficks_law}) and $g_{w}^{\eta}$ is a generic source term.

%
\begin{itemize} 

\item Mass balance equation of the \emph{planktonic} microorganisms component ($m$) in the water phase ($w$) \cite{Li2011,Zhang1995}

Making $g_{w}^{\eta} = \phi S_{w} \left(g_{m} -d_{m} -k_{d} \right)c_{w}^{m} +k_{d} \rho _{m} \sigma _{r} -{R_{\mathit{surf}}  \mathord{\left/ {\vphantom { }} \right. \kern-\nulldelimiterspace} Y_{\mathit{surf}/m} }$ in Eq. (\ref{Eq_D4}), it is obtained:

\begin{equation} \label{Eq_D4b} 
\begin{array}{l} {\left(\phi S_{w} c_{w}^{m} \right)_{t} -\nabla \cdot \left(\phi S_{w} \underline{\underline{D}}_{\,w}^{m }  \cdot \nabla  c_{w}^{m }\right)+\nabla \cdot \left(\underline{u}_{\,w} c_{w}^{m} \right)=} \\ {\phi S_{w} \left(g_{m} -d_{m} -k_{d} \right)c_{w}^{m} +k_{d} \rho _{m} \sigma _{r} -{R_{\mathit{surf}}  \mathord{\left/ {\vphantom { }} \right. \kern-\nulldelimiterspace} Y_{\mathit{surf}/m} }  }
 \end{array} 
\end{equation}

\item Mass balance equation of the nutrients component ($n$) in the water phase ($w$) \cite{Zhang1995}

Making 

\begin{equation} \label{Eq_omega_nutrients)} 
\begin{array}{l} {g_{w}^{\eta} =} \\ {- g_{m}^{} \left(\phi S_{w} c_{w}^{m} +\sigma \rho _{m} \right)/Y_{m/n} -{R_{\mathit{surf}}  \mathord{\left/ {\vphantom { }} \right. \kern-\nulldelimiterspace} Y_{\mathit{surf}/n} } } \\ {- m_{n} \left(\phi S_{w} c_{w}^{m} +\sigma \rho _{m} \right) } \end{array}
\end{equation}

in Eq. (\ref{Eq_D4}), it is obtained: 

\begin{equation} \label{Eq_D4c} 
\begin{array}{l} {\left(\phi S_{w} c_{w}^{n} \right)_{t} -\nabla \cdot \left(\phi S_{w} \underline{\underline{D}}_{\,w}^{n}  \cdot \nabla  c_{w}^{n }\right)+\nabla \cdot \left(\underline{u}_{\,w} c_{w}^{n} \right)=} \\ {-g_{m}^{} \left(\phi S_{w} c_{w}^{m} +\sigma \rho _{m} \right)/Y_{m/n} -{R_{\mathit{surf}}  \mathord{\left/ {\vphantom { }} \right. \kern-\nulldelimiterspace} Y_{\mathit{surf}/n} } } \\ {-m_{n} \left(\phi S_{w} c_{w}^{m} +\sigma \rho _{m} \right) } \end{array}
\end{equation}

\item Mass balance equation of the surfactant component ($\mathit{surf}$) in the water phase ($w$) 

Making  $g_{w}^{\eta} =R_{\mathit{surf}}$ in Eq. (\ref{Eq_D4}), it is obtained:  

 \begin{equation} \label{Eq_D4d} 
\begin{array}{l} 
{\left(\phi S_{w}c_{w}^{\mathit{surf}}\right)_{t} -\nabla \cdot \left(\phi S_{w} \underline{\underline{D}}_{\;w}^{\mathit{surf}}  \cdot \nabla  c_{w}^{\mathit{surf}}\right) +\nabla \cdot \left(\underline{u}_{\,w} c_{w}^{n} \right) = R_{\mathit{surf}}} \\ 
 \end{array} 
\end{equation} 

In the same manner applying the mass balance equation for a component (${\eta}$) in a solid phase ($s$) a generic equation was obtained

 \begin{equation} \label{Eq_D5} 
 \left(c_{s}^{\eta} \right) _t=g_{s}^{\eta}  
\end{equation}

\item Mass balance equations of the \emph{sessile} microorganisms component ($m$) in the solid phase ($s$) \cite{Li2011,Zhang1995}

Expressing the time derivative of Eq. (\ref{Eq_D5}) $\left(c_{s}^{\eta} \right) _t = \left(\rho _{m} \sigma _{r} \right) _t$, and making 

\begin{equation}
g_{w}^{\eta} = \left(g_{m} -k_{d} -d_{m} \right)\rho _{m} \sigma _{r} +k_{c1} \phi c_{w}^{m}
\end{equation}

a mass balance equation for reversible clogging is obtained

 \begin{equation} \label{Eq_D5b} 
 \left(\rho _{m} \sigma _{r} \right) _t=\left(g_{m} -k_{d} -d_{m} \right)\rho _{m} \sigma _{r} +k_{c1} \phi c_{w}^{m}  
\end{equation}

And in a similar way  expressing $\left(c_{s}^{\eta} \right) _t = \left(\rho _{m} \sigma _{i} \right) _t$, and making

\begin{equation}
g_{w}^{\eta} =\left(g_{m} -d_{m} \right)\rho _{m} \sigma _{i} +k_{c2} \phi c_{w}^{m} 
\end{equation}  

a mass balance equation for irreversible clogging is obtained

\begin{equation}\label{Eq_D5c}
 \left(\rho _{m} \sigma _{i} \right) _t=\left(g_{m} -d_{m} \right)\rho _{m} \sigma _{i} +k_{c2} \phi c_{w}^{m}  
\end{equation}

where $c_{w}^{\eta}$ is the concentration, $\underline{\underline{D}}_{\; w}^{\eta}  $ is the hydrodynamic dispersion tensor,  $D^{*\eta}_{\;\;w}$ is the molecular diffusion coefficient of microorganisms, nutrients or surfactant, respectively ($\eta = m,\,n,\,\mathit{surf}$), in the water phase ($w$); $d_{m} $ is the decaying rate (death) of microorganisms; $k_{d} $ is the declogging rate; $\rho _{m} $ is the microbial density at surface conditions $g_{m} $ is the microbial rate growth; $m_{n} $ is the energy coefficient for maintaining life through substrate consumption; $\sigma,\; \sigma _{r},\;\sigma _{i}$ are the total, reversible and irreversible volume fraction (respectively) occupied by the settled microorganisms over total pore volume; $k_{c1} $ is the reversible clogging rate; $k_{c2} $ is the irreversible clogging rate; $Y_{m/n}$ is the yield coefficient of microorganism per unit of nutrients; $Y_{{\mathit{surf} \mathord{\left/ {\vphantom {\eta m}} \right. \kern-\nulldelimiterspace} \eta} } $ is the yield coefficient of surfactant per unit of microorganism or nutrients ($\eta = m,\,n$); $R_{\mathit{surf}}$ is a source term for the surfactant, and is given by \cite{Zhang1995}
\begin{equation}
R_{\mathit{surf}} =\mu _{\mathit{surf}}^{\max } \left(\frac{c_{w}^{n} -c_{w}^{nC} }{K_{\mathit{surf}/n} +c_{w}^{n} -c_{w}^{nC} } \right)\left(\phi S_{w} c_{w}^{m} +\sigma \rho _{m} \right),
\end{equation}

where $\mu _{\mathit{surf}}^{\max } $ is the maximal surfactant formation rate; $K_{\mathit{surf}/n} $ is the saturation constant for surfactant over nutrient; $c_{w}^{nC} $ is the critical nutrient ($n$) concentration   for the surfactant ($\mathit{surf}$) formation.

\end{itemize}
\subsubsection{Initial and boundary conditions}\label{InitialBoundary}
\begin{itemize} 
\item {Initial conditions:}
\begin{equation}
\begin{array}{r}
p_{o}\left( {t_0 } \right) = p_o^{0}, \;S_{w} \left( {t_0 } \right) = S_{w}^{0};   \\
\sigma\left( {t_0 } \right)  = \sigma_{0}, {c^\eta_w} \left( {t_0 } \right)   = c^{\eta0}_w, \;\eta=m,n,\mathit{surf} ;
\end{array}
\end{equation}
where $p_o^0$ is the initial oil pressure,  $k_{0} $ is the initial absolute permeability;  $\phi_{0}$ is the initial porosity of porous medium; $\sigma_{0}$ is the initial volume fraction occupied by the settled microorganisms over total pore volume.

\item {Boundary conditions}
\begin{enumerate}
\item Inlet conditions (constant rate)
\begin{equation}
\begin{array}{r}
 {\underline{u}}_{\,o} \cdot {\underline{n}} = {\underline{u}}_{w} \cdot {\underline{n}} =  {\underline{u}}^{in}_{w} \cdot {\underline{n}}; \;\;\;   c^{\mathit{surf}_{in} }_w=0; \\
 - \left[c^\eta_w{\underline{u}}^{in}_{\,w}  -\phi S_w{\underline{\underline D}}^\eta_{\,w}\cdot\nabla c^\eta_w\right] \cdot {\underline{n}}=c^{\eta _{in} }_w{\underline{u}}^{in}_{w} \cdot {\underline{n}}, \;\eta=m,n,\mathit{surf} ;
 \end{array}
\end{equation}

\item Outlet conditions (constant pressure, $p_o^{out}$ at the top of the core)
\begin{equation}
p_{o} =p^{out}_{o},\; \frac{{\partial S_{w} }}{{\partial \underline{n}}} = 0;
\frac{{\partial c^\eta_w }}{{\partial \underline{n}}} = 0, \;\eta=m,n,\mathit{surf} ;
\end{equation}

\item No flow conditions for all not specified boundaries.

\end{enumerate}

\end{itemize}

\subsubsection[Complementary relationships]{Complementary relationships}

\begin{itemize} 
\item It is used the Fick's law: $\underline{\tau}_w^\eta=\phi S_w \underline{\underline D}\, w^\eta\cdot \nabla c_w^\eta ,\, \eta = m,\, n$, where $c_w^\eta$ is the concentration of the component $\eta$ in water (mass of the component $\eta$ per volume of water), and

\begin{equation}\label{Eq_Ficks_law}
\left(D_{w}^{\eta } \right)_{ij} =\left(\alpha _{T} \right)_{w}^{\eta } \left|\underline{{\mathrm v}}_{w} \right|\delta _{ij} +\left(\left(\alpha _{L} \right)_{w}^{\eta } -\left(\alpha _{T} \right)_{w}^{\eta } \right)\frac{{\mathrm v}_{wi} {\mathrm v}_{wj} }{\left|\underline{{\mathrm v}}_{w} \right|} +\tau \left(D_{m} \right)_{w}^{\eta } \delta _{ij} 
\end{equation}

are the hydrodynamic dispersion tensor components ($i,j=1,2,3$) \cite{Bear1988}, where $\left(\alpha _{T} \right)_{w}^{\eta } $ and $\left(\alpha _{L} \right)_{w}^{\eta } $ are the transversal and longitudinal dispersivity of component $\eta $ in water, respectively. $\left|\underline{{\mathrm v}}_{w} \right|$ is the Euclidean norm of the water velocity vector, $\delta _{ij} $ is the Kronecker delta function, $\tau <1$ is the tortuosity and $\left(D_{m} \right)_{w}^{\eta } $ is the molecular diffusion coefficient of the component $\eta $ in water.

\item It is used the growth Monod equation: \newline
$ g_m  = g_m^{\max } \left( \frac{c^n_w }{K_{m/n}  + c^n_w } \right)$, where $g_m^{\max }$ is the maximum specific growth rate, $K_{m/n} $ is the Monod constant for nutrients, $c^n_w $ is the nutrients concentration in water.

\item A Brooks-Corey \cite{Brooks1964} model for capillary pressure is used (Eq. \ref{Eq_A3}).

\begin{equation}\label{Eq_A3}
p_{cow}(S_w)=p_t\left( S_e \right)^{\left(-1/\theta\right)}
\end{equation}

and

\begin{equation}\label{Eq_A3a}
S_e=\frac{S_w-S_{wr}}{1-S_{wr}-S_{or}}
\end{equation}

where $S_e$ is the normalized water saturation, $p_t$ is the threshold input pressure and $\theta $ is the capillary pressure exponent for the Brooks-Corey model.

\item A modified Brooks-Corey  \cite{Lake1989} model for relative permeability curves is used:

\begin{equation}\label{Eq_A2}
\begin{array}{l} 
{k_{rw} = k_{rw}^{0} S_e^{n_w}} \\ 
{k_{ro} = k_{ro}^{0} \left( 1-S_e  \right)^{n_o}} 
\end{array}
\end{equation}

where $S_e$ is the normalized water saturation, $k_{rw}^{0}$, $k_{ro}^{0}$ are the endpoints and $n_w$, $n_o$ are the exponents for water $w$ and oil $o$ relative permeability curves, respectively.

\item Porosity modification:
The porosity modification due to the clogging/ declogging processes was taken in account by the following expression given in \cite{Chang1992}:
\begin{equation}
\begin{array}{r}
\phi  = \phi _0  - \sigma 
\end{array}
\label{porosity} 
\end{equation}
where $\phi$ - is the actual porosity, $\phi _0$ - is the initial porosity and $\sigma$ - the volume fraction occupied by sessile microorganisms.

\item {Permeability modification:}
The permeability modification is expressed as a porosity function by the Kozeny-Carman equation \cite{Carman1956}:
\begin{equation}
\begin{array}{r}
k = k_0 \frac{{\left( {1 - \phi _0 } \right)^2 }}{{\phi _0^3 }}\frac{{\phi _{}^3 }}{{\left( {1 - \phi } \right)^2 }}
\end{array}
\label{permeability} 
\end{equation}
where $k$ and $k_0$ are the actual and initial permeability, respectively.
\item {Interfacial tension:}
The change of interfacial tension due to the surfactant concentration has been measured in the laboratory, and its behavior is represented approximately by a best fit curve obtained with experimental data:
\begin{equation}\label{Eq_3_3_3}
\sigma _{ow}  = \frac{1}{{1000}}\sqrt {\frac{{10}}{{\left( {c_{w}^{surf}  + c_{w}^{crit} } \right)}}} 
\end{equation}

where $c_{w}^{surf}  = 0 \Rightarrow \sigma _{ow}^{\max }  = \frac{1}{{1000}}\sqrt {\frac{{10}}{{\left( {c_{w}^{crit} } \right)}}}  = {\mathrm{35}}.{\mathrm{59 }}\left[ {{\mathrm{mN}}/{\mathrm{m}}} \right]$

The interfacial tension oil-water $\sigma_{ow}$ is given in ${{\mathrm{mN}}/{\mathrm{m}}}$ and the surfactant concentration $c_{w}^{surf}$ in $\mathrm{kg/m^{3}}$
The original fitting is shown in Fig. \ref{FigFitSigmaOW}.
\begin{figure}[htb]
    \centering
    \includegraphics[width=0.6\textwidth]{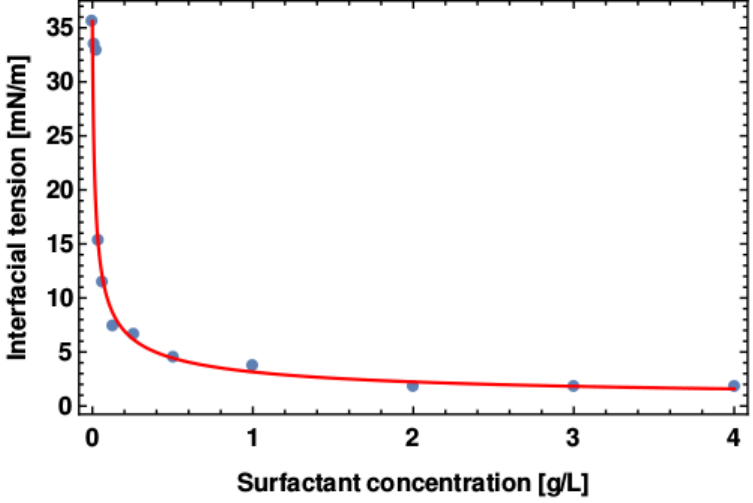}
  \caption{Best fit curve of experimental data for the interfacial tension as a function of surfactant concentration.}
\label{FigFitSigmaOW}           
\end{figure} 
\item {Trapping number:}

Following the formulation given by Penell et al. \cite{li2007experimental,pennell1996influence} the trapping number of the oil phase is defined as follows: 

\begin{equation}\label{Eq_3_3_4}
N_{To} =\sqrt{N_{CA}^2+2N_{CA}N_{B}\sin(\alpha)+N_B^2}
\end{equation}
where $N_{To}$ is the oil phase trapping number, $N_{CA} = (u_w^{in} \mu_w) /(\sigma_{ow} \cos(\theta_{ow}))$ is the capillary number, $N_B=\Delta \rho \gamma_g k k_{rw}/(\sigma_{ow} \cos(\theta_{ow}))$ is the Bond number, $\theta_{ow}$ is the contact angle between the aqueous/non-aqueous interface and the porous medium, $\alpha$ is the angle of flow relative to the horizontal, and $\gamma_g$ is the gravitational constant. 

Since in our experiments we are considering only vertical flow, then  $\alpha=\pi/2$, and consequently Eq. (\ref{Eq_3_3_4}) becomes $N_{To} =\mid N_{CA}+N_B\mid$. Moreover, the Bond number is two orders of magnitude smaller than the capillary number so it can be neglected  ($N_B=0$). Since the velocity of the water phase ($u_w$) varies relatively little, to avoid instabilities, it was set equal to the injection velocity ($u^{in}_w$), the trapping number finally results in

\begin{equation}\label{Eq_3_3_4a}
N_{To} =N_{CA} = \frac{{{u^{in}_w}{\mu _w}}}{{{\sigma _{ow}}\cos \theta_{ow} }}
\end{equation}



\item {Change in the residual oil saturation:}

The change in the residual oil saturation is a function of the trapping number, which can be expressed as a linear interpolation as follows:

\begin{equation}\label{Eq_3_3_5}
S_{or} =S_{or}^{low} -\left(S_{or}^{low} -S_{or}^{high} \right)\frac{\left(N_{To} -N_{To}^{low} \right)}{\left(N_{To}^{high} -N_{To}^{low} \right)}
\end{equation}
where $S_{or}$ is the residual oil phase  saturation and
\[ S_{or} =S_{or}^{low} \quad if \quad S_{or} >S_{or}^{low} \quad  \mathrm{and}  \quad S_{or} =S_{or}^{high} \quad if\quad S_{or} <S_{or}^{high} \]

Here the index $low$ is associated with a low trapping number and the index $high$ is associated with a high trapping number.

\item {Change of relative permeability}: 
 The change of relative permeability for the modified Brooks-Corey model (Eq. \ref{Eq_3_3_6}) consist of the modification of the \emph{endpoint} $k_{ro}^{0}$ and the exponent $n_{o}$ as follows:

\begin{equation}\label{Eq_3_3_6}
k_{ro}^{0} =k_{ro}^{0\, low} +\frac{S_{or}^{low} -S_{or} }{S_{or}^{low} -S_{or}^{high} } \left(k_{ro}^{0\, high} -k_{ro}^{0\, low} \right)
\end{equation}

\begin{equation}\label{Eq_3_3_7a}
n_{o} =n_{o}^{low} +\frac{S_{or}^{low} -S_{or} }{S_{or}^{low} -S_{or}^{high} } \left(n_{o}^{high} -n_{o}^{low} \right)
\end{equation}

\item Note that the recovery mechanism connected with the surfactant concentration has been taken in account through an empiric function (Eq. \ref{Eq_3_3_3}) of the interfacial tension, which modifies the trapping number (Eq. \ref{Eq_3_3_4}), and then the residual oil saturation is modified (Eq. \ref{Eq_3_3_5}), which in turn modifies the oil-water capillary pressure, the endpoint (Eq. \ref{Eq_3_3_6}) and the exponent (Eq. \ref{Eq_3_3_7a}) of the oil relative permeability. All this can be expressed schematically as in Eq.(\ref{Eq_3_3_7}) or as in Fig. \ref{Figura_ch_p_diag}.
\begin{equation}\label{Eq_3_3_7}
C_{\mathit{surf}}\Rightarrow \sigma _{ow} \Rightarrow N_{T_o}\Rightarrow S_{or}\Rightarrow  \lbrace p_{cow}, k_{ro}^{0}, n_{o}  \rbrace \Rightarrow  \lbrace k_{ro}  \rbrace
\end{equation}

\begin{figure}[htb]
    \centering
    \includegraphics[width=1.0\textwidth]{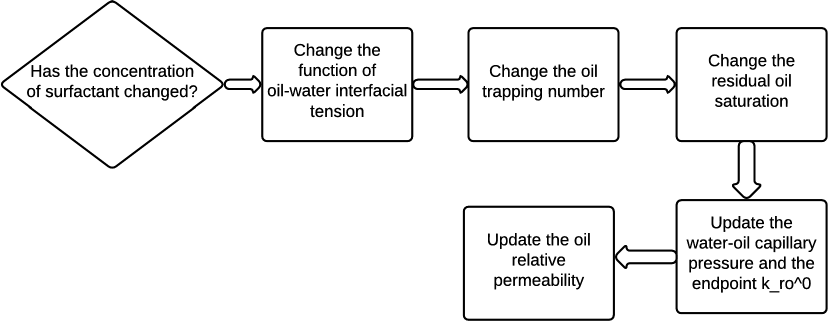}
  \caption[Flow chart of changes in the interfacial tension phenomenon, and the values of the relevant parameters]{Flow chart of changes in the interfacial tension phenomenon, and the change in the corresponding values of the relevant parameters}
\label{Figura_ch_p_diag}           
\end{figure}

\end{itemize}
\subsection{Numerical model}

The numerical model consists of making the appropriate choice of the numerical methods in terms of precision and the efficiency for the solution of the mathematical model. \\

A brief review of the state of the art literature concerning the numerical implementation of multiphasic fluid flow model reflects that the finite difference (FD) and finite volume (FV) methods are the general framework for numerical simulation in very large problems \cite{Chen2006}; however, the basic mixed finite element (MFE) method \cite{raviart1977mixed} has shown to be superior for accurate flux calculation in heterogeneous media in comparison to conventional FD and FV methods. 

On the other hand, a mixed finite element approach requires a special Raviart-Thomas mixed space for base and weighting functions, which makes more difficult its implementation. In view of the scale and resolution requirements for our flow model, we decided that it could be acceptable to perform the implementation making use of the standard finite element framework provided in COMSOL Multiphysics \cite{COMSOL2008}. In particular, the numerical implementation of previously derived model was accomplished using the PDE mode for time dependent analysis in the coefficient form.

\begin{figure}[h!]
    \centering
    \includegraphics[width=0.35\textwidth]{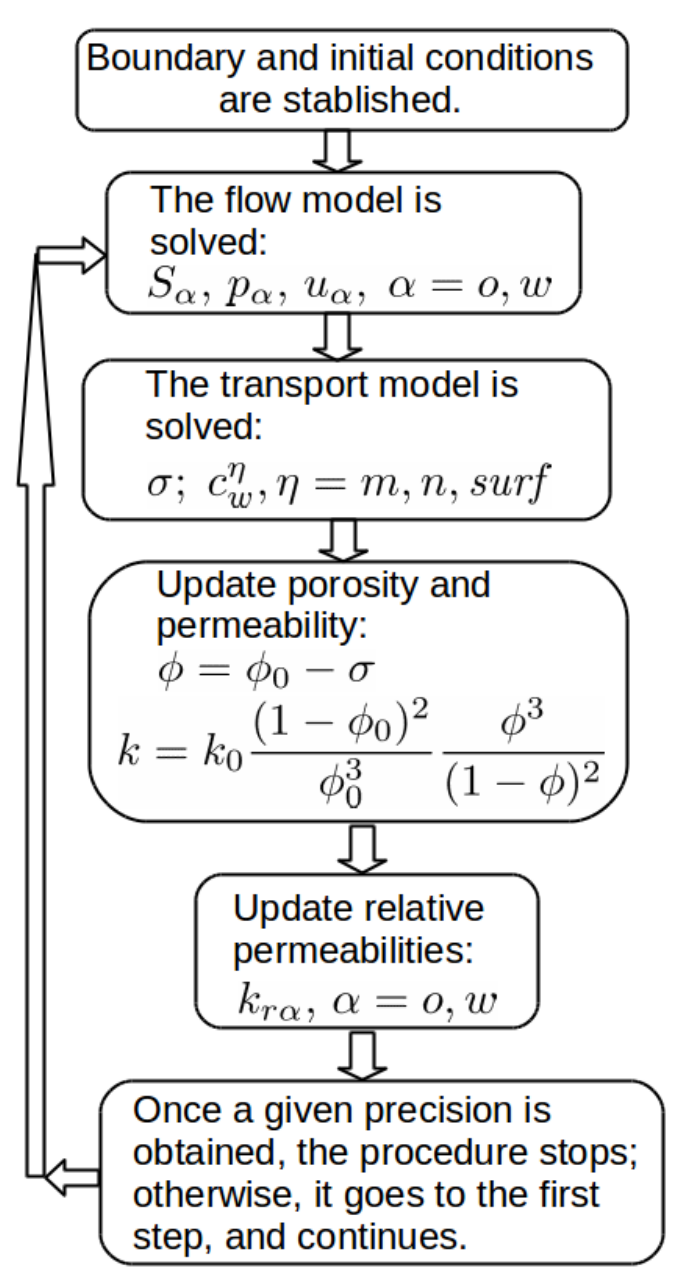}
  \caption[Flow chart of the numerical model]{Flow chart of the numerical model}
\label{FigFlowChartNumModel}           
\end{figure}

In this case, the resulting problem is a system of nonlinear partial differential equations with initial and boundary conditions. For the numerical solution the following methods were applied:
\begin{itemize} 

\item For the time derivative, it was used a second order backward finite differences discretization, resulting in a totally implicit scheme in time.

\item For the rest of the differential operators, concerning the spatial derivatives, it was applied a standard Galerkin finite element discretization, where Lagrange quadratic polynomials were used as weighting and basis functions, which in this work imply a convergence of order two \cite{zienkiewicz1977finite}.

\item It was used a regular mesh of tetrahedral elements in 3D.

\item For the linearization of the nonlinear system of equations, an iterative Newton-Raphson method was applied.

\item For the solution of the resulting algebraic system of linear equations, it was used a variant of the direct LU method for sparse, unsymmetrical matrices. 

\item The general procedure for coupling flow and transport is sequential and is executed iteratively in the following manner:
\begin{enumerate}
\item The flow model is solved, and from them it is obtained: saturations $S_{\alpha}$, pressures $p_{\alpha}$ and the velocities of the phases $\underline{u}_{\,\alpha }$ for $\alpha=o,w$,
\item The transport model is solved and from them it is obtained: the concentration  components  $c_{w}^{\eta },\,\eta =m,n,\mathit{surf}$,  and $\sigma$,
\item The porosity $\phi$ and permeability $k$ are modified according to equations (\ref{porosity}) and (\ref{permeability}).
\item The relative permeabilities $k_{r\alpha }$ for $\alpha=o,w$ are modified according to equations (\ref{Eq_3_3_5}-\ref{Eq_3_3_7}).
\item Once a given precision is obtained, the procedure stops; otherwise, it goes over the first step, and continues.
\end{enumerate}

\end{itemize}
All the above description is summarized in Fig. \ref{FigFlowChartNumModel}

\subsection{Computational model}

The computational implementation was done in COMSOL multiphysics, using the general coefficient form \cite{COMSOL2008,Diaz2015b}) for all the equations.
This computational implementation in COMSOL, of the flow and transport models, has been thoroughly validated in other previous works \cite{Diaz2008,Diaz2012,Diaz2012b,Diaz2015,Diaz2015b}, with data taken from the literature.

The system of equations (\ref{OilEq1}-\ref{WaterEq1}) can be rewritten in matrix form as follows:
\begin{equation}\label{18a)}
\begin{aligned} 
\begin{pmatrix}
{0} & {0} \\
{0} & {\phi }\\
\end{pmatrix}
\frac{\partial }{\partial t}
\begin{pmatrix}
{p_{o} } \\
{S_{w} }\\
\end{pmatrix}\\
\begin{split}
+\nabla \cdot \left\{- \left(\begin{array}{cc} {\lambda } & {-\lambda _{w} \frac{dp_{cow} }{dS_{w} } } \\ {\lambda _{w} } & {-\lambda _{w} \frac{dp_{cow} }{dS_{w} } } \end{array}\right)\underline{\underline{k}} \cdot \nabla \left(\begin{array}{c} {p_{o} } \\ {S_{w} } \end{array}\right)\right.\\
\left. +\left(\begin{array}{l} {-\left(\lambda _{o} \rho _{o} +\lambda _{w} \rho _{w} \right) \gamma_g \underline{\underline{k}} \cdot \nabla z} \\ {-\lambda _{w} \rho _{w}  \gamma_g \underline{\underline{k}} \cdot \nabla z} \end{array}\right)\right\}
\end{split}
\\
+ \begin{pmatrix}
0 & 0 \\
0 & \partial \phi / \partial t\\
\end{pmatrix}
\begin{pmatrix}
{p_{o} } \\
{S_{w} }\\
\end{pmatrix}
= \begin{pmatrix}
q_o+q_w-\partial \phi / \partial t\\
q_w\\
\end{pmatrix}
\end{aligned}  
\end{equation} 

The previous matrix representation can be translated in straightforward manner to the standard COMSOL notation in coefficient form \cite{COMSOL2008,Diaz2015b}):

\begin{equation} \label{18b)} 
{\rm e}_{{\rm a}} \frac{\partial ^{2} \mathbf{u}}{\partial t^{2} } {\rm +d}_{{\rm a}} \frac{\partial \mathbf{u}}{\partial t} + \nabla \cdot \left(-{\rm c}\nabla \mathbf{u}-{\rm \alpha }\mathbf{u}+{\rm \gamma }\right)+{\rm \beta }\cdot \nabla \mathbf{u}+{\rm a}\mathbf{u}={\rm f} 
\end{equation} 

where

\begin{eqnarray}\label{18c)} 
\mathbf{u}& \equiv & \left(\begin{array}{c} {p_{o}} \\ {S_{w} } \end{array}\right),\,{\rm d}_{{\rm a}} \equiv \left(\begin{array}{cc} {0} & {0} \\ {0} & {\phi } \end{array}\right),\,{\rm c}\equiv \left(\begin{array}{cc} {\lambda } & {-\lambda _{w} \frac{dp_{cow} }{dS_{w} } } \\ {\lambda _{w} } & {-\lambda _{w} \frac{dp_{cow} }{dS_{w} } } \end{array}\right)\underline{\underline{k}},\nonumber\\
{\rm \gamma }& \equiv &\left(\begin{array}{l} {-\left(\lambda _{o} \rho _{o} +\lambda _{w} \rho _{w} \right) \gamma_g\underline{\underline{k}} \cdot \nabla z} \\ {-\lambda _{w} \rho _{w}  \gamma_g \underline{\underline{k}} \cdot \nabla z} \end{array}\right), \nonumber\\
 {\rm a}&\equiv &\left(\begin{array}{cc} {0} & {0} \\ {0} & {{\partial \phi \mathord{\left/ {\vphantom {\partial \phi  \partial t}} \right. \kern-\nulldelimiterspace} \partial t} } \end{array}\right),\, {\rm f}\equiv \left(\begin{array}{c} {q_{o} +q_{w} -{\partial \phi \mathord{\left/ {\vphantom {\partial \phi  \partial t}} \right. \kern-\nulldelimiterspace} \partial t} } \\ {q_{w} } \end{array}\right)\nonumber
\end{eqnarray}
and ${\rm e}_{{\rm a}} {\rm ,\alpha ,\beta }\equiv 0$.

To complete the model only remain to define suitable constitutive laws for relative permeabilities ($k_{rw} $, $k_{ro} $) and oil-water capillary pressure $p_{cow} $ and to prescribe proper initial and boundary conditions.

For defining initial and boundary conditions in those problems the following notations will be introduced:

\noindent \textbf{Initial conditions}
\begin{equation} \label{20)} 
p\left(t_{0} \right)=p_{0} ,{\rm \; \; }S_{w} \left(t_{0} \right)=S_{w0}  
\end{equation} 
\textbf{Boundary conditions}

Using the standard COMSOL notation in coefficient form for boundary conditions

\noindent 
\begin{equation} \label{20a)} 
\begin{array}{l} {\mathbf{n}\cdot \left({\rm c}\nabla \mathbf{u}+{\rm \alpha }\mathbf{u}-{\rm \gamma }\right)+{\rm q}\mathbf{u}={\rm g}-{\rm h}^{T} {\rm \mu }} \\ {{\rm h}\mathbf{u}=\rm r} \end{array} 
\end{equation} 

\begin{equation} \label{21)} 
\begin{array}{l} {{\rm r}^{in} \equiv \left(\begin{array}{c} {p^{in} } \\ {S_{w}^{in} } \end{array}\right),{\rm g}^{in} \equiv \left(\begin{array}{c} {g_{p}^{in} } \\ {g_{S_{w} }^{in} } \end{array}\right){\rm ,}{\rm h}^{in} \equiv \left(\begin{array}{cc} {{\rm h_{11}^{in}}} & {{\rm h_{12}^{in}}} \\ {{\rm h_{21}^{in}}} & {{\rm h_{22}^{in}}} \end{array}\right)} \\ 
{{\rm r}^{out} \equiv \left(\begin{array}{c} {p^{out} } \\ {S_{w}^{out} } \end{array}\right),{\rm g}^{out} \equiv \left(\begin{array}{c} {g_{p}^{out} } \\ {g_{S_{w} }^{out} } \end{array}\right){\rm ,}{\rm h}^{out} \equiv \left(\begin{array}{cc} {{\rm h_{11}^{out}}} & {{\rm h_{12}^{out}}} \\ {{\rm h_{21}^{out}}} & {{\rm h_{22}^{out}}} \end{array}\right)}  \end{array}  
\end{equation} 
where ${\rm q}\equiv 0$, but ${\rm h}$ depends on the type of specific  boundary conditions. Here $in$ and $out$ stand for \emph{inlet} and \emph{outlet} boundary conditions, recpectively.

%
\section{Validation of the flow model}\label{SeccValidationFlowModel}

In the next two subsections the flow model described above will be tested for two specific problems in 1-D. 

\subsection{Buckley-Leverett Problem}
We first will verify the implemented numerical flow model with known analytical solutions. To this end, we solve the Buckely--Leverett problem in a homogeneous medium with different fluid properties and zero capillary pressure \cite{Buckley1942}. 

We consider a 1-D horizontal homogeneous domain of length 300 m, initially saturated with oil. Water is injected with a constant flow rate at one end to displace oil to the other end, where the pressure is kept constant. 

The relative permeability constitutive equations are given by:
\begin{equation} \label{22)} 
k_{rw} =S_{e}^{\omega } {\rm ;\; \; \; }k_{ro} =\left(1-S_{e} \right)^{\omega } ; 
\end{equation} 
where $\omega =1$ is for the linear case and $\omega =2$ is for the quadratic case,  whereas $S_{e} $ is the effective or normalized saturation, which is defined as:
\begin{equation} \label{23)} 
S_{e} =\frac{S_{w} -S_{wr} }{1-S_{wr} -S_{or} }  
\end{equation} 
where $S_{wr} $ and $S_{or} $ are the residual saturations for water and oil, respectively.

In relation to the general model description given in COMSOL notation in equations (\ref{18b)}, \ref{18c)}), only the matrix ${\rm c}$ and vector ${\rm f}$ are modified
\begin{equation} \label{24)} 
{\rm c}\equiv \left(\begin{array}{cc} {k\lambda } & {0} \\ {k\lambda _{w} } & {\varepsilon } \end{array}\right);{\rm \; \; }{\rm f}\equiv \left(\begin{array}{c} {0} \\ {0} \end{array}\right) 
\end{equation} 

Since $p_{cow} \equiv 0$ and $q_{w} ,q_{o} \equiv 0$. Note that it was introduced a small artificial diffusion coefficient ($\varepsilon $) in the saturation equation to stabilize the numerical solution, due to its hyperbolic nature, numerical instabilities can be appeared.

\begin{table}
\begin{tabular}{p{2.0in} p{0.5in} p{0.7in}} \hline 
Property & Units & Value \\ \hline 
Domain length (\textit{$L$}) & m & 300 \\ 
Absolute Permeability ($k$) & m${}^{2}$ & 1.00E-15 \\ 
Porosity ($\phi $) &  & 0.2 \\ 
Water viscosity ($\mu _{w} $) & Pa.s & 1.00E-03 \\ 
Oil viscosity ($\mu _{o} $) & Pa.s & 1.00E-03 \\ 
Residual water saturation ($S_{wr} $) &  & 0 \\ 
Residual oil saturation ($S_{or} $) &  & 0.2 \\ 
Injection velocity ($u_{w}^{in}$) & m.s${}^{-1}$ & 3.4722E-07 \\ 
Production pressure ($p^{out} $) & MPa & 10 \\ 
Artificial diffusion coefficient ($\varepsilon $) & ~ & 1.00E-7 \\\hline   
\end{tabular}
\caption[Buckley-Leverett problem data. The viscosity of oil and water were taken initially as equal for a basis case, and then played with their proportion, as in\cite{Hoteit2008}, (see Table \ref{TableBL2})]{Buckley-Leverett problem data. The viscosity of oil and water were taken initially as equal for a basis case, and then played with their proportion, as in\cite{Hoteit2008}, (see Table \ref{TableBL2})}\label{TableBL1}
\end{table}

For this problem the initial conditions are
\begin{equation} \label{25)} 
p\left(t_{0} \right)=p_{0} \equiv 10{\rm MPa},{\rm \; \; }S_{w} \left(t_{0} \right)=S_{w0} \equiv 0 
\end{equation} 
and the corresponding boundary conditions are
\begin{equation} \label{26)} 
{\rm r}^{in} \equiv \left(\begin{array}{c} {0} \\ {S_{w}^{in} } \end{array}\right),{\rm g}^{in} \equiv \left(\begin{array}{c} {u_{w}^{in} } \\ {0} \end{array}\right){\rm ,}{\rm h}^{in} \equiv \left(\begin{array}{cc} {0} & {0} \\ {0} & {1} \end{array}\right) 
\end{equation} 
\begin{equation} \label{27)} 
{\rm r}^{out} \equiv \left(\begin{array}{c} {p^{out} } \\ {0} \end{array}\right),{\rm g}^{out} \equiv \left(\begin{array}{c} {0} \\ {0} \end{array}\right){\rm ,}{\rm h}^{out} \equiv \left(\begin{array}{cc} {1} & {0} \\ {0} & {0} \end{array}\right) 
\end{equation} 
Where $S_{w}^{in} \equiv 0.8$, $g_{p}^{in} \equiv 3.47{\rm E-}7{\rm m}\cdot {\rm s}^{{\rm -1}} $, $p^{out} \equiv 10{\rm MPa}$.

The relevant data are taken from \cite{Hoteit2008} and are provided in Table \ref{TableBL1}.

The simulations were carried out for three cases with different water-oil viscosity ratios combining two types of relative permeability models (linear and quadratic) for seven time periods, (see Table \ref{TableBL2}). 

\begin{table}
\begin{tabular}{p{0.7in}  p{0.8in}  p{0.7in}  p{0.9in}} \hline 
Cases & Relative \newline permeability\newline model ($\omega $) & Viscosity \newline ratio ($\mu _{w}/\mu _{o}$) & Simulation\newline periods\newline ($t_{\max } $) [days] \\ \hline
(a) & 1 & 2 & 300-900 \\ 
\textbf{}(b) & 1 & 2/3 & 300-900 \\ 
\textbf{}(c) & 2 & 2/3 & 300-900 \\ \hline 
\end{tabular}
\caption[Simulated cases for Buckley-Leverett problem]{Simulated cases for Buckley-Leverett problem}\label{TableBL2}
\end{table}
\subsection{Water Flooding Case Study}

The second problem is about to reproduce the flow behavior in a water flooding experiment through a sandstone core under laboratory conditions. The intention is to couple this flow model with multicomponent transport equations to study Enhaced Oil Recovery processes \cite{lopez2008transport}. Data of this problem is given in Table \ref{TableBL3}.

In this case, the relative permeability constitutive equations are based on the Brooks-Corey  model \cite{Brooks1964}:
\begin{equation} \label{28)} 
k_{rw} =S_{e}^{\frac{2+3\theta }{\theta } } {\rm ;\; \; \; }k_{ro} =\left(1-S_{e} \right)^{2} \left(1-S_{e}^{\frac{2+\theta }{\theta } } \right); 
\end{equation} 
where $\theta $ characterizes the pore size distribution.

While oil-water capillary pressure is defined by the Brooks-Corey \cite{Brooks1964} model given in Eq. (\ref{Eq_A3}).

%

We impose the following initial conditions
\begin{equation} \label{31)} 
p\left(t_{0} \right)=p_{0} \equiv 10{\rm MPa},{\rm \; \; }S_{w} \left(t_{0} \right)=S_{w0} \equiv 0.2 
\end{equation} 
and boundary conditions
\begin{equation} \label{32)} 
{\rm r}^{in} \equiv \left(\begin{array}{c} {0} \\ {S_{w}^{in}} \end{array}\right),{\rm g}^{in} \equiv \left(\begin{array}{c} {u_{w}^{in} } \\ {0} \end{array}\right){\rm ,}{\rm h}^{in} \equiv \left(\begin{array}{cc} {0} & {0} \\ {0} & {1} \end{array}\right) 
\end{equation} 
\begin{equation} \label{33)} 
{\rm r}^{out} \equiv \left(\begin{array}{c} {p^{out} } \\ {0} \end{array}\right),{\rm g}^{out} \equiv \left(\begin{array}{c} {0} \\ {0} \end{array}\right){\rm ,}{\rm h}^{out} \equiv \left(\begin{array}{cc} {1} & {0} \\ {0} & {0} \end{array}\right) 
\end{equation} 
Where $g_{p}^{in} \equiv {\rm 5}.{\rm 3E}-0{\rm 7\; }{\rm m}\cdot {\rm s}^{{\rm -1}} $, $p^{out} \equiv 10{\rm MPa}$.

\begin{table}
\begin{tabular}{p{2.0in} p{0.4in} p{0.6in} } \hline 
Property & Units & Value \\ \hline 
Domain length ($L$) & m & 0.25 \\ 
Absolute permeability ($k$) & m${}^{2}$ & 8.25E-13 \\  
Porosity ($\phi $) &  & 0.2 \\  
Water viscosity ($\mu _{w} $) & Pa.s & 1.00E-03 \\  
Oil viscosity ($\mu _{o} $) & Pa.s & 1.00E-02 \\ 
Residual water saturation ($S_{wr} $) &  & 0.2 \\ 
Residual oil saturation ($S_{or} $) &  & 0.15 \\ 
Injection velocity ($u_{w}^{in} $) & m.s${}^{-1}$ & 5.3E-07 \\ 
Production pressure ($p^{out} $) & MPa & 10 \\  
Brooks-Corey parameter ($\theta $)  &  & 2 \\ 
Entry threshold pressure ($p_{t} $) & ~MPa & 1.00E-2 \\ \hline 
\end{tabular}
\caption{Water coreflooding experimental data}\label{TableBL3}
\end{table}

\begin{figure}[h!]
\centering
\captionsetup{width=0.43\textwidth}
\begin{minipage}{.5\textwidth}
  \centering
  \includegraphics[width=.9\linewidth]{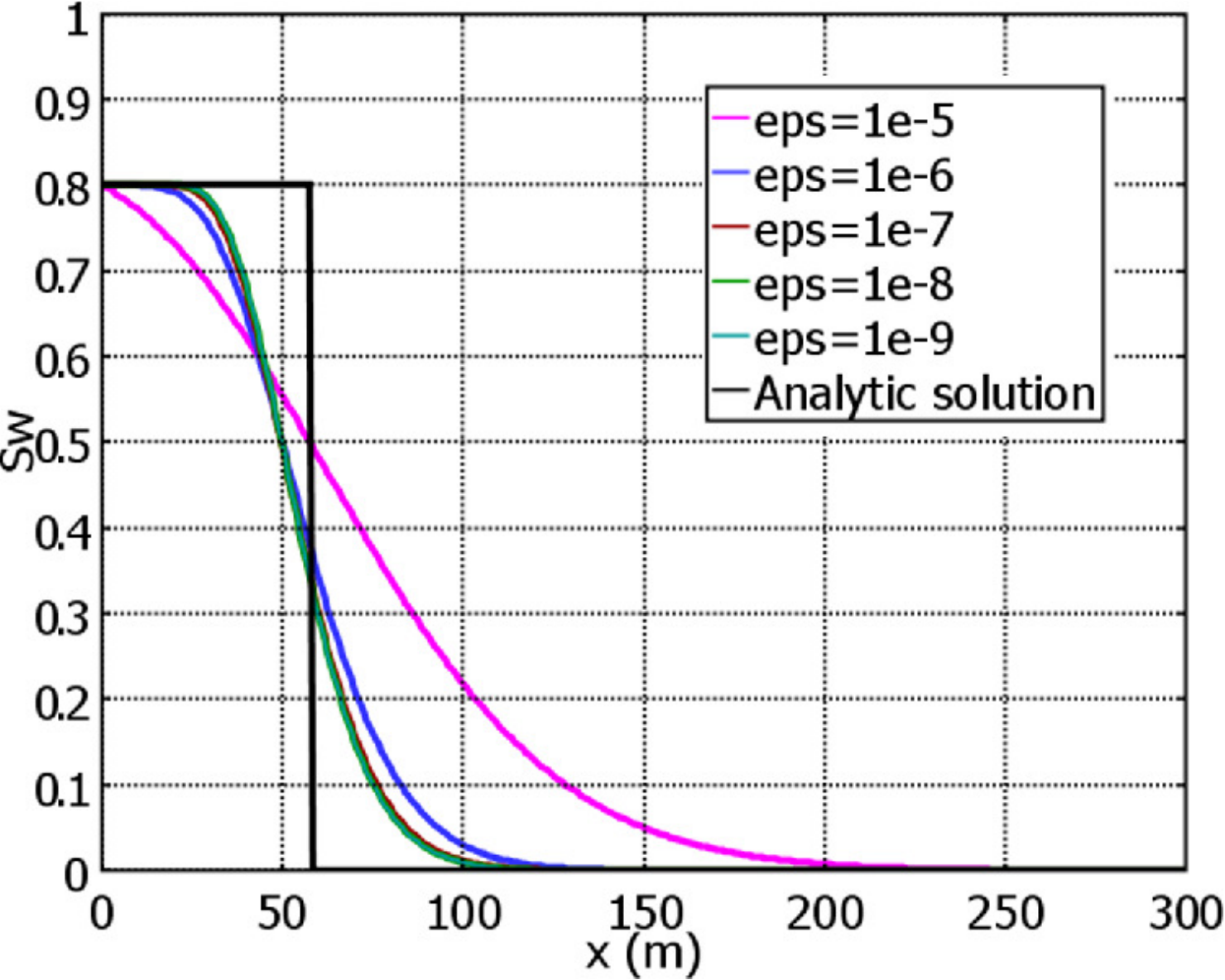}
  \captionof{figure}{Numerical solutions of the Buckley--Leverett problem with linear relative permeabilities and viscosity ratio $\mu _{w}/\mu _{o} =1$ for a period of 300 days, varying artificial diffusion coefficient ($\varepsilon $).}
  \label{FigBL1}
\end{minipage}%
\begin{minipage}{.5\textwidth}
  \centering
  \includegraphics[width=.9\linewidth]{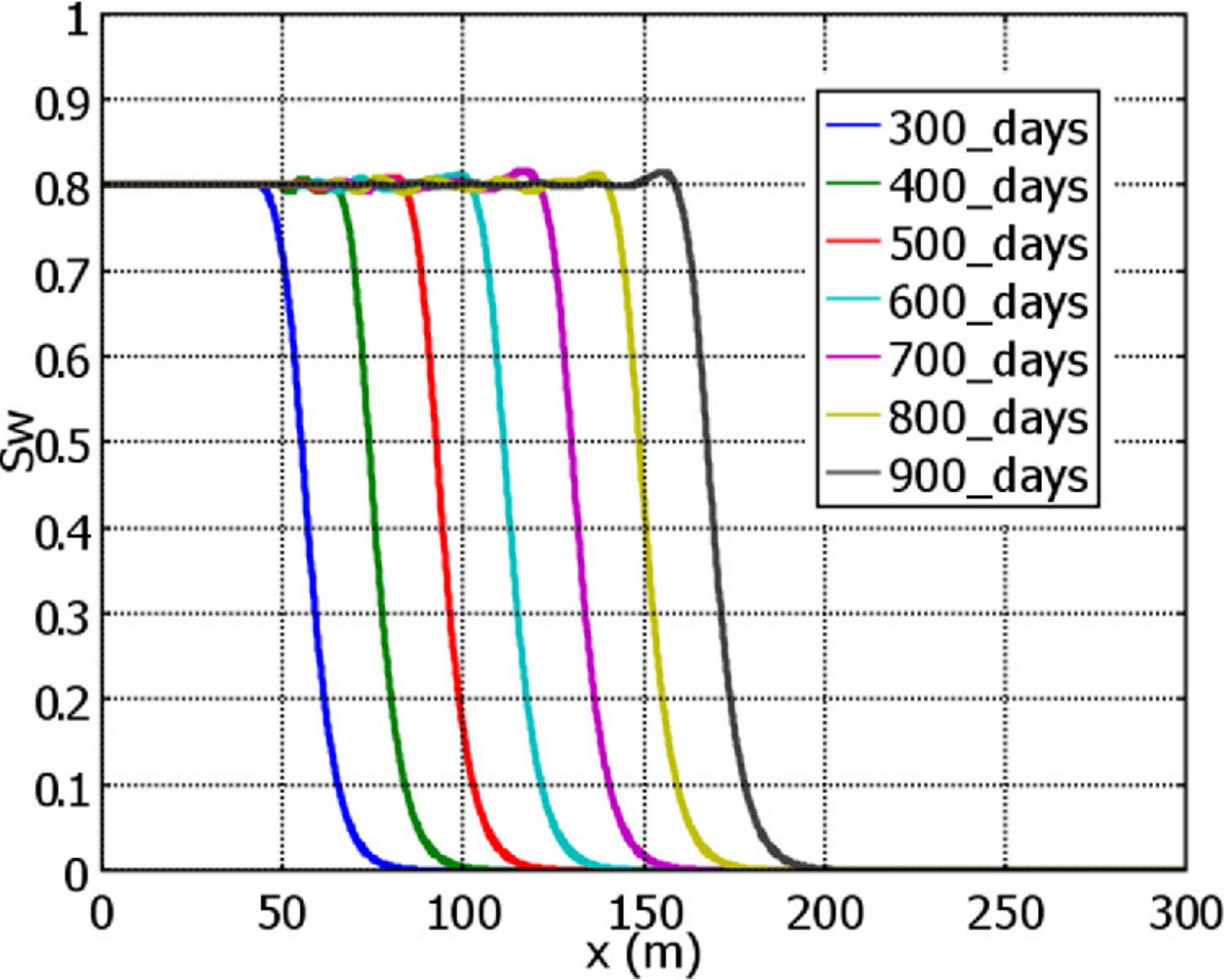}
  \captionof{figure}{Numerical solutions of the Buckley--Leverett problem with linear relative permeabilities for case (a) and viscosity ratio $\mu _{w}/\mu _{o} =2$ for time periods from 300 to 900 days.}
  \label{FigBL2}
\end{minipage}
\end{figure}

In the Fig. \ref{FigBL1} the numerical solutions of the Buckley--Leverett problem with linear relative permeabilities and viscosity ratio $\mu _{w}/\mu _{o}=1$ for a time period of 300 days with different artificial diffusion coefficients are shown.  We can observe that it is attained the best trade of in terms of efficiency and accuracy for an artificial diffusion coefficient value $\varepsilon =1e-7$.

\begin{figure}[htb]
\centering
\captionsetup{width=0.43\textwidth}
\begin{minipage}{.5\textwidth}
  \centering
  \includegraphics[width=.9\linewidth]{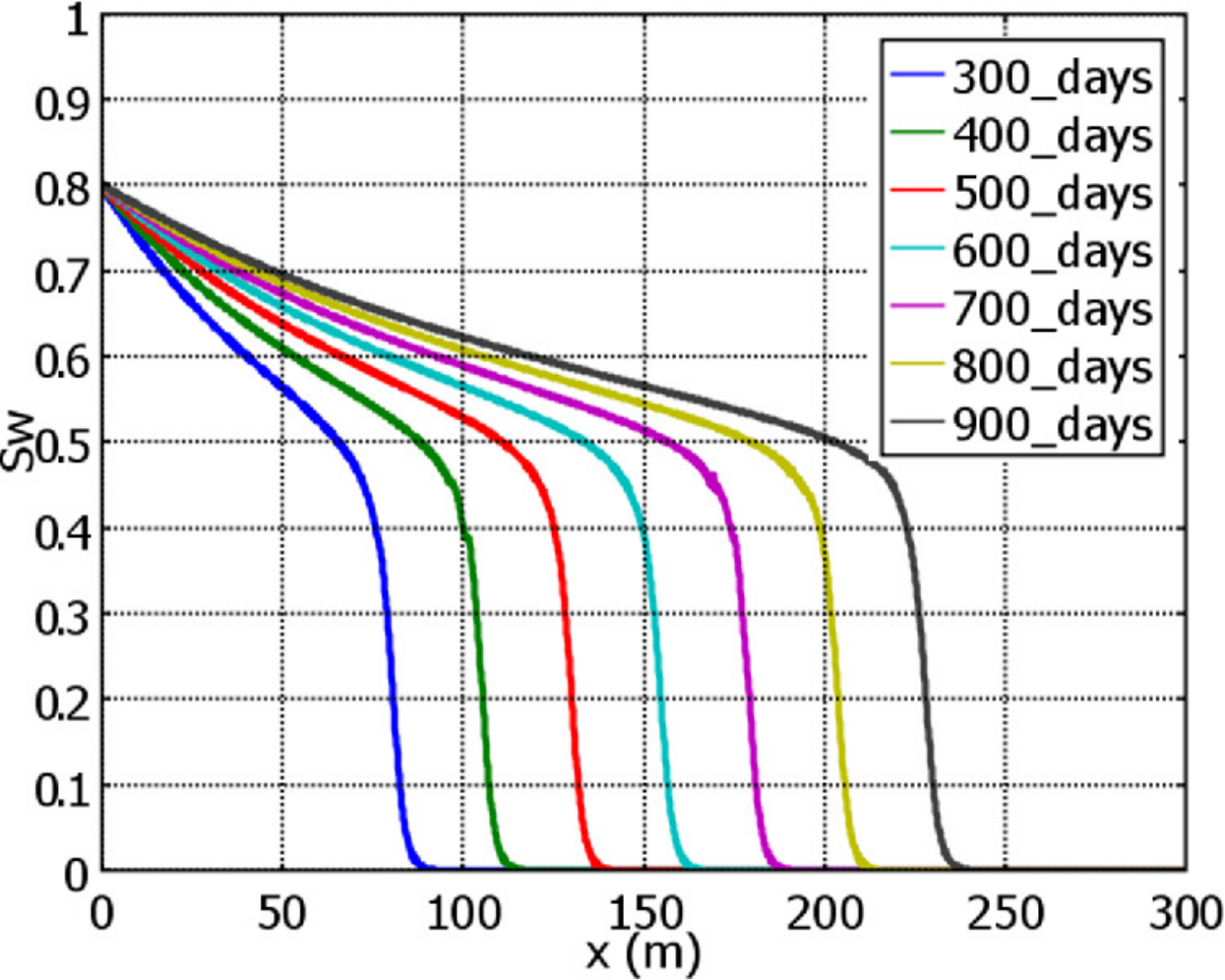}
  \captionof{figure}{Numerical solutions of the Buckley--Leverett problem case (b) with linear relative permeabilities and viscosity ratio $\mu _{w}/\mu _{o} =2/3$ for a period of 300 days, varying artificial diffusion coefficient ($\varepsilon $).}
  \label{FigBL3}
\end{minipage}%
\begin{minipage}{.5\textwidth}
  \centering
  \includegraphics[width=.9\linewidth]{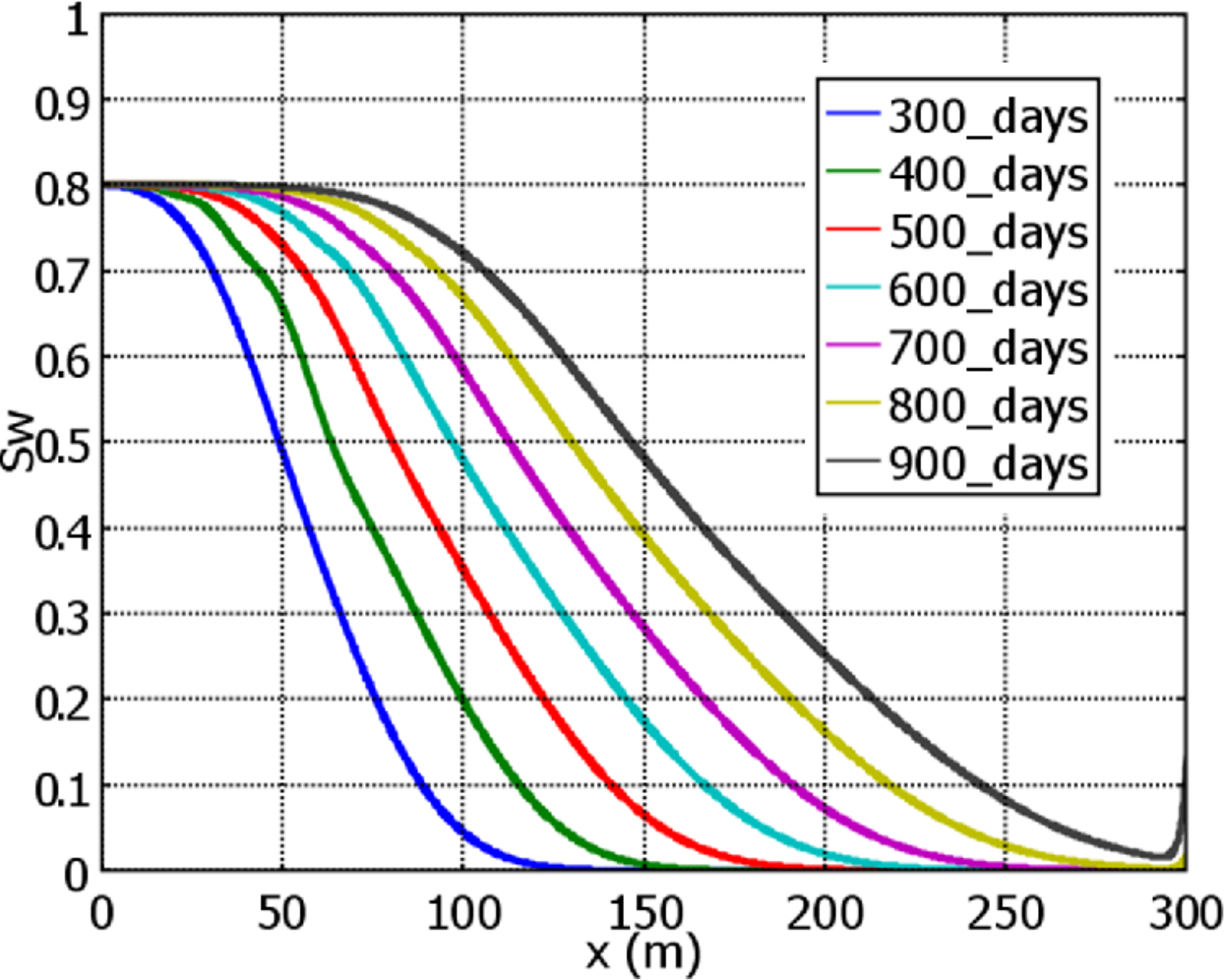}
  \captionof{figure}{Numerical solutions of the Buckley--Leverett problem, case (c) with quadratic relative permeabilities and viscosity ratio $\mu _{w}/\mu _{o} =2/3$ for time periods from 300 to 900 days.}
  \label{FigBL4}
\end{minipage}
\end{figure}

Figs. \ref{FigBL2}-\ref{FigBL4} show a quite well qualitative reproduction of the analytic solution behavior for the Buckley--Leverett problems for cases (a)-(c) described in table \ref{TableBL2}, respectively. These problems were numerically solved with the optimal artificial diffusion coefficient value $\varepsilon =1e-7$ previously obtained.

\begin{figure}[h!]
    \centering
    \includegraphics[width=0.4\textwidth]{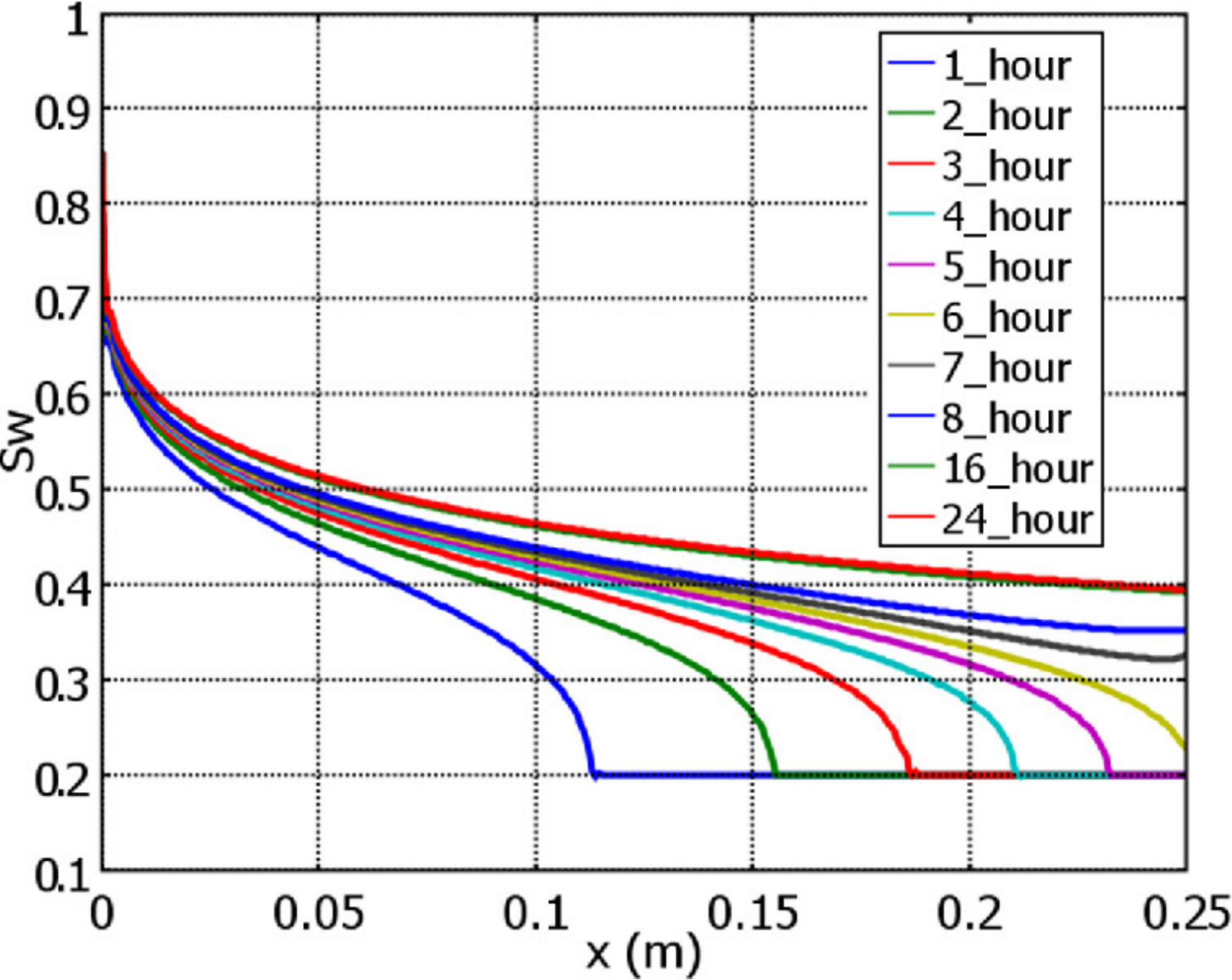}
  \caption[Numerical simulation of the water coreflooding experiment for a time period of 24 hours.]{Numerical simulation of the water coreflooding experiment for a time period of 24 hours.}
\label{FigBL5}           
\end{figure}

The numerical simulation of the water coreflooding experiment through a sandstone core during a time period of 24 hours is shown in Fig. \ref{FigBL5}. It can be observed the formation of a water front displacing the oil through the porous medium which is recovered at the production end of the core. 

The main result of the present work is the implementation of a biphasic (water-oil) flow model in porous media, including capillary pressure, which coupled to multiphase and multicomponent transport equations could be useful to study Enhaced Oil Recovery processes at laboratory scale.

Even more, applying a flow model coupled with transport equations can serve to study the impact in the flow conditions due to the porosity and permeability alterations by transport processes, such as adsorption of some fluent components.



%
%
\clearpage

\section{Validation of the clogging/declogging modeling}
\label{validation}
The validation of the flow model has been already presented  in several previously published works and for further details the reader is referred to them \cite{Diaz2008a, Diaz2012, Diaz2012b, Diaz2015,Diaz2015b}. The transport model equations (Eqs. \ref{Eq_D4b},\ref{Eq_D5b},\ref{Eq_D5c}) were tested for validation of clogging/declogging processes against the experimental data published by \cite{Hendry1997} and compared with the results obtained in \cite{Li2011}  and \cite{Kim2006}. Firstly, our model parameters were fitted  with the experimental Hendry's data and  then were compared with the fitting parameters obtained by \cite{Li2011}  and \cite{Kim2006}. 

\subsection[Experiment Description]{Experiment Description}
The experimental setup  consists of the water injection with microorganisms from the bottom into a vertical column packed with silica sand \cite{Hendry1997}. Where column length is $L=40$ cm and the diameter is $d=5$ cm, respectively. The initial column porosity $\phi_0$ is 0.4, and the dispersivities are ${\alpha_{L}}_{w}^{m}= {\alpha_{T}}_{w}^{m} = 0.27$ cm. The microorganisms, having a density $\rho _ m = 1.085\times 10^{6}$ $\mathrm{mg/L}$, are injected with a constant concentration of $c^{m_{in} }_w = 4.32$ $\mathrm{mg/L}$ at a constant water velocity of $u_w = 2.17\times 10^{-4}$ $\mathrm{cm/s}$. 

\subsection[Modeling Considerations]{Modeling Considerations}

For comparison purposes \cite{Kim2006} and \cite{Li2011} works are considered. In all cases
molecular diffusion is omitted, i.e., $D^{*m}_{\;\;w} =0$ and bacterial diffusion processes are set as a single parameter of dispersivity. The bacterial settling velocity is included by \cite{Kim2006} and \cite{Li2011}, but in this work this was neglected since it is too small, about two orders of magnitude with respect to water phase velocity $u_w$.

The biological kinetic parameters controlling growth/death and clogging/ declogging processes such as maximal growth rate ($g_{m}^{\max } $), decaying rate ($d_{m}$), reversible clogging rate ($k_{c1}$), irreversible clogging ($k_{c2}$), declogging rate ($k_{d}$),  given in Table \ref{Table_6_1}, are considered for fitting the experimental results of Fig. \ref{Figura_6_1}.

The combination of the described adjustments in the values of the rates  of growth, decaying, reversible/irreversible clogging and declogging, renders a total of three possible configurations, those given in the Table \ref{Table_6_1} (with their corresponding resultant RMS), which are depicted in the Fig. \ref{Figura_6_1}.
\clearpage
\begin{table}
\centering
  \begin{tabular}{ c  c  c  c  c  }
\hline  
  Property & \cite{Kim2006} &\cite{Li2011} &  This work \\
    \hline
    $g_{m}^{max}$ &  $1\times 10^{-6} $ & 0   &  $1\times 10^{-6} $\\
    $d_m$ & $1\times 10^{-7}$& 0 &  $1\times 10^{-7}$ \\
    $k_{c1}$ & $2.28 \times 10^{-5} $  & $2.28 \times 10^{-5} $ &$5.7 \times 10^{-5} $ \\
    $k_d$ & $3.56\times 10^{-7} $  & $3.56\times 10^{-7} $ &  $6.408\times 10^{-7} $\\
    $k_{c2}$ & $1.72\times 10^{-6} $  & $1.72\times 10^{-6} $ & $3.44\times 10^{-6} $ \\
    RMS & $1.9\times 10^{-2}$ & $2.1\times 10^{-2}$ & $2.1\times 10^{-3}$\\
    \hline
   \end{tabular}
       \caption[Fitting parameters for the experimental results of \cite{Hendry1997}.]{Fitting parameters for the experimental results of \cite{Hendry1997}. All have units as [$\mathrm{s^{-1}}$], except the resultant measure of error, RMS, which is adimensional (proportion of effluent/injected microorganisms).}\label{Table_6_1}
\end{table}

\subsection[Results and Discussion]{Results and Discussion}
A one-dimensional transport problem of clogging and declogging was solved in the test of the model in \cite{Kim2006} and \cite{Li2011}, using for comparison the experimental data in \cite{Hendry1997}. In the model developed for this work, there is a complete 3D system, and the clogging/declogging equations were considered as two separate transport equations. So, the model derived for this work is general enough that it can take as particular cases those models of \cite{Kim2006} and \cite{Li2011}.

\begin{figure}[htb]
    \centering
    \includegraphics[width=0.9\textwidth]{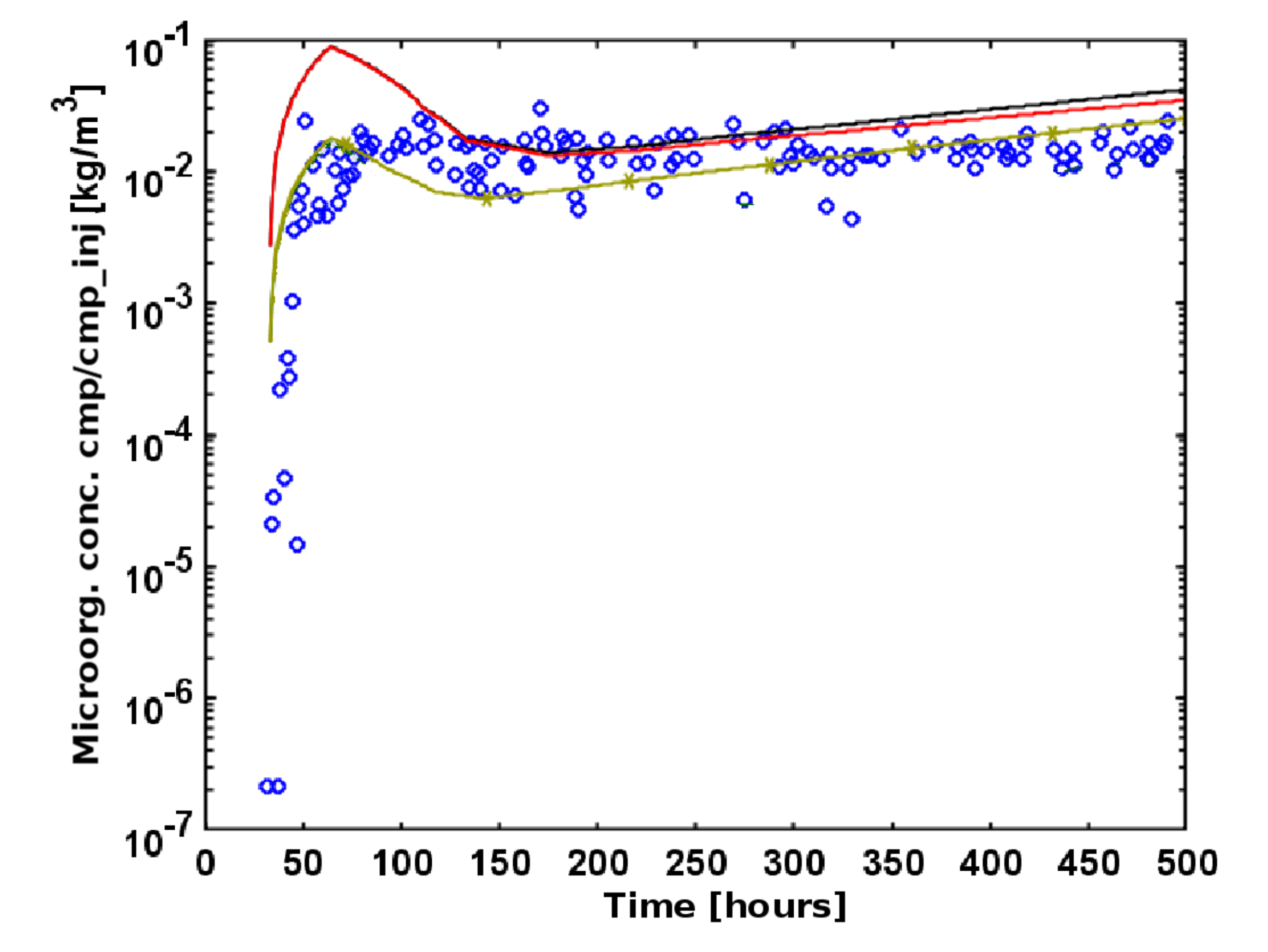}
  \caption[Experimental data and numerical simulations of breakthrough planktonic microorganisms concentration with respect to time. The blue circles are the \cite{Hendry1997} experimental values; the red line is the \cite{Kim2006} fitting; the black line is the \cite{Li2011} fitting; the light green with asterisks is our own fitting with $k_d=6.408\times 10^{-7}$.]{Experimental data and numerical simulations of breakthrough planktonic microorganisms concentration with respect to time. The blue circles are the \cite{Hendry1997} experimental values; the red line is the \cite{Kim2006} fitting, with $\text{RMS}= 1.9\times 10^{-2}
  $; the black line is the \cite{Li2011} fitting, with $\text{RMS}= 2.1\times 10^{-2}$; the light green with asterisks is our own fitting with $k_d=6.408\times 10^{-7}$, and $\mathrm{RMS}=2.1\times 10^{-3}$.}
\label{Figura_6_1}           
\end{figure}

In the Fig. \ref{Figura_6_1} it can be appreciated the comparison of the numerical simulation with the different values for the experimental parameters. Our model does cross among the data from \cite{Hendry1997}, in a qualitative manner similar to that of \cite{Li2011} and \cite{Kim2006}.

 Tests were performed, first with the suggested parameters by \cite{Li2011}, but as the fitting of the experimental data with \cite{Hendry1997} did not quite agree with the results from \cite{Li2011}, the parameters were changed, specially those for clogging and declogging, until  more adequate results were obtained. 

Concerning Fig. \ref{Figura_6_1}, the following observations can be made:
\newline
The difference between the black line and the red one is the fact that the red one has values different of zero for the growth and decaying rates ($g_{m}^{max}$ and $d_m$, respectively), but that does not seem to make a significant qualitative or quantitative difference (the fitting curves overlapped), with respect to that instance which does not consider either rate, so for our own adjustments it was decided to simply leave those processes active (with the values from \cite{Kim2006}), and only adjust the  reversible/irreversible clogging rates and the declogging rate. In a first attempt, it was only considered the adjustment of the clogging rates, and the result is the dark green line, which fits better the group of experimental values near the concentration of $10^{-2}$. Then it was manipulated the value of the declogging rate (light green line with asterisks), and it was obtained the fit which seems to embrace the most of the experimental values, compared with the other adjustments, including those of \cite{Kim2006} or those of \cite{Li2011}. From here, and from the resultant RMS, it can be concluded that the computational model has been validated,  since the fit is approximately very close, qualitatively and quantitatively speaking, to the experimental values of \cite{Hendry1997}. 

\section{Case study}
\label{caseDP1}
 
\subsection[Description of the Experiment]{Description of the Experiment}

 A laboratory experiment of microbial recovery of hydrocarbons was performed  where it was used as porous media a sandstone Berea core, while the oil and the microbial culture employed came from the field Agua Fr\'ia. Hereafter this experiment will be called DP1, which has been previously reported in the work of \cite{Castorena2012a} and a similar set up \cite{Castorena2012c}. The microbial culture used in the  recovery experiment  is comprised of microorganisms thermophiles, barophiles, halotolerants, acidotolerants, and anaerobes; containing the following microorganisms among others: \emph{Thermoanaerobacter ethanolicus}, \emph{Thermoanaerobacter uzonensis}, \emph{Thermoanaerobacter inferii}, \emph{Geothermobacterium sp}, \emph{Methanobacterium subterraneum}, \emph{Methanobacterium formicicum}, \emph{Methanolinea tarda} and /or \emph{Methanoculleus sp} \cite{Olguin2014}. The  mixed culture  were obtained  from oil samples of the aforementioned Mexican oil field.

\begin{figure}[htb]
    \centering
    \includegraphics*[width=5.17in, height=3.24in, keepaspectratio=false, trim=0.00in 0.00in 0.00in 0.00in]{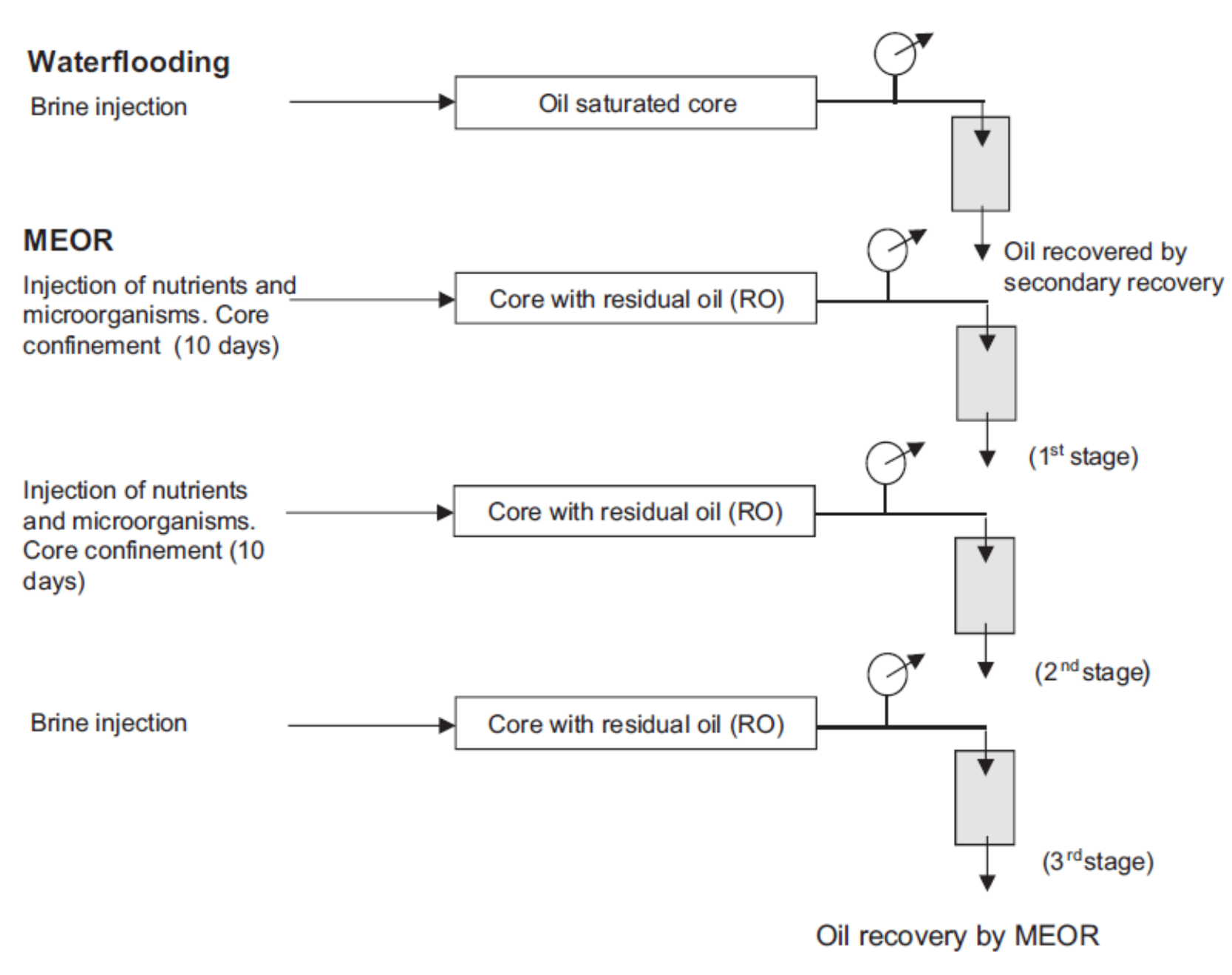}
  \caption[Schematic representation of the experiment DP1.]{Schematic representation of the experiment DP1.}
\label{Fig_4_Experiment}           
\end{figure}

The experiment consisted on initially saturating the Berea core with water and oil, to perform oil displacement through water injection in a first stage of \emph{secondary recovery}, and after that perform three stages of water injection with microorganisms and nutrients that will be named \emph{MEOR} stages, alternating with periods with no internal flow that will be called \emph{confinement} stages. So, first there is a secondary recovery, with oil displacement due to brine injection, until it was obtained the residual oil saturation lasting 68 hours. Then  MEOR1 stage, the first injection of microorganisms and nutrients, lasting 66 hours. Then the Confinement1 stage, where the system is closed and it is maintained under conditions of pressure and temperature previously established, lasting 240 hours. Then MEOR2, the second injection of microorganisms and nutrients, lasting 72 hours, followed by the second confinement stage of 240 hours. Finally, the MEOR3 stage, the third injection of microorganisms and nutrients, lasting 62 hours (see Fig. \ref{Fig_4_Experiment}).

 The results of the process of hydrocarbon recovery in the core showed that there was recovery in each of the stages, being the most significant the secondary recovery, followed by the first stage of microbial recovery  (see Fig. \ref{Figura_7_1_61}). It will be explained in detail the secondary recovery stage, and the first stage of microbial recovery. A summary of all the stages is described in Table \ref{Tabla_7_1_2b}. In Fig. \ref{Figura_7_1_10} it is observed the profile of the pressure drop, which behaves smoothly and continuously. While formally the oil displacement experiment by water lasted around 68 hours,  by observing the graph in Fig. \ref{Figura_7_1_5} one can see that the oil recovery curve keeps a similar behavior until 68 hours, which shows the typical conduct of a recovery process, where the porous medium is strongly wettable by water. It can be observed that the oil recovery reaches a value of approximately 72 ml.

\begin{table}[htp]
\begin{center}
\begin{tabular}{p{0.7in} p{0.7in} p{0.7in} p{0.7in} p{0.8in} p{0.7in}}\hline 
\textbf{\newline Stages of \newline recovery} & \textbf{Oil \newline recovery (experiment) [ml]} & \textbf{Oil \newline recovery (numerical) [ml]} &\textbf{Oil \newline recovery\newline \newline \newline [\%]} &\textbf{Root \newline mean\newline square\newline error\newline [ml]} & \textbf{Duration\newline \newline \newline \newline [hr]} \\ \hline 
Secondary \newline recovery & 72.2 & 73.23 & 49.45 & 2.13  & 68 \\ 
MEOR1 & 33.9 & 33.01 & 23.22 & 2.08  & 66  \\  
MEOR2 & 9.4 & 9.7 & 6.44  & 1.58  & 72 \\ 
MEOR3 & 0.21 & 3.17 & 0.14 & 2.34  & 62 \\ 
\emph{MEOR} & \emph{43.51} & \emph{45.89}  & \emph{29.8} & \emph{2.00}  & \emph{200}  \\ 
\textbf{Total} & \textbf{115.71} &\textbf{123.84} & \textbf{79.25} &\textbf{2.04} & \textbf{268}  \\ \hline 
\end{tabular}
\end{center}
\caption{Hydrocarbon recovery for each stage of the DP1 experiment.}
\label{Tabla_7_1_2b}
\end{table}

 The total microbial recovery of 29.8 \% of oil was close to the achieved values by other authors \cite{Bryant1988} with light oils ($31^{\circ}$ API, 32\% of recovery) and greater that the obtained with heavy oils. It is worth noticing that the oil used in this work has $21^{\circ}$ API, which would classify it as a medium oil \cite{Olguin2014}.

\subsection[Numerical Simulation of the Experiment]{Numerical Simulation of the Experiment}

\subsubsection[Simulation of the secondary recovery stage]{Simulation of the secondary recovery stage}

\begin{figure}[htb]
    \centering
    \includegraphics*[width=5.17in, height=3.24in, keepaspectratio=false, trim=0.00in 0.00in 0.00in 0.00in]{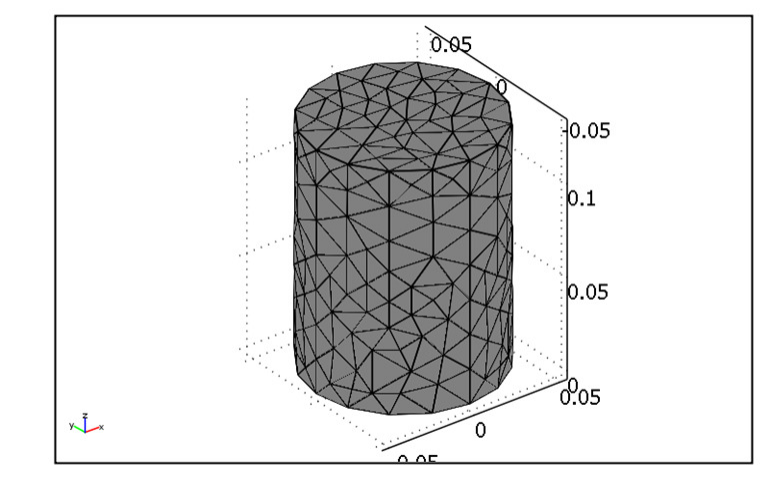}
  \caption[Representation of a Berea core where the computational domain and mesh is observed.]{Representation of a Berea core where the computational domain and mesh is observed.}
\label{Figura_7_1_4}           
\end{figure}

 The computational mesh used for the secondary recovery stage is shown in  Fig. \ref{Figura_7_1_4}. 
 The total number of elements of the mesh is 2,358 for a total of 7,464 degrees of freedom (for the case when only the flow model is active).

\clearpage

 The input data for the computational model corresponding to the secondary recovery stage are listed in the Table \ref{Tabla_7_1_3}.

\begin{center}
\begin{longtable}{p{0.8in} p{0.7in} p{0.8in}} \hline 
\textbf{Property} &  \textbf{Value} & \textbf{Units} \\ \hline 
 $L$ & 0.13 & [m] \\ 
$d$& 0.1016 & [m] \\ 
$k_0$  & 1.51E-13 & [$\mathrm{m^2}$] \\ 
$\phi_{0}$ & 0.1978 & [$\mathrm{m^3/m^3}$] \\ 
$V_p$ & 2.08E-04 & [m$^3$] \\  
$\mu _{w}$ & 4.00E-04 &  [$\mathrm{Pa.s}$]\\  
$\mu _{o}$ & 1.31E-02 & [$\mathrm{Pa.s}$] \\  
$\rho_{w}$ & 9.89E+02 & [kg/m$^3$] \\  
$\rho_{o}$ & 8.72E+02 & [kg/m$^3$] \\  
$Q$ & 1.39E-09 & [$\mathrm{m^3/s}$] \\  
$u_{w}^{in}$ & 1.71E-07 & [m.s${}^{-1}$] \\  
$p^{out}$ & 551,580 &  [Pa] \\  
$p_{o}^{0}$ & 551,580 & [$\mathrm{Pa}$] \\  
$S_{w0}$ & 0.299 & [$\mathrm{m^3/m^3}$] \\  
$S_{o0}$ & 0.701 & [$\mathrm{m^3/m^3}$] \\  
$t$ & 244,800 & [s] \\ \hline 
\caption{Data for the secondary recovery stage of the recovery experiment  DP1. Further details can be found in the section for Nomenclature and Units.}\label{Tabla_7_1_3}
\end{longtable}
\end{center}

 For the fitting of the production of oil and water for the secondary recovery stage, it was designed the following methodology:

\begin{itemize}

\item The relative permeabilities of oil and water were obtained from the experimental values of the capillary pressure as in \cite{Rodriguez1973}.
\item Once obtained the values of relative permeability using the procedures of \cite{Rodriguez1973} above mentioned, those values were fit to the modified Brooks-Corey model.

\item  Obtention of the fit with the inflection point of the corresponding curve to the change from piston flow to simultaneous flow of oil and water. This depends on the residual oil and water saturation, and the  ``endpoints'' $k_{ro}^{0}$ and $k_{rw}^{0} $ of the formula of the modified Brooks-Corey model. By changing the residual saturations it is necessary  to find again the parameters $p_{t}$ and  $\theta $ which are adjusted to the capillary pressure curve in the Brooks-Corey formula.

\item  Once achieved the inflection point, the adjustment of the point of inflection requires the adjustment of the parameters $k_{rw}^{0}$, $n_{w} $, $k_{ro}^{0} $ and $n_{o} $.

\end{itemize}

\begin{figure}[htb]
    \centering
    \includegraphics*[width=12cm,keepaspectratio]{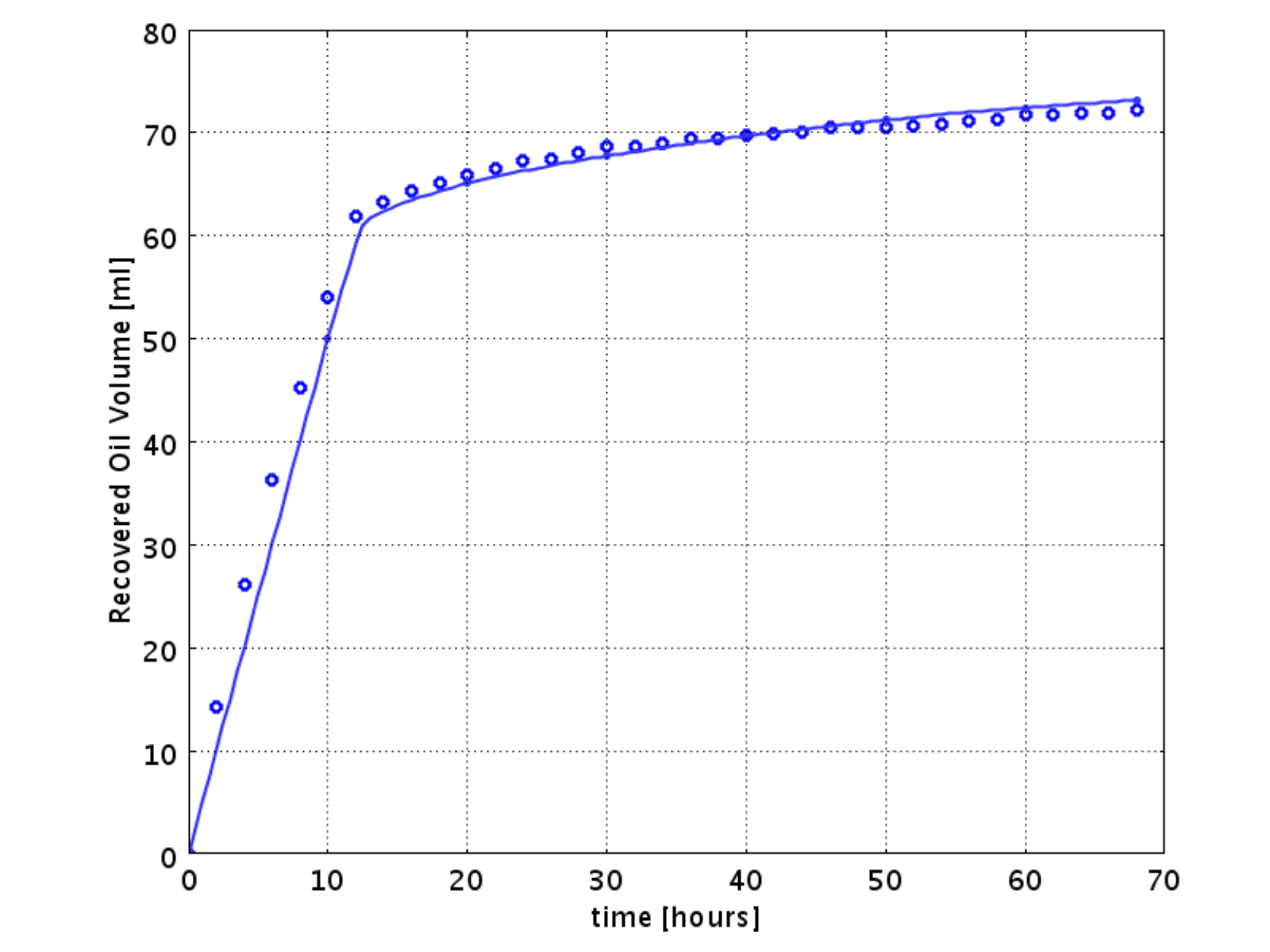}
  \caption[Comparison of the experimental curve (circles) and the simulated one (continuous line) of the recovered oil volume for the secondary recovery stage.]{Comparison of the experimental curve (circles) and the simulated one (continuous line) of the recovered oil volume for the secondary recovery stage.}
\label{Figura_7_1_5}           
\end{figure}

The secondary recovery stage lasted 68 hours with an oil recovery of 72.20 ml, corresponding to 49.52\% of the initial volume. Fig. \ref{Figura_7_1_5} shows a typical recovery curve in a porous medium strongly wettable by water. In this figure a comparison of the experimental recovery values and the simulated ones for the secondary recovery stage is displayed. One can observe an adequate fitting (with an error of less than 10\%) of the recovery curve. In this stage the biphasic flow model was used employing the fitting parameters of capillary pressure and relative permeability, given in Table \ref{Tabla_7_1_4}, which best reproduced the production curves of oil and water. It can be seen that during the first 12 hours the recovery curve behaves as a straight line, which means that it is recovered only oil, followed by an asymptotic behavior of oil and water recovery.

\begin{table}[htb]
\begin{center}
\begin{tabular}[m]{c c c c} \hline 
\textbf{Parameter} & \textbf{Value Secondary Recovery} & \textbf{Value MEOR1}\\ \hline 
$\theta $ & 2.1 & - \\
$p_{t} $ & 2650 [Pa]& - \\
$S_{or} $ & 0.2  & - \\ 
$k_{ro}^{0} $ & 0.12 & - \\
$n_{o} $ & 3.1 & - \\ 
$S_{wr} $ & 0.299 & 0.299\\ 
$k_{rw}^{0} $ & 0.05 & 0.05 \\
$n_{w} $ & 7  & 7 \\  
$S_{or}^{low}$ & -  & 0.2 \\
$k_{ro}^{0\, low}$ & - & 0.12 \\
$n_{o}^{low} $ & - & 3.1 \\
$S_{or}^{high}$ &  - & 0.05 \\
$k_{ro}^{0\, high}$ & - & 0.18 \\
$n_{o}^{high} $ & - & 1.7 \\
\hline
\end{tabular}
\end{center}
\caption{Parameters for the fitting of the capillary pressure and relative permeability curves.}
\label{Tabla_7_1_4}
\end{table}

\begin{figure}[h!]
\centering
\captionsetup{width=0.43\textwidth}
\begin{minipage}{.5\textwidth}
  \centering
  \includegraphics[width=.9\linewidth]{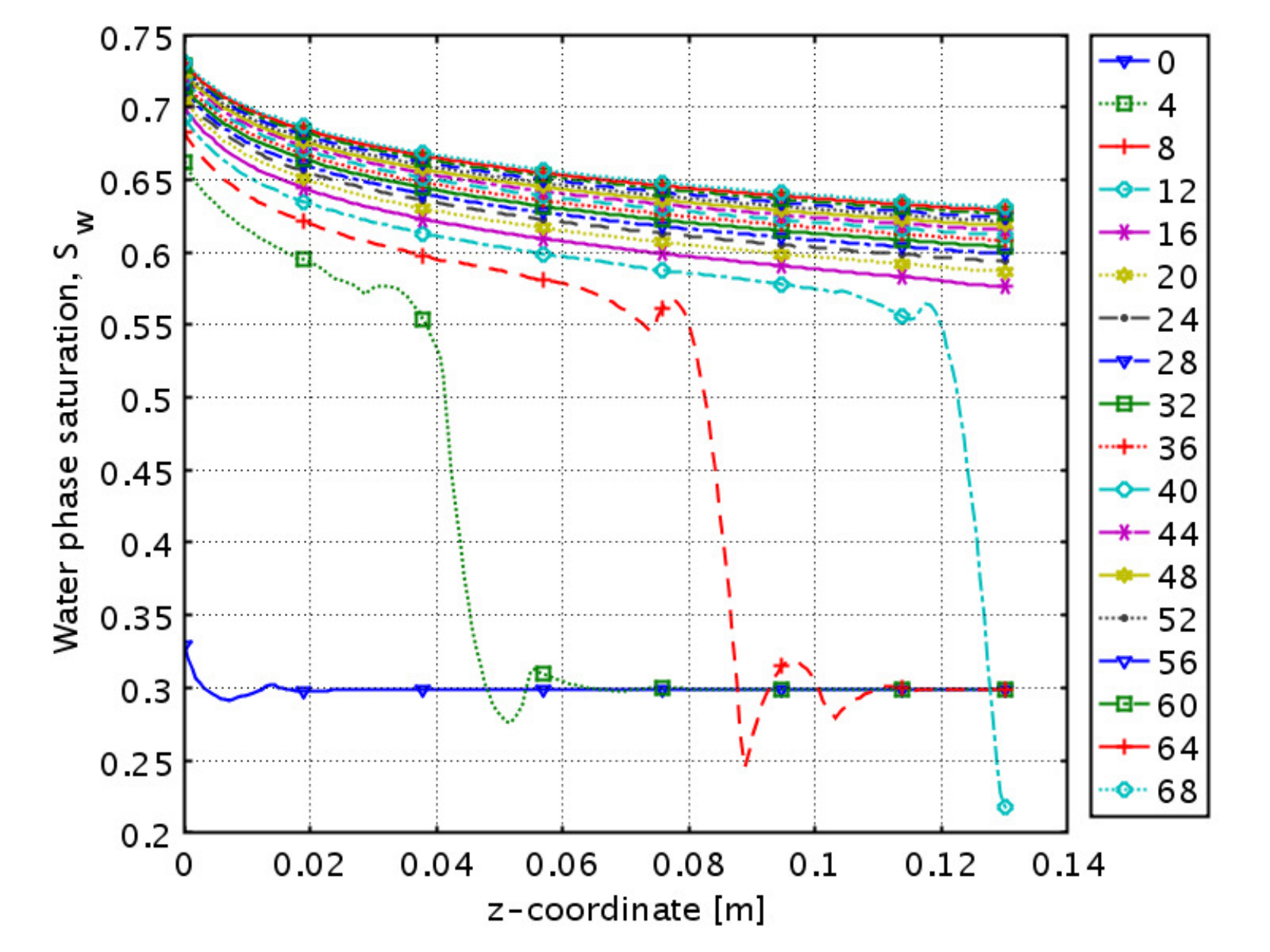}
  \captionof{figure}{Evolution of the water saturation ($S_w$) during the secondary recovery stage for a period of 68 hours. The right column is in hours.}\label{Figura_7_1_7}
\end{minipage}%
\begin{minipage}{.5\textwidth}
  \centering
  \includegraphics[width=.9\linewidth]{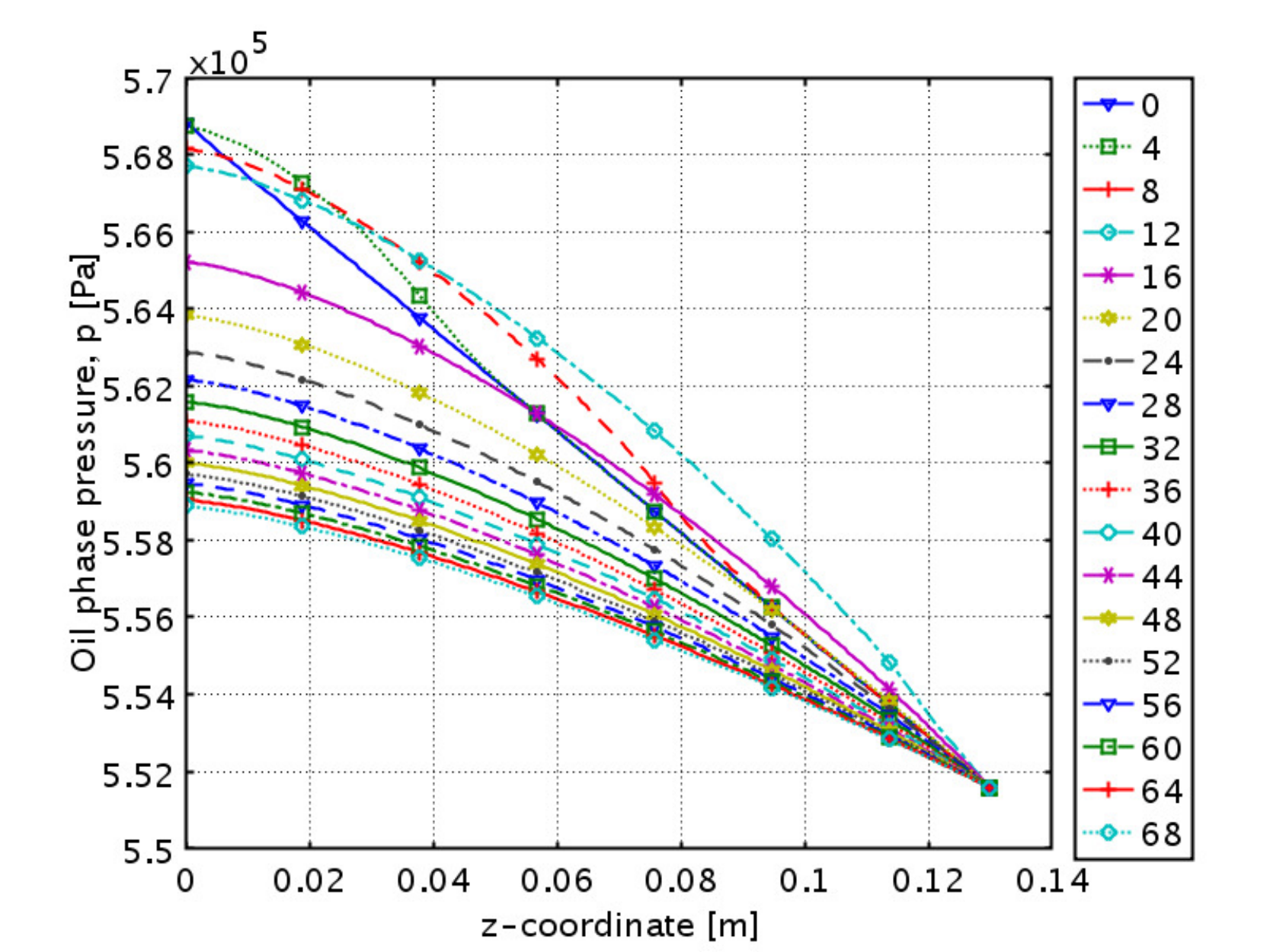}
  \captionof{figure}{Evolution of the oil pressure ($p_o$) during the secondary recovery stage for a period of 68 hours. The right column is in hours.}\label{Figura_7_1_8}
\end{minipage}
\end{figure}

 In Fig. \ref{Figura_7_1_7} it is shown the evolution of the water saturation profile along the core during the water flooding. One can observe the formation of a water front through the porous medium displacing the oil, which in turn is being recovered at the other end of the core.  The water front breaks through the upper extreme in around twelve hours, which matches the breaking point of the recovery curve in Fig. \ref{Figura_7_1_5}. From this point the flow develops an oil saturation close to the residual one.

\begin{figure}[h!]
\centering
\captionsetup{width=0.43\textwidth}
\begin{minipage}{.5\textwidth}
  \centering
  \includegraphics[width=.9\linewidth]{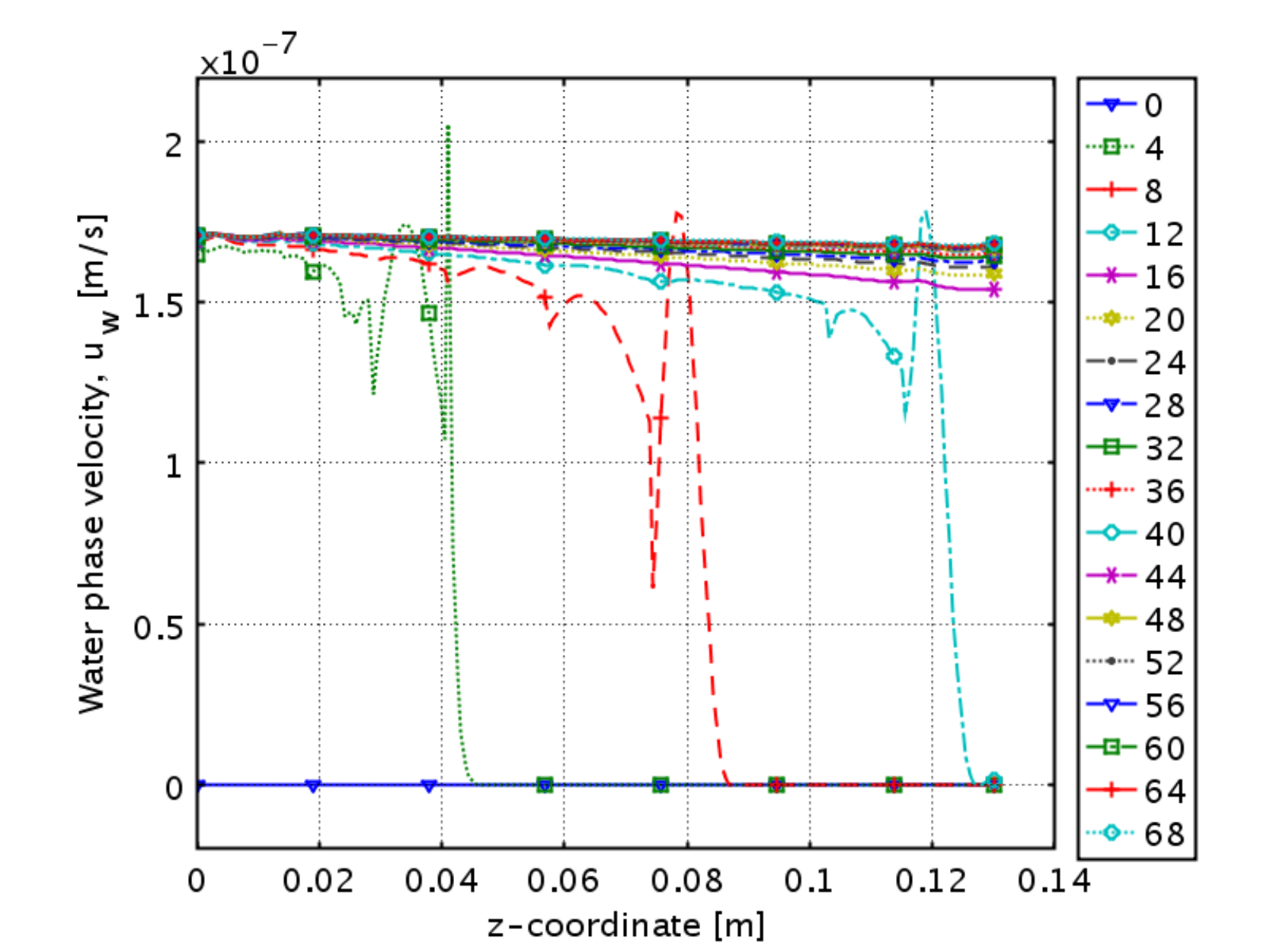}
  \captionof{figure}{Evolution of the water velocity ($u_w$) during the secondary recovery stage for a period of time of 68 hours. The right column is in hours.}\label{Figura_7_1_9}
\end{minipage}%
\begin{minipage}{.5\textwidth}
  \centering
  \includegraphics[width=.9\linewidth]{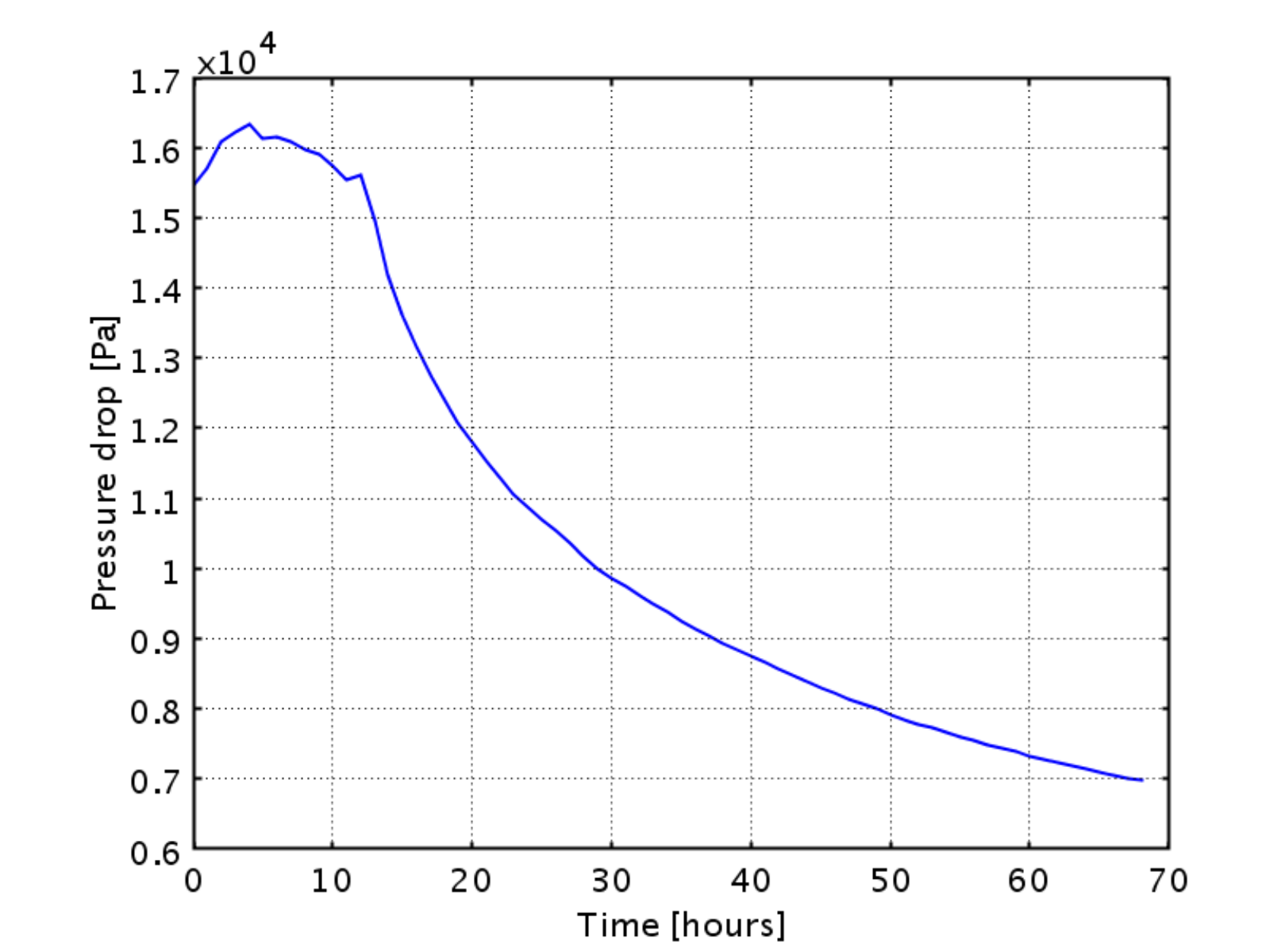}
  \captionof{figure}{Evolution of the oil pressure drop during the secondary recovery stage for a period of time of 68 hours.}\label{Figura_7_1_10}
\end{minipage}
\end{figure}

 In the Fig. \ref{Figura_7_1_8} and \ref{Figura_7_1_9} it is shown the evolution of the oil pressure and the water velocity along the z axis of the core for a period of 68 hours. One can observe that the velocity presents some numerical instability while the displacement front does not break through the production end of the core; this numerical instability can be seen as well in Fig. \ref{Figura_7_1_10}, where it is observed that the profile of the pressure drop behaves more smoothly and continuously, just after around 12 hours.


%
%
 \clearpage

\subsubsection[MEOR1 stage simulation]{MEOR1 stage simulation}

 The data used for the first microbial recovery stage MEOR1 are listed in Table \ref{Tabla_7_1_5}, while Table \ref{Tabla_7_1_4} shows the parameters of the rendered fit for the relative permeability curves.

\begin{center}
\begin{longtable}{p{0.8in} p{0.7in} p{0.8in}} \hline 
\textbf{Property} &  \textbf{Value} & \textbf{Units} \\ \hline 
$c_{w_{iny} }^{m} $ & 2.00E-02 & [kg.m$^{-3}$] \\ 
$D^{*m}_{\;\;w}$ & 1.5E-09 & [$\mathrm{m^2/s}$]\\ 
${\alpha_{L}}_{w}^{m} $ & 0.01 & [m] \\ 
${\alpha_{T}}_{w}^{m} $ & 0.01 & [m] \\ 
$\rho _{m}$ & 1600 & [$\mathrm{kg/m^3}$]\\  
$g_{m}^{\max}$& 1.85E-05 & [$\mathrm{s^{-1}}$] \\  
$K_{m/n} $ & 3.20E-01 & [$\mathrm{kg/m^3}$] \\ 
$k_{c1} $ & 2.28E-05 & [$\mathrm{s^{-1}}$]\\ 
$k_{c2} $& 1.72E-06 & [$\mathrm{s^{-1}}$]\\ 
$k_{d} $& 3.56E-05 & [$\mathrm{s^{-1}}$]\\ 
$d_{m} $ & 6.10E-06 & [$\mathrm{s^{-1}}$]\\ 
$\mu _{\mathit{surf}}^{\max}$ & 5.00E-06 & [$\mathrm{s^{-1}}$]\\ 
$K_{\mathit{surf}/n} $ & 1.00 & [$\mathrm{kg/kg}$]  \\ 
$Y_{\mathit{surf}/m} $ & 7.94 & [$\mathrm{kg/kg}$]\\ 
$c_{w}^{nC}$ & 0.00 & [$\mathrm{kg/m^3}$]\\ 
$c_{w_{inj} }^{n} $ & 4.80 & [kg.m$^{-3}$]\\ 
$D^{*n}_{\;\;w}$ & 1.5E-09 & [$\mathrm{m^2/s}$]\\ 
${\alpha_{L}}_{w}^{n} $ & 0.01 & [m]\\ 
${\alpha_{T}}_{w}^{n} $ & 0.01 & [m] \\ 
$\rho _{n}$ & 1.4E+03 & [$\mathrm{kg/m^3}$]\\  
$Y_{m/n} $ & 7.00E-02 & [$\mathrm{kg/kg}$]\\ 
$Y_{\mathit{surf}/n} $ & 1.41 & [$\mathrm{kg/kg}$]\\ 
$m_{n} $ & 1.00E-08 & [$\mathrm{kg/kg}$]\\ 
$D^{*\mathit{surf}}_{\;\;w}$ & 1.5E-09 & [$\mathrm{m^2/s}$]\\ 
${\alpha_{L}}_{w}^{\mathit{surf}} $ & 0.01 & [m]\\ 
${\alpha_{T}}_{w}^{\mathit{surf}} $& 0.01 & [m]\\ 
$\rho _{\mathit{surf}}$& 1.60E+03 & [$\mathrm{kg/m^3}$]\\ 
$t$ & 237,600 & [s] \\ \hline 
\caption{Data for the MEOR1 stage of the recovery experiment DP1. $Q$, $u_{w}^{in}$, and $p^{out}$ were kept the same as in the secondary recovery stage. Further details can be found in the section for Nomenclature and Units.}\label{Tabla_7_1_5}
\end{longtable}
\end{center}

\begin{figure}[htb]
    \centering
    \includegraphics*[width=12cm,keepaspectratio]{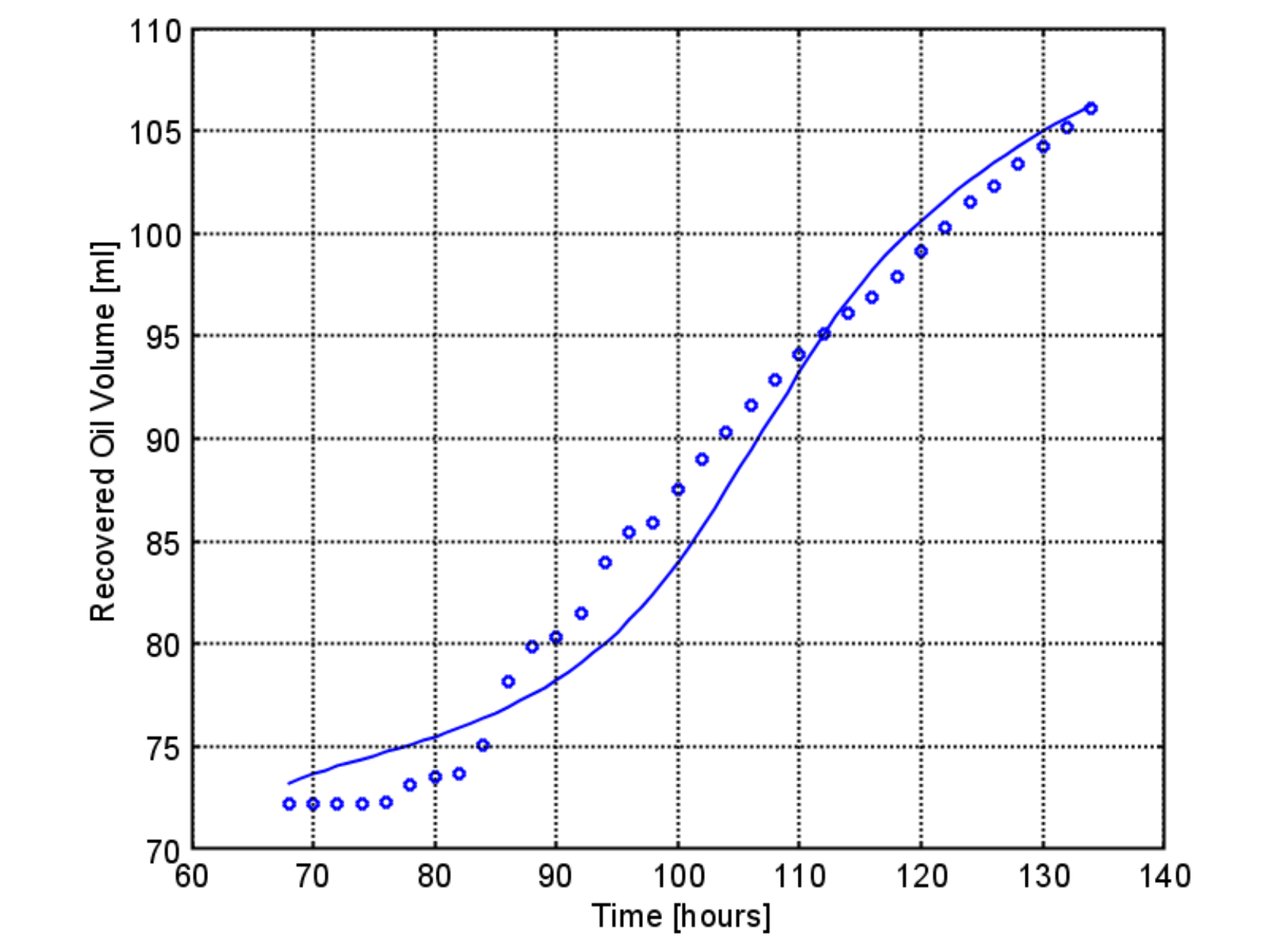}
  \caption{Comparison of the experimental (circles) and simulated (continuous line) of recovered oil volume for the MEOR1 stage.}
\label{Figura_7_1_11}           
\end{figure}

Fig. \ref{Figura_7_1_11} shows a comparison of the experimental recovery values and the simulated ones for the MEOR1 stage. One can observe an adequate fitting of the recovery curve.

 \clearpage
\begin{figure}[h!]
\centering
\captionsetup{width=0.43\textwidth}
\begin{minipage}{.5\textwidth}
  \centering
  \includegraphics[width=.9\linewidth]{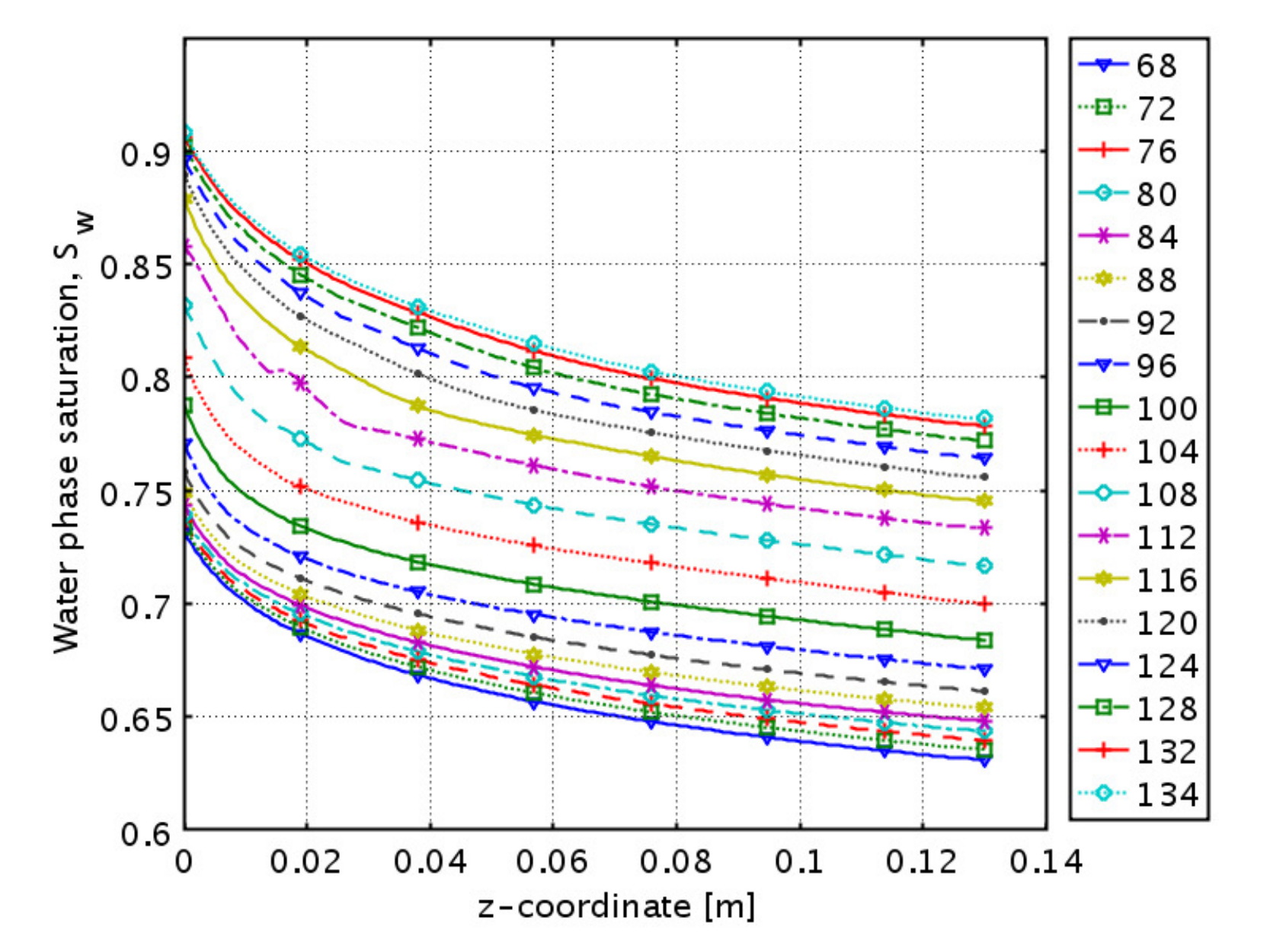}
  \captionof{figure}{Evolution of the water saturation ($S_w$) during the MEOR1 stage for a period of 66 hours. The right column is in hours.}\label{Figura_7_1_12}
\end{minipage}%
\begin{minipage}{.5\textwidth}
  \centering
  \includegraphics[width=.9\linewidth]{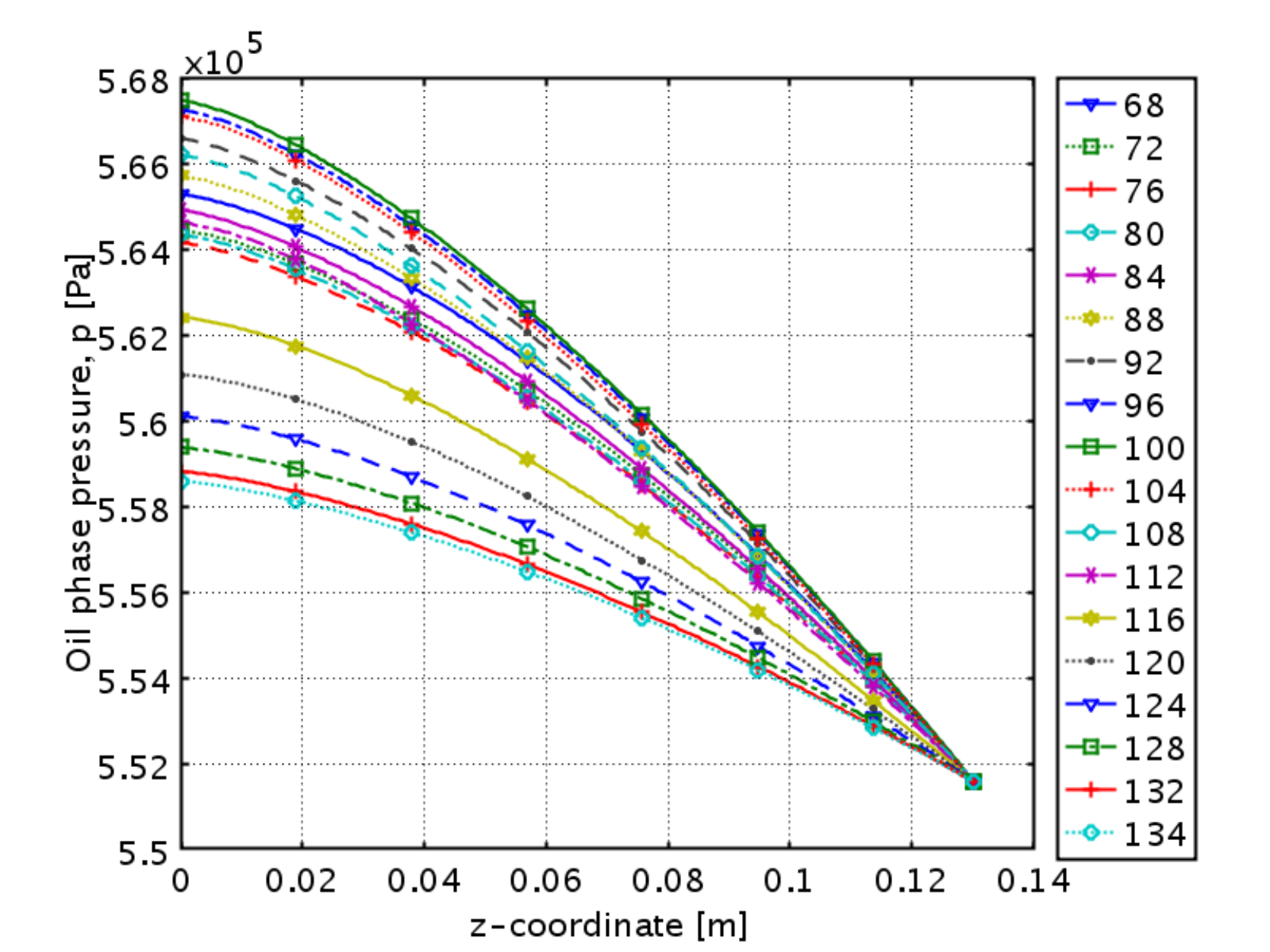}
  \captionof{figure}{Evolution of the oil pressure ($p_o$) during the MEOR1 stage for a period of 66 hours. The right column is in hours.}\label{Figura_7_1_13}
\end{minipage}
\end{figure}
\begin{figure}[h!]
\centering
\captionsetup{width=0.43\textwidth}
\begin{minipage}{.5\textwidth}
  \centering
  \includegraphics[width=.9\linewidth]{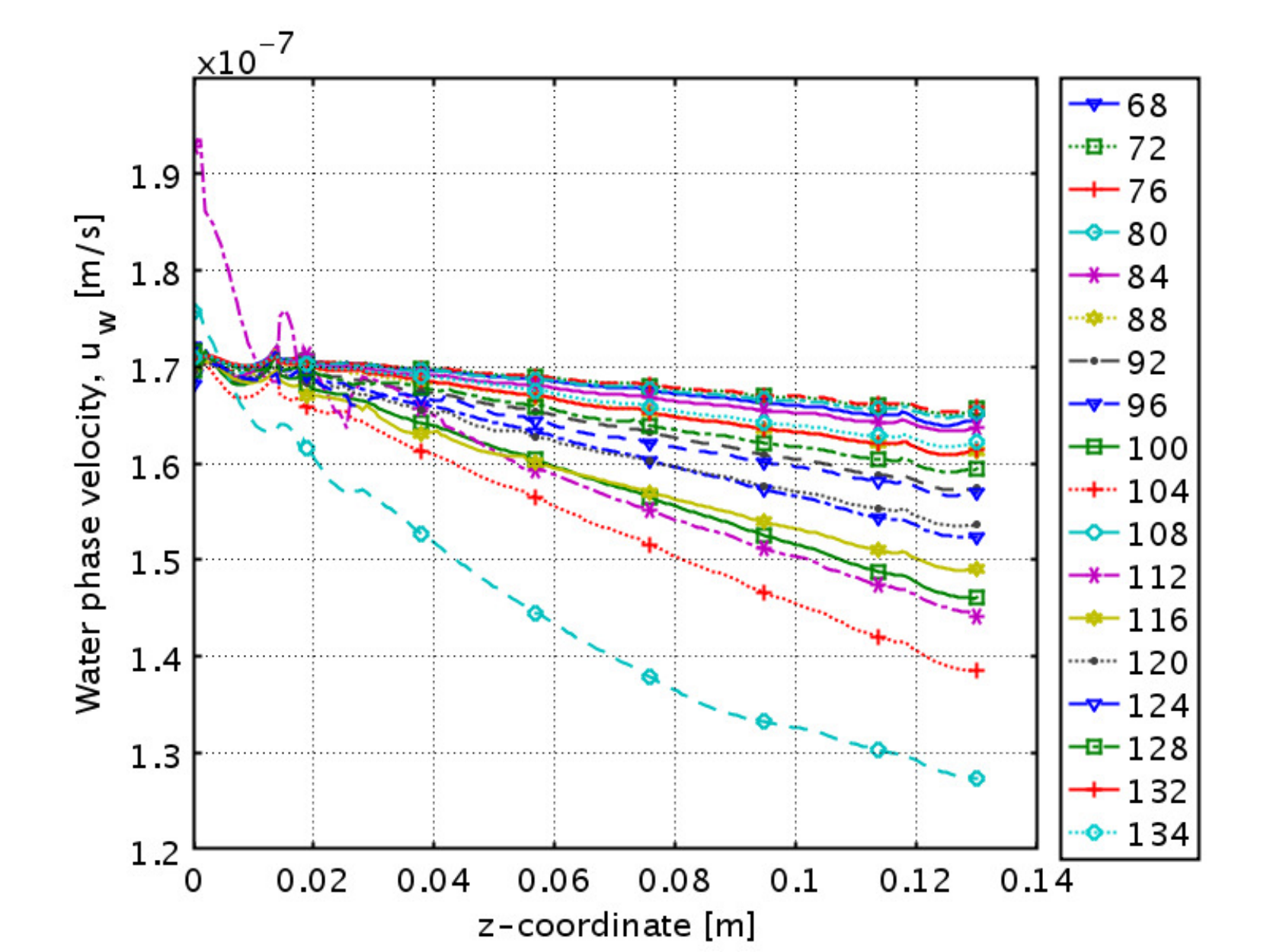}
  \captionof{figure}{Evolution of the water velocity ($u_w$) during the MEOR1 stage for a period of 66 hours. The right column is in hours.}\label{Figura_7_1_14}
\end{minipage}%
\begin{minipage}{.5\textwidth}
  \centering
  \includegraphics[width=.9\linewidth]{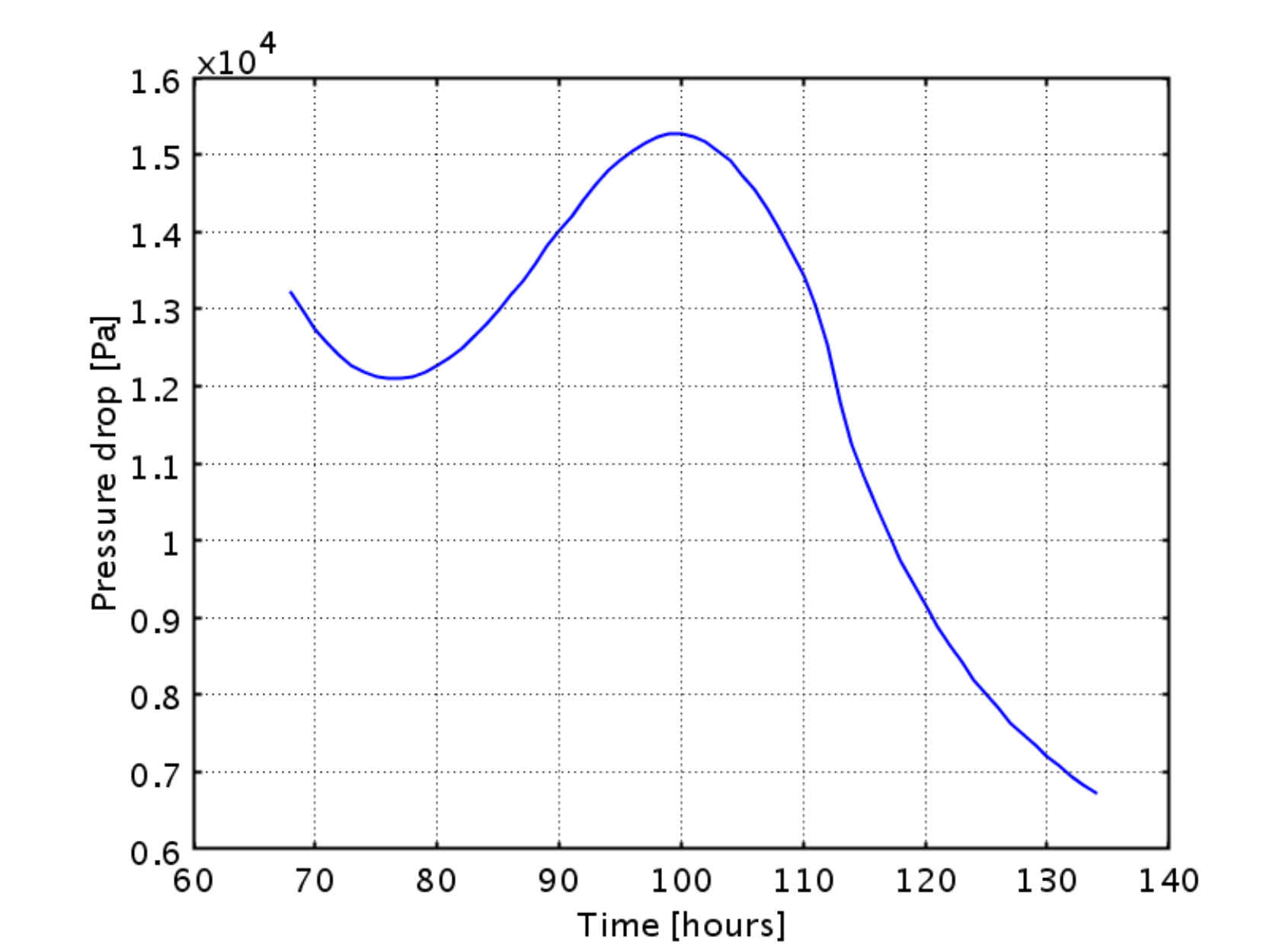}
  \captionof{figure}{Evolution of the oil pressure drop during the MEOR1 stage for a period of 66 hours.}\label{Figura_7_1_15}        
\end{minipage}
\end{figure}
%
\begin{figure}[h!]
\centering
\captionsetup{width=0.43\textwidth}
\begin{minipage}{.5\textwidth}
  \centering
  \includegraphics[width=.9\linewidth]{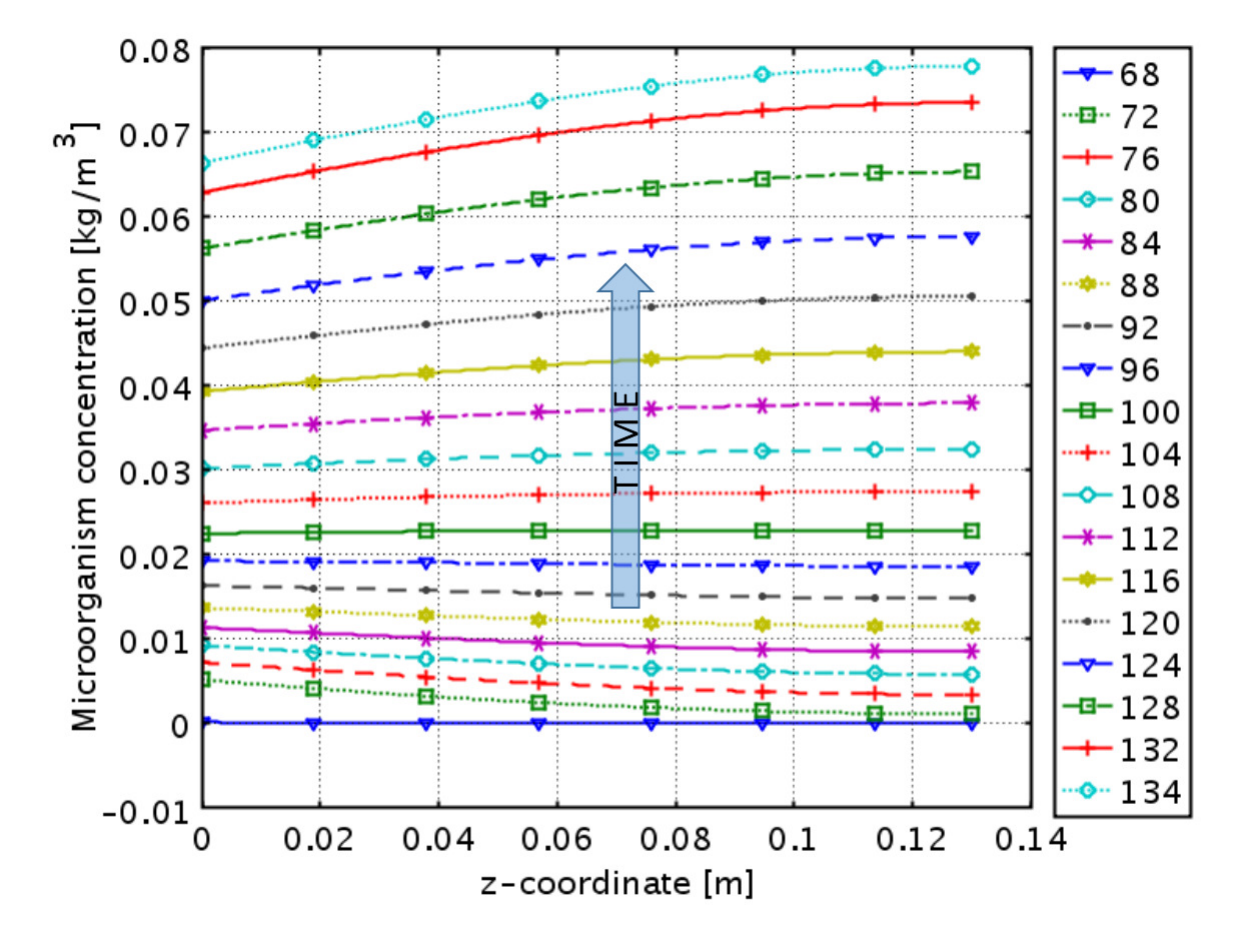}
  \captionof{figure}{Distribution of planktonic microorganisms  ($c_{w}^{m}$) along the core for different times during the MEOR1 stage. The right column is in hours.}\label{Figura_7_1_16}         
\end{minipage}%
\begin{minipage}{.5\textwidth}
  \centering
  \includegraphics[width=.9\linewidth]{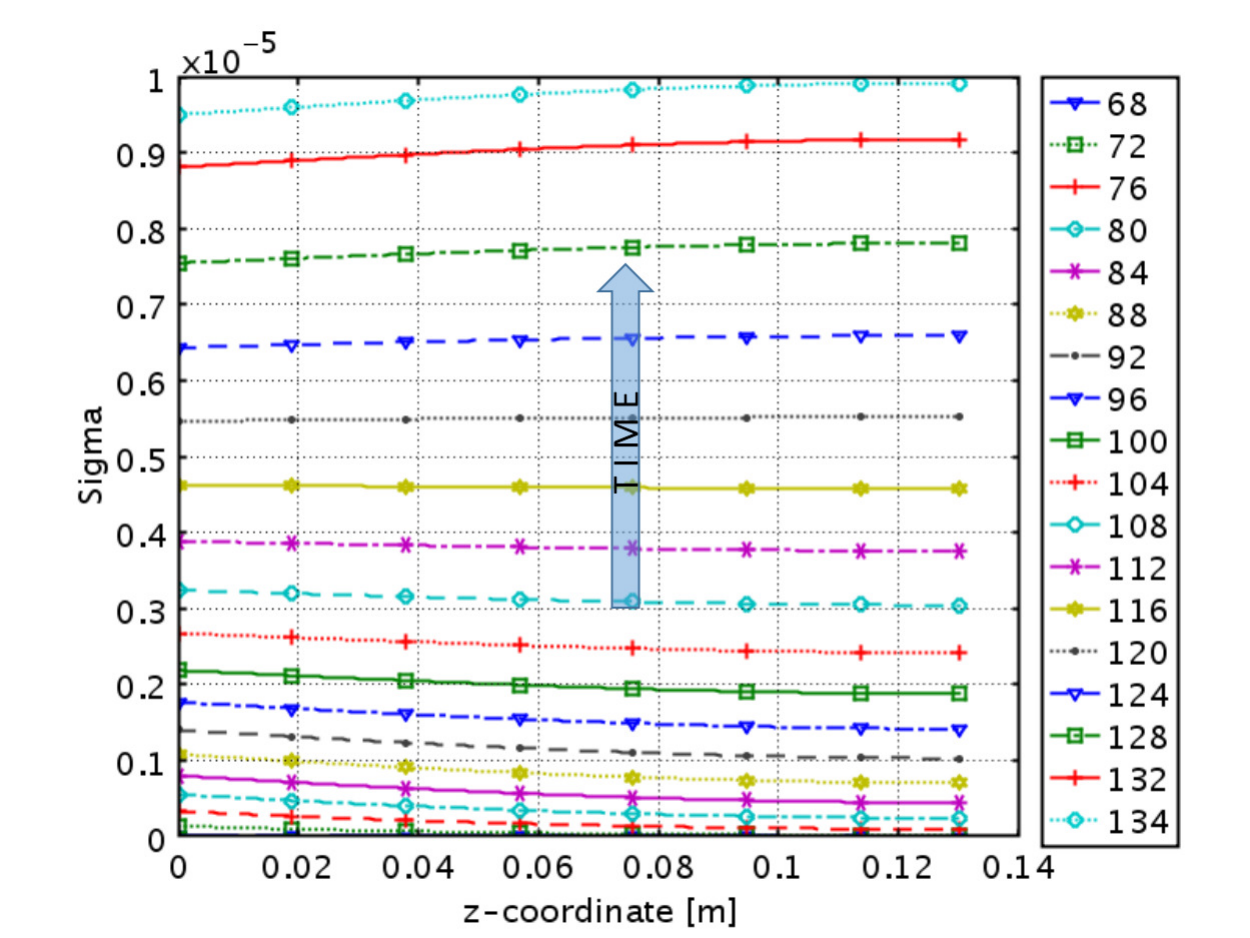}
  \captionof{figure}{Distribution of sessile microorganisms ($\sigma$) along the core for different times during the MEOR1 stage. The right column is in hours.}\label{Figura_7_1_17}          
\end{minipage}
\end{figure}
%
\begin{figure}[h!]
\centering
\captionsetup{width=0.43\textwidth}
\begin{minipage}{.5\textwidth}
  \centering
  \includegraphics[width=.9\linewidth]{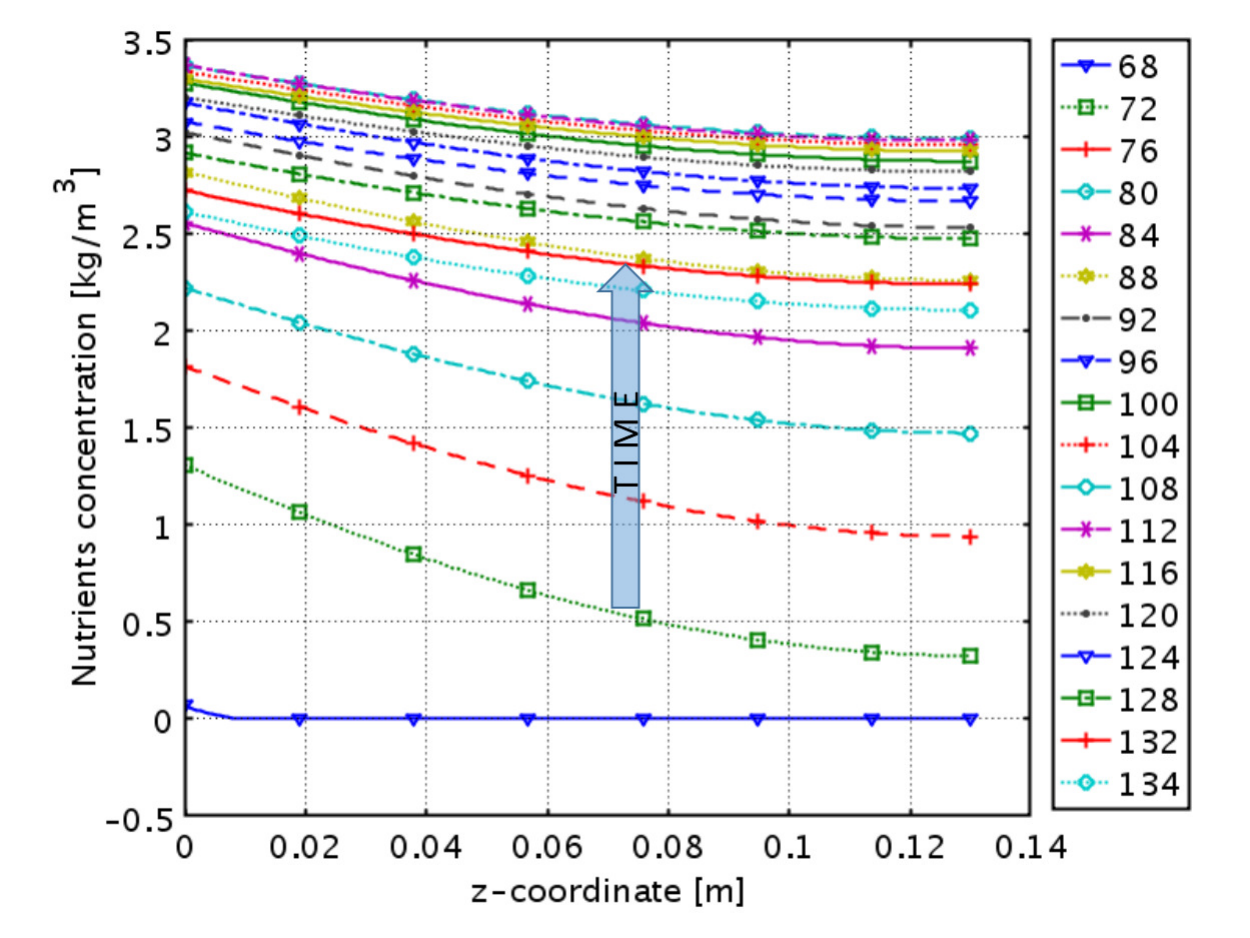}
  \captionof{figure}{Distribution of nutrients ($c_{w}^{n}$) along the core for different times during the MEOR1 stage. The right column is in hours.}\label{Figura_7_1_18}          
\end{minipage}%
\begin{minipage}{.5\textwidth}
  \centering
  \includegraphics[width=.9\linewidth]{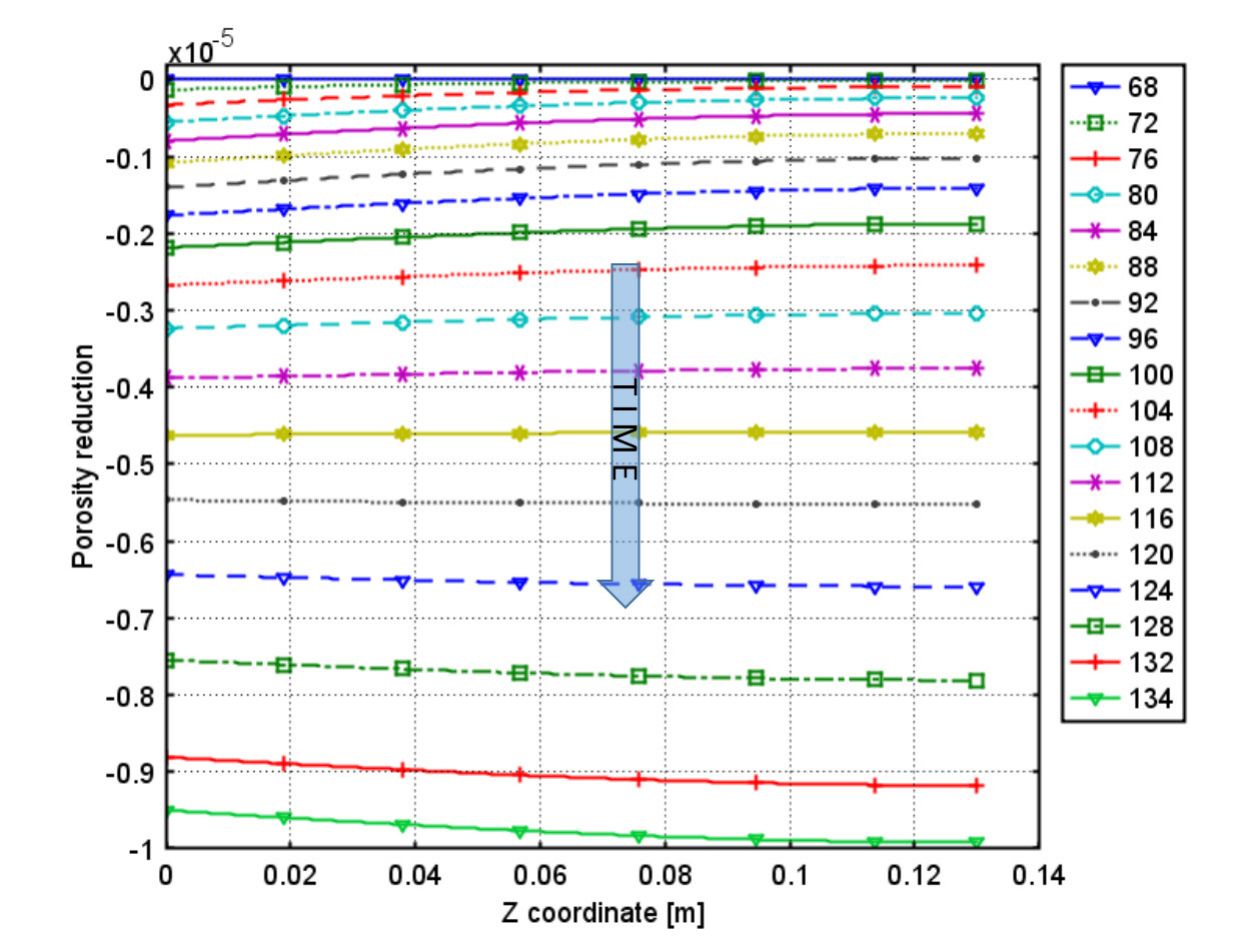}
  \captionof{figure}{Variation of porosity reduction ($\phi-\phi_0$) along the core for different times during the MEOR1 stage. The right column is in hours.}\label{Figura_7_1_20}          
\end{minipage}
\end{figure}
%

\begin{figure}[h!]
\centering
\captionsetup{width=0.43\textwidth}
\begin{minipage}{.5\textwidth}
  \centering
  \includegraphics[width=.9\linewidth]{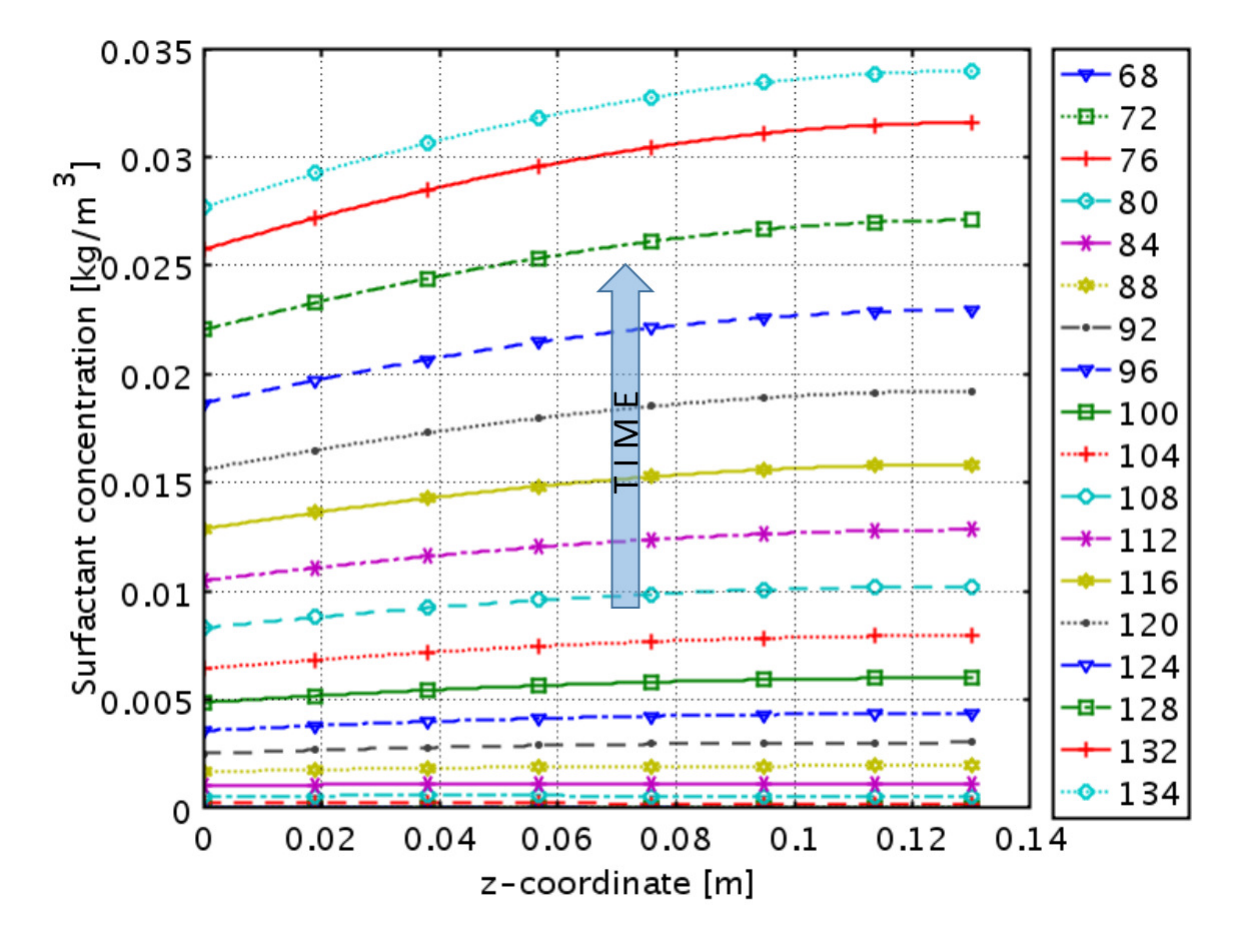}
  \captionof{figure}{Distribution of surfactant ($c_{w}^{surf}$) along the core for different times during the MEOR1 stage. The right column is in hours.}\label{Figura_7_1_19}        
\end{minipage}%
\begin{minipage}{.5\textwidth}
  \centering
  \includegraphics[width=.9\linewidth]{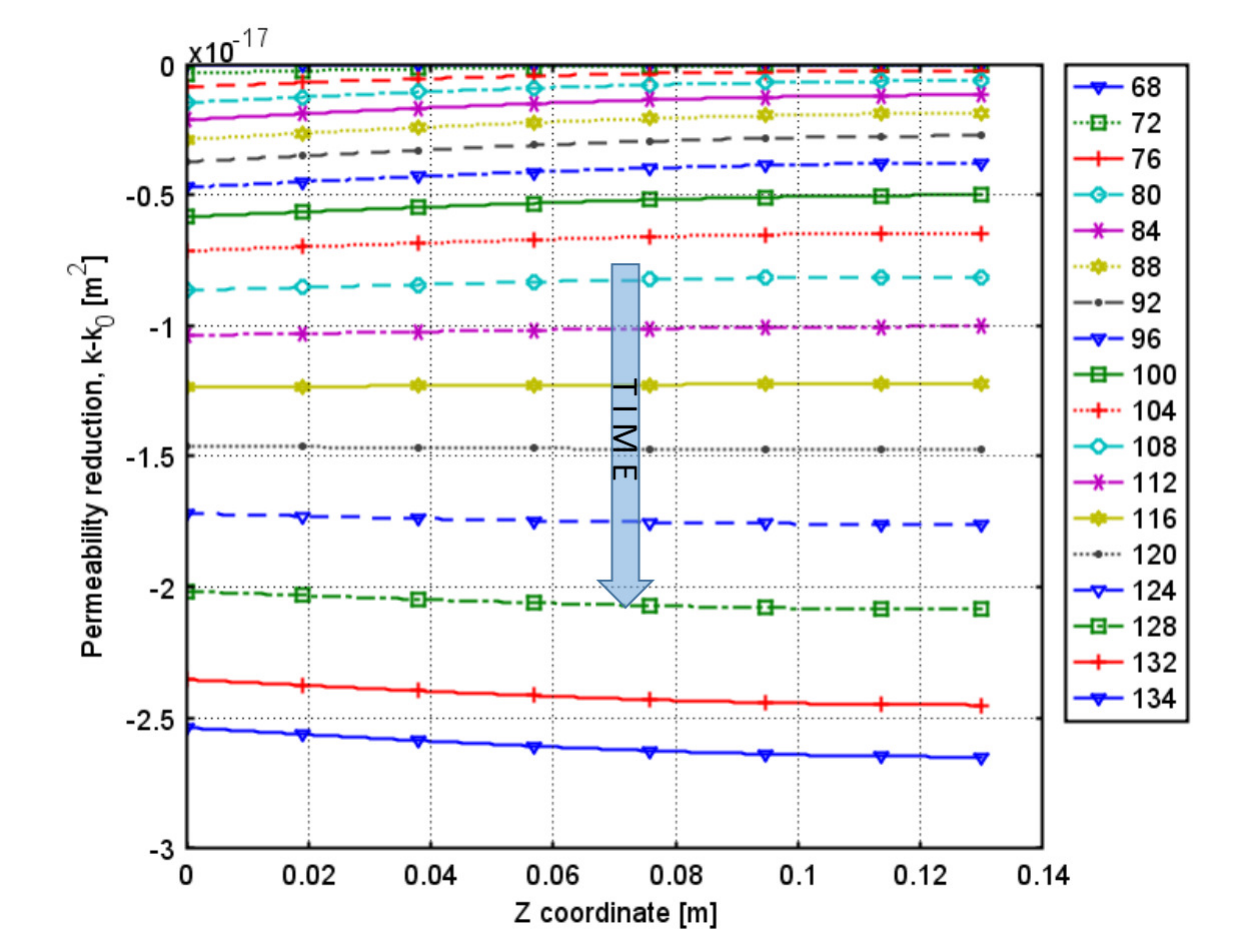}
  \captionof{figure}{Variation of absolute permeability reduction  ($k-k_0$) along the core for different times during the MEOR1 stage. The right column is in hours.}\label{Figura_7_1_21}        
\end{minipage}
\end{figure}
%
\begin{figure}[h!]
\centering
\captionsetup{width=0.43\textwidth}
\begin{minipage}{.5\textwidth}
  \centering
  \includegraphics[width=.9\linewidth]{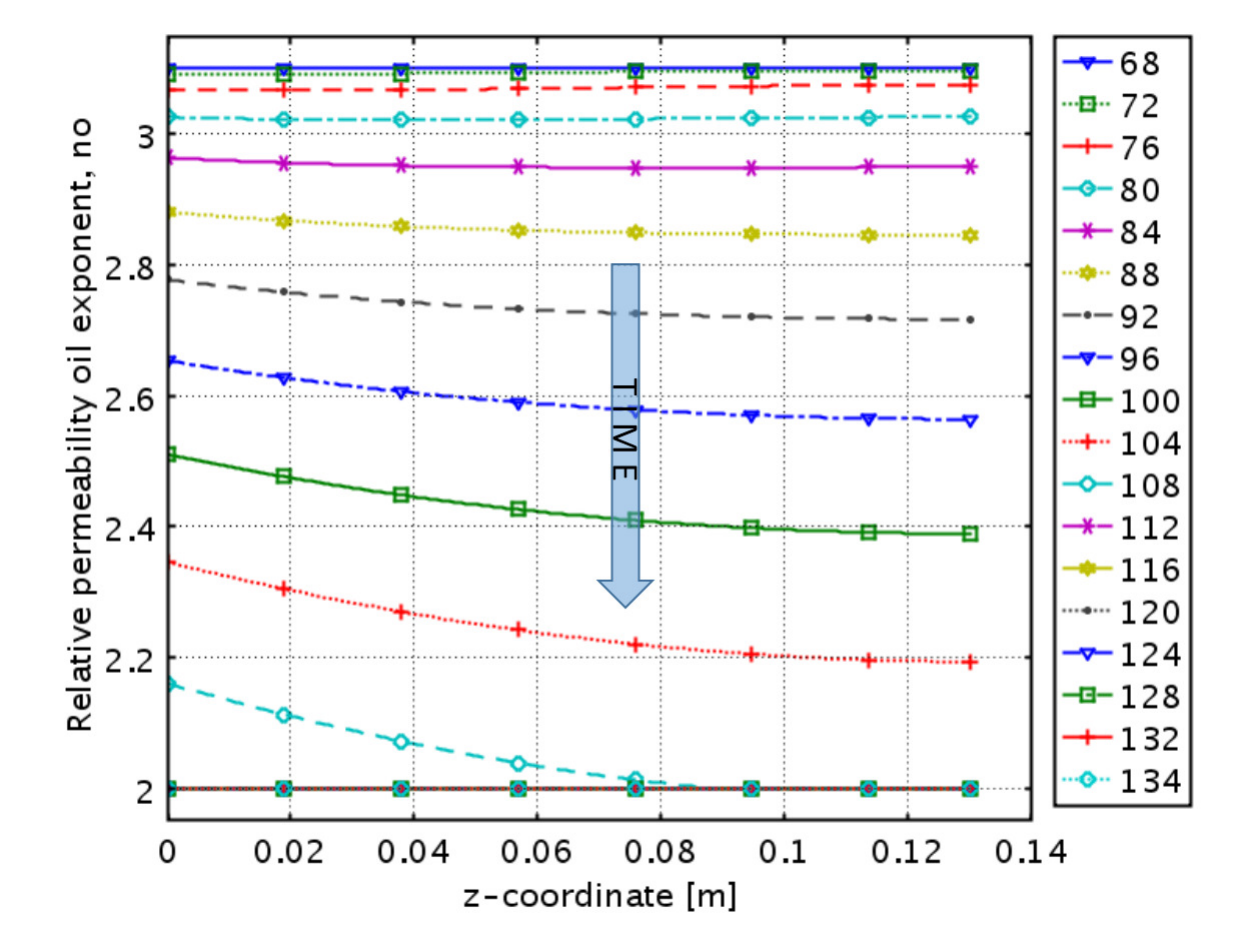}
  \captionof{figure}{Variation of the oil exponent of relative permeability ($n_o$) along the core for different times during the MEOR1 stage. The right column is in hours.}
      \label{Figura_7_1_22c}           
\end{minipage}%
\begin{minipage}{.5\textwidth}
  \centering
  \includegraphics[width=.9\linewidth]{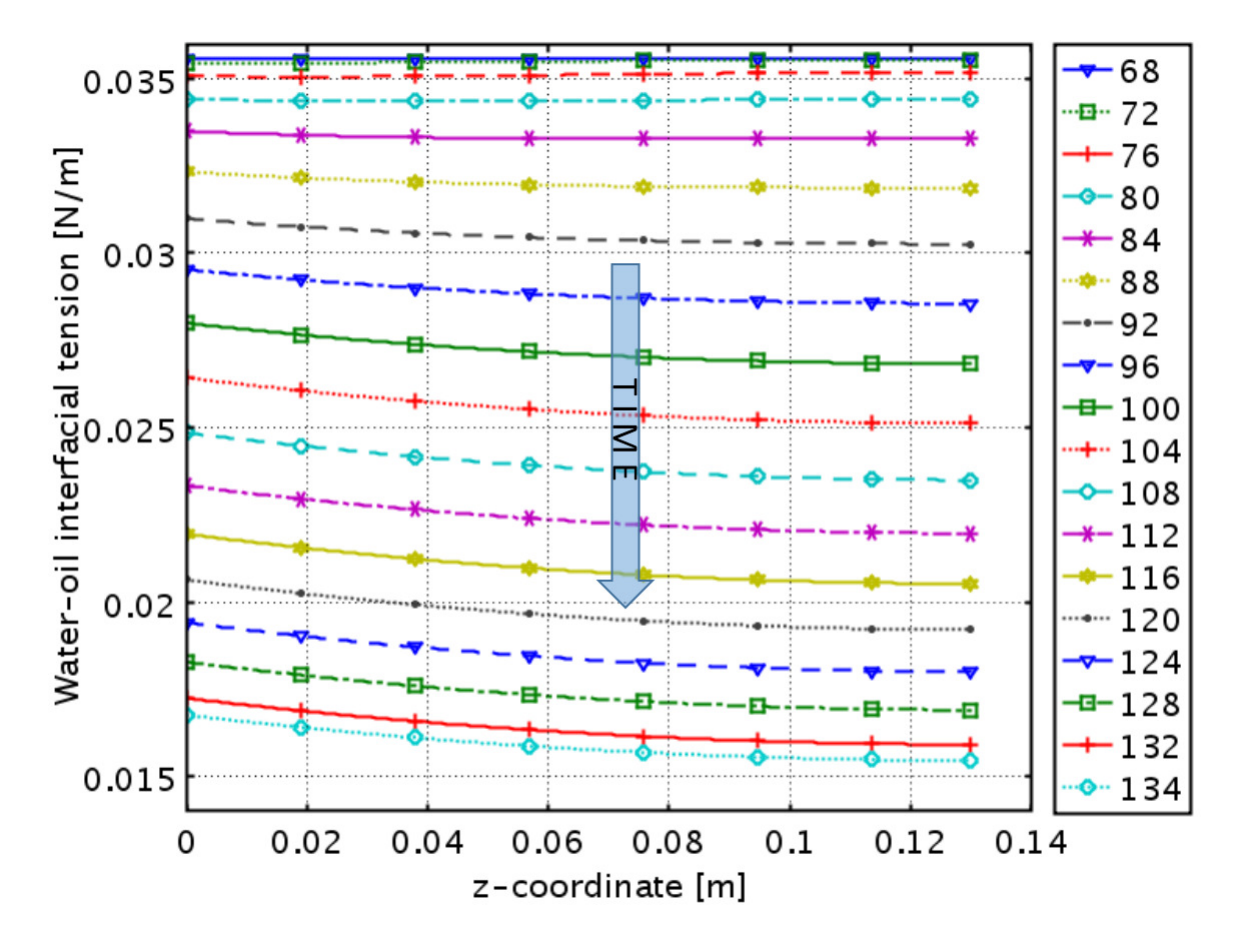}
  \captionof{figure}{Variation of oil-water interfacial tension ($\sigma_{ow}$) along the core for different times during the MEOR1 stage. The right column is in hours.}
      \label{Figura_7_1_22}           
\end{minipage}
\end{figure}
%
\begin{figure}[h!]
\centering
\captionsetup{width=0.43\textwidth}
\begin{minipage}{.5\textwidth}
  \centering
  \includegraphics[width=.9\linewidth]{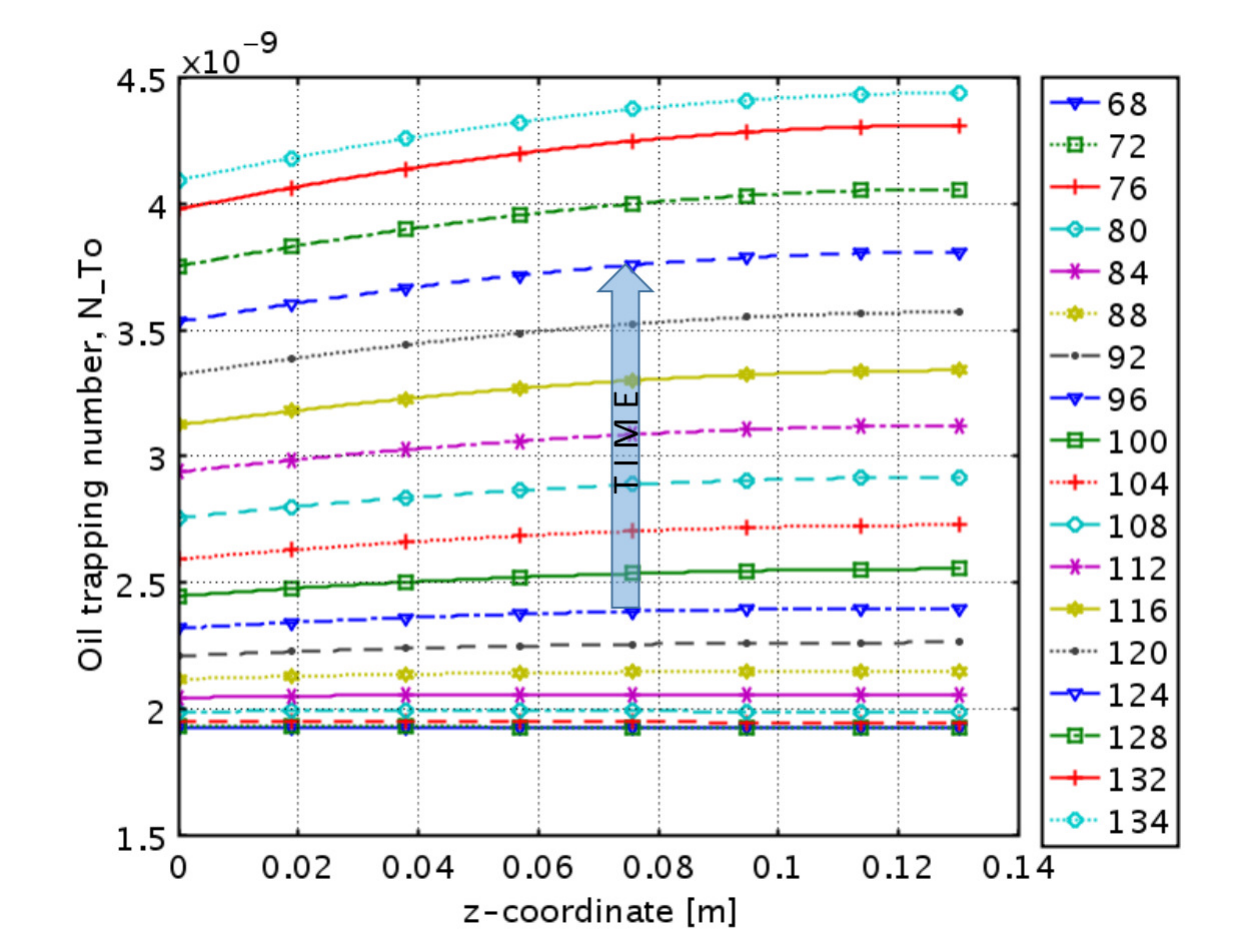}
  \captionof{figure}{Variation of the oil phase trapping number ($N_{To}$) along the core for different times during the MEOR1 stage. The right column is in hours.}
      \label{Figura_7_1_22a}           
\end{minipage}%
\begin{minipage}{.5\textwidth}
  \centering
  \includegraphics[width=.9\linewidth]{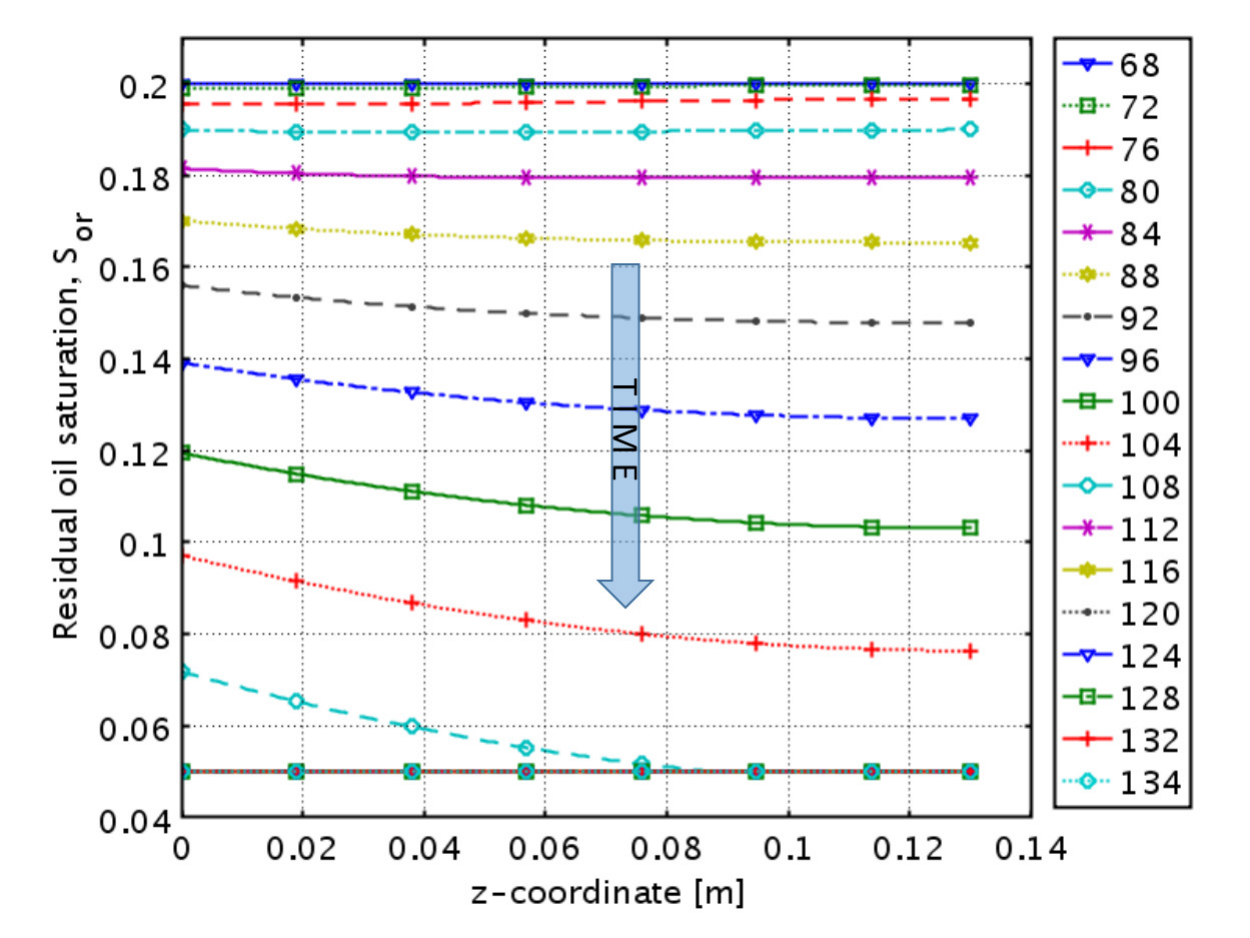}
  \captionof{figure}{Variation of the residual oil saturation ($S_{or}$) along the core for different times during the MEOR1 stage. The right column is in hours.}
      \label{Figura_7_1_22b}           
\end{minipage}
\end{figure}
%
\begin{figure}[h!]
\centering
\captionsetup{width=0.43\textwidth}
\begin{minipage}{.5\textwidth}
  \centering
  \includegraphics[width=.9\linewidth]{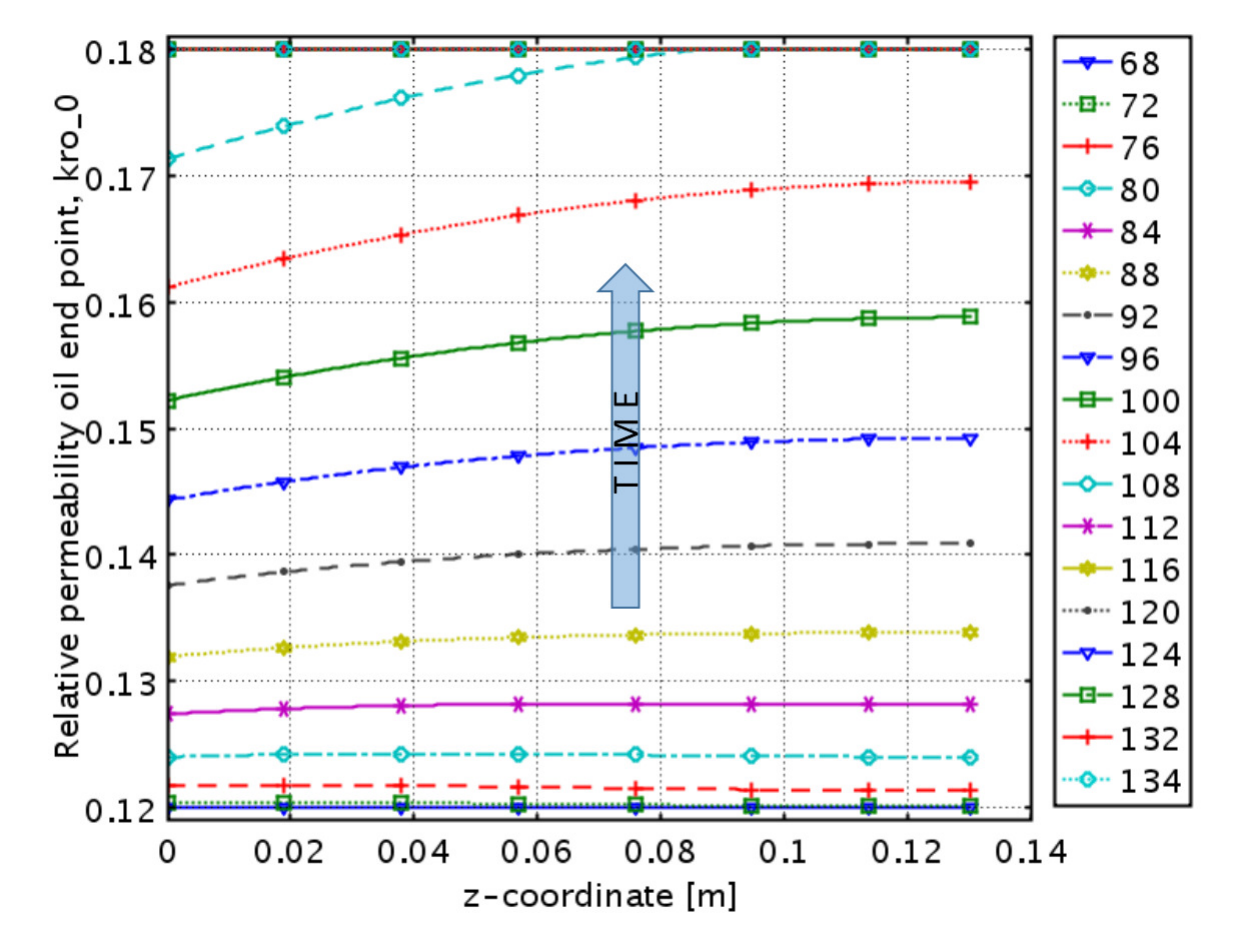}
  \captionof{figure}{Variation of the end point for the oil relative permeability ($k_{ro}^0$) along the core for different times during the MEOR1 stage. The right column is in hours.}
      \label{Figura_7_1_22d}           
\end{minipage}%
\begin{minipage}{.5\textwidth}
  \centering
  \includegraphics[width=.9\linewidth]{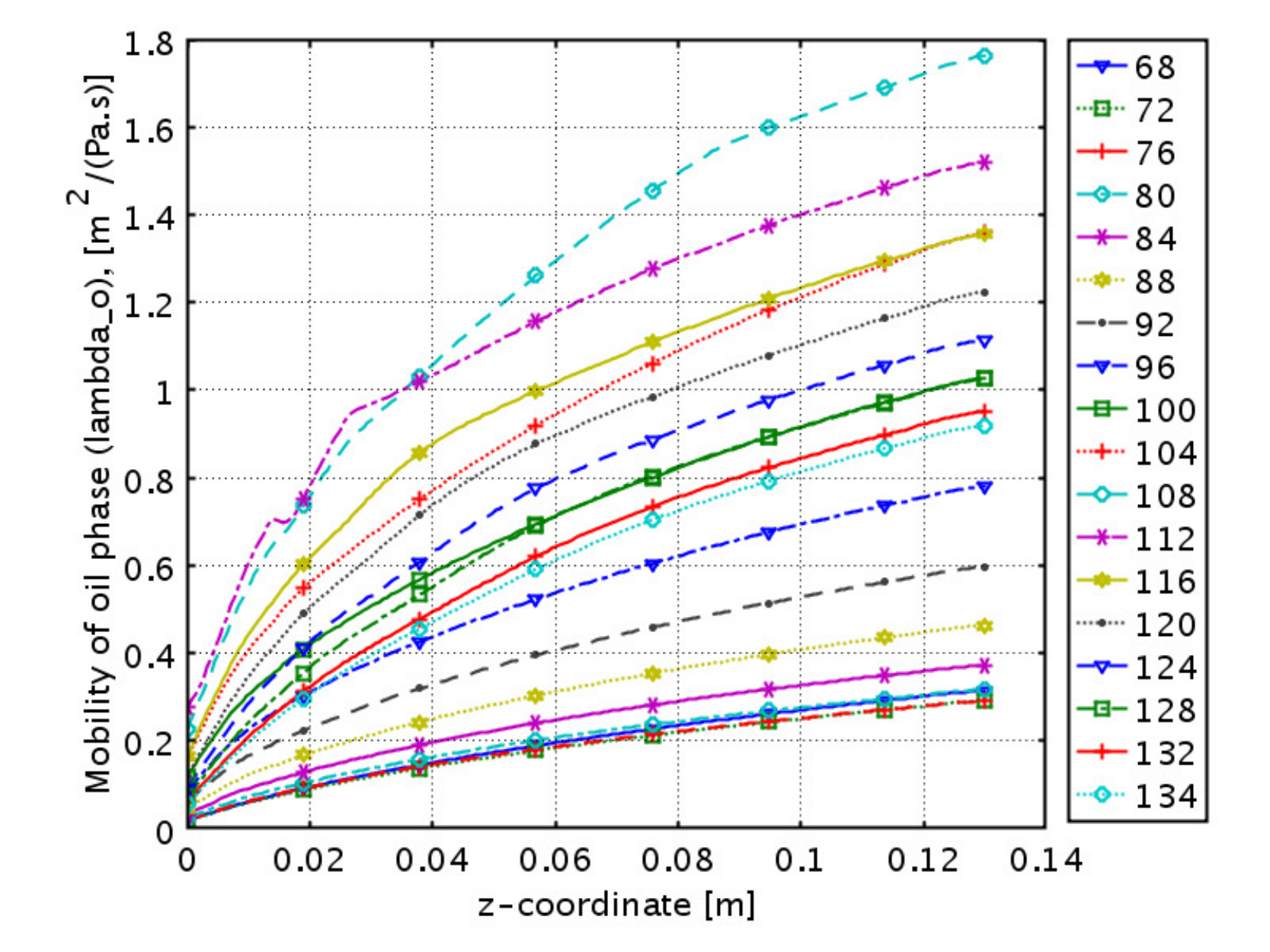}
  \captionof{figure}{Variation of the oil phase mobility ($\lambda_o$) along the core for different times during the MEOR1 stage. The right column is in hours.}
      \label{Figura_7_1_22e}           
\end{minipage}
\end{figure}
%
 
\clearpage


For MEOR1, Fig. \ref{Figura_7_1_12} shows the variation of the water saturation along the core, for different times, rising in a consistent and smooth manner, as the evolution of the oil pressure (Fig. \ref{Figura_7_1_13}). The evolution of the water velocity along the core for different times remains within a factor of about 0.75, from the maximum to the minimum, showing just a slight drop towards the end (Fig. \ref{Figura_7_1_14}).  


It is worth noticing that in Fig. \ref{Figura_7_1_15} there is an increase of the oil pressure drop from the time 80 hours up to around 100 hours, which is congruent with the behavior of Fig. \ref{Figura_7_1_17}, where there is an increase and later a drop in the concentration of sessile organisms for a similar period of time.
 
 The Fig.  \ref{Figura_7_1_16} and \ref{Figura_7_1_17} show the variation of the distribution of planktonic and sessile microorganisms, respectively, during the MEOR1 stage. The curve of Fig. \ref{Figura_7_1_16} forms a maximum displacing with the flow.

 Fig. \ref{Figura_7_1_18} shows the temporal variation of the spatial distribution of nutrients during the 66 hours. In here it is appreciated that according to the elapsed time, the nutrients concentration at the beginning (bottom end) of the core corresponds to the continuous injection of nutrients at that concentration, but decreases with time, as it is consumed by the microorganisms. This nutrient usage is not total, since at 66 hours there is a nutrient concentration distribution at opposite side.
 
 Fig. \ref{Figura_7_1_20} shows a slight, but consistent drop in porosity, behavior that is matched backwards in Fig. \ref{Figura_7_1_17} (the distribution of sessile microorganisms along the core increases its concentration as time goes by), and even Fig. \ref{Figura_7_1_14}, where it is seen that for the last period, the water velocity decreases even further than for all previous times. 
 
 
 There is an increase of the concentration of surfactant along the core as time increases (Fig. \ref{Figura_7_1_19}), which is consistent with the increase of the concentration of both planktonic and sessile microorganisms (Figs. \ref{Figura_7_1_16} and \ref{Figura_7_1_17} ), which produce it. The rise of the concentration of microorganisms also match the drop of absolute permeability (Fig. \ref{Figura_7_1_21} ). The increase of surfactant concentration is reflected also in the drop of the exponent of relative permeability, and the drop of the oil-water interfacial tension (Figs. \ref{Figura_7_1_22c} and \ref{Figura_7_1_22}). 
 
As the trapping number increases (Fig. \ref{Figura_7_1_22a}), there is a decrease in the residual oil saturation (Fig. \ref{Figura_7_1_22b}), and there is an increase of the end point for the oil relative permeability (Fig. \ref{Figura_7_1_22d}), and also an increase in the oil phase mobility (Fig. \ref{Figura_7_1_22e}) up to a maximum of $1.8\, \mathrm{[m^2/(Pa\cdot s)]}$ at 108 hour,  later coming to a middle point, at about $0.9\, \mathrm{[m^2/(Pa\cdot s)]}$, suggesting a connection with the increase of surfactant concentration along the core, for increasing times (Fig. \ref{Figura_7_1_19}).
%

%

Similar observations can be made for the other MEOR stages, for all the mentioned properties.

\section{Discussion}
\label{discussion}

 The main objective pursued with the adjustment of recovery histories is to obtain the controlling parameters of the recovery mechanism when they are implemented in a process of microbial recovery. The procedure here used is fairly reliable during the first stage of imbibition since there was more control over the data obtained in the laboratory. However, in the case of the parameters which were adjusted for the second stage of microbial recovery, they posses a greater uncertainty because many of them were taken from the published literature, so the results here obtained must be considered preliminary and an additional analysis process is required.

 However, from the experimental results it is observed in the graphic of brine injection without microorganisms (Fig. \ref{Figura_7_1_5}), a recovery of 46\%, considering that the remaining oil is residual. With the brine injection with microorganisms and nutrients it is obtained an additional recovery of 22.3 \%, which is pretty significant and in the second inoculation injection it is obtained an additional 6.2 \% of oil recovery. This would not be possible without the effect of the technique of enhanced recovery over the properties of oil and the rock.  This is shown with the necessity to change the petrophysical parameters (Table \ref{Tabla_7_1_4}) in order to reproduce the additional recovery, as it is shown with comparison curves of experimentally recovered oil volume (circles)  and the numeric ones (blue line) for the MEOR1 stage (Fig. \ref{Figura_7_1_11}), and in the Table \ref{Tabla_7_1_4}.
 
   The petrophysical parameter which influenced the most was the decrease of the residual oil saturation, allowing the additional recoveries. 

\begin{figure}[htb]
    \centering
    \includegraphics*[width=12cm,keepaspectratio]{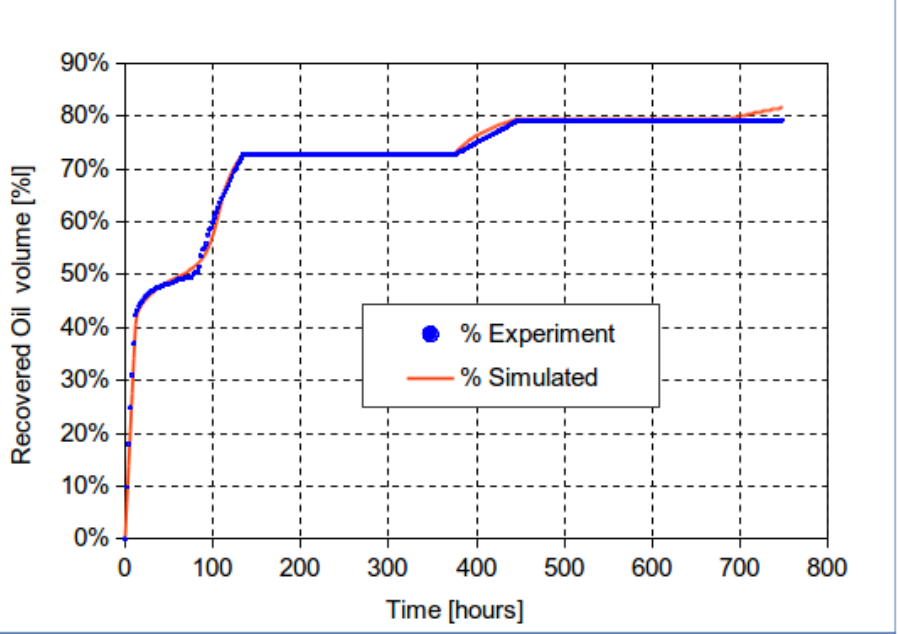}
  \caption{Comparison of the experimental (blue circles) with simulated (orange line) curves  of recovered oil during the recovery stages of the experiment DP1.}
\label{Figura_7_1_61}           
\end{figure}


In Fig. \ref{Figura_7_1_61} it is graphed the experimental (blue circles) and simulated (orange line) of the recovered oil volume during all the stages of the experiment DP1. One can confirm that the root mean square error does not go over the 2.5 ml for any of the stages, while the global root mean square error is of 2.04 ml representing the 1.76 $\%$ of the total recovered oil, see Table \ref{Tabla_7_1_2b}.

It is important to outline that the resulting error while fitting the experimental oil recovery curve to the simulated values  obtained with the model, was considering as factors impinging on the recovery,  the processes of clogging/declogging and the production of surfactant generated by the microorganisms in the microbial recovery stages, being this last the main factor. Other factors, such as, for example the produced MEOR gas (CO$_2$ y CH$_4$), were not considered in the modeling of the DP1 experiment, because the total amount of produced gases is relatively small and the production curves of those gases are not available, having only the total produced volume on each stage. The inclusion of those additional factors might allow a better fit for the recovery curve.\\

\section{Conclusions}
\label{conclusion}

A systematic methodology for a very general 3D flow and transport model in porous media for microbial enhanced oil recovery (MEOR) process under laboratory conditions comprising  conceptual, mathematical, numerical and computational stages was applied. In particular, for mathematical modeling was used the axiomatic continuum modeling approach, for numerical modeling a finite element method and COMSOL Multiphysics\circledR $\;$ Software for computational implementation.  

The transport model was validated against the experimental data from \cite{Hendry1997} and compared with the results from \cite{Kim2006} and \cite{Li2011}.

The model was successfully applied to simulate a case study which has been previously reported in the work of \cite{Castorena2012a} consisting of a laboratory experiment of microbial oil recovery process \cite{Castorena2012c}. The experimental results were accurately predicted by the simulations. 

Due to the model generality it can be easily extended and applied to other cases. The constitutive relationships such as porosity, permeability may be set up as distributions; capillary pressure, relative permeability, interfacial tension, and other relationships may be changed as well. This general model may be used for the case when the MEOR product is other than surfactant (or in addition of), and other components may be added as well.  

\section*{Acknowledgments}
 This research was supported by IMP-D.00417 Project ``Recuperaci\'on Mejorada de Hidrocarburos V\'ia Microbiana'' of the Instituto Mexicano del Petr\'oleo.

\bibliographystyle{siamplain} 
\bibliography{MEOR_ijnam}

%
%
%
%
%
\section[Nomenclature and units]{Nomenclature and units}

${\alpha_{L}}_{w}^{\eta} $, ${\alpha_{T}}_{w}^{\eta} $: Longitudinal and transversal dispersivity coefficients for the $\eta$ component($\eta = m,\,n,\,\mathit{surf}$) [m] 

$\gamma _g$: Gravitational acceleration constant [$\mathrm{m/s^2}$] 

$\theta $: Capillary pressure exponent for the Brooks-Corey model Log[Pa]

$\theta_{ow}$: Contact angle between the aqueous/non-aqueous interface and the porous medium [$\mathrm{rad}$]

$\lambda_\alpha ;\; \lambda$: Phase or total mobility, respectively $\alpha = w,o$ [$\mathrm{m^2}$/($\mathrm{Pa.s}$)]

$\mu _{\alpha} $: Phase viscosity $\alpha = w,o$ [$\mathrm{Pa.s}$]
 

$\mu _{\mathit{surf}}^{\max } $: Maximal production rate specific for the surfactant ($\mathit{surf}$) [$\mathrm{s^{-1}}$]

$\rho _{\eta} $: Density for the $\eta$ component($\eta = m,\,n,\,\mathit{surf}$)  at surface conditions [$\mathrm{kg/m^3}$]

$\rho_{\alpha}$: Phase density ($\alpha = w,\,o$) [kg/m$^3$]

$\sigma,\; \sigma _{r},\;\sigma _{i}$: Total, reversible and irreversible (respectively) volume fraction occupied by the sessile microorganisms over total pore volume [$\mathrm{m^3/m^3}$]

$\sigma_{ow}$: Oil-water interfacial tension, [${{\mathrm{mN}}/{\mathrm{m}}}$]

$\tau $: Porous medium tortuosity ($>$1) [$\mathrm{m/m}$]

$\phi,\; \phi_{0}$: Current and initial porosity of porous medium [$\mathrm{m^3/m^3}$]

$c^{\eta _{in} }_w$: Injected concentration of microorganisms or nutrients, respectively ($\eta = m,\,n$) [kg.m$^{-3}$]

$c_{w}^{\eta}$: Concentration of microorganisms, nutrients or surfactant, respectively ($\eta = m,\,n,\,\mathit{surf}$), in the water phase ($w$) [$\mathrm{kg/m^3}$]

$c_{w}^{nC} $: Critical substrate concentration ($n$)  for the surfactant formation [$\mathrm{kg/m^3}$]

$d$: Core diameter [m]

$D^{*\eta}_{\;\;w}$: Molecular diffusion coefficient for the $\eta$ component \iwp   ($w$), [$\mathrm{m^2/s}$]
 
$\underline{\underline{D}}_{\; w}^{\eta}  $: Hydrodynamic dispersion tensor for the component microorganisms, nutrients or surfactant, respectively, ($\eta = m,\,n,\,\mathit{surf}$), in the water phase ($w$), [$\mathrm{m^2/s}$]

$d_{m} $: Decaying rate (death) of microorganisms, [$\mathrm{s^{-1}}$]

$f_\alpha$: Phase fractional flow ($\lambda_\alpha / \lambda , \,\alpha = w,o$)

$g_{m} $: Microbial rate growth [$\mathrm{s^{-1}}$]
 
$g_{m}^{\max } $: Maximal growth rate for microorganisms [$\mathrm{s^{-1}}$] 

$k,\; k_{0}$: Current and initial absolute permeability [$\mathrm{m^2}$]

$k_{c1} $: Reversible clogging rate [$\mathrm{s^{-1}}$]

$k_{c2} $: Irreversible clogging rate [$\mathrm{s^{-1}}$]

$k_{d} $: Declogging rate [$\mathrm{s^{-1}}$]
 
$K_{m/n} $: Saturation constant for the growth of microorganisms ($m$) through substrate consumption ($n$) [$\mathrm{kg/m^3}$]

$k_{r\,\alpha}$: Phase relative permeability ($\alpha = o,\, w$) [$\mathrm{m^2}$]

$k_{r\,\alpha}^0$: End point for the phase relative permeability, at residual saturation of the other phase, in the modified Brooks-Corey model ($\alpha = o,\, w$) [$\mathrm{m^2/m^2}$]

$k_{ro}^{0\:low/high}$: End point for the oil relative permeability for a low/high trapping number in the modified Brooks-Corey model [$\mathrm{m^2/m^2}$]

$K_{\mathit{surf}/n} $: Saturation constant for the surfactant formation ($\mathit{surf}$) through substrate consumption ($n$) [$\mathrm{kg/kg}$]

$L$: Core Length [m]

$m_{n} $: Energy coefficient for maintaining life through substrate consumption [$\mathrm{kg/kg}$]
 
$n_{\,\alpha}$: Exponent for the phase relative permeability in the modified Brooks-Corey model ($\alpha = o,\, w$)

$N_{CA}$: Capillary number [$\mathrm{(m^2.Pa)/mN}$]

$N_B$: Bond number [$\mathrm{(m^2.Pa)/mN}$]

$N_{To}$: Oil phase trapping number [$\mathrm{(m^2.Pa)/mN}$]

$p_{cow}$: Oil-water capillary pressure [$\mathrm{Pa}$]

$p_o$: Oil pressure [$\mathrm{Pa}$]

$p_o^{0}$: Initial oil pressure [$\mathrm{Pa}$]

$p^{out}$: Production pressure [Pa]

$p_{t}$: Entry capillary pressure for the Brooks-Corey model [$\mathrm{Pa}$]

$Q$: Injection/production rate [$\mathrm{m^3/s}$]

$S_{e}$: Effective or normalized saturation [1]

$S_{\alpha}$: Phase saturation $\alpha = w,o$ [$\mathrm{m^3/m^3}$]

$S_{\alpha \, r}$:  Residual phase saturation ($\alpha = o,\, w$) [$\mathrm{m^3/m^3}$]

$S_{\eta 0}$: Initial phase saturation ($\eta = w,\,o$) [1]

$S_{or}^{low/high}$: Residual oil phase saturation for a low/high trapping number  [$\mathrm{m^3/m^3}$]

$t$: Stage duration [s] 

$\underline{u} $: Total Darcy velocity [$\mathrm{m/s}$]

$\underline{u}_\alpha $: Phase Darcy velocity $\alpha = w,o$ [$\mathrm{m/s}$]

$u_{w}^{in}$: Injection velocity [m.s${}^{-1}$]

$V$: Global block volume  [$\mathrm{m^3}$]

$V_p$: Initial pore volume [m$^3$]

$Y_{\eta/n} $: Yield coefficient of the microbial cell mass, or surfactant mass; mass of produced cells or surfactant ($\eta = m,\,\mathit{surf}$) per unit of removed substrate mass (nutrient)  [$\mathrm{kg/kg}$]

\end{document}